\documentclass[12pt]{amsart}

\usepackage{amssymb}
\usepackage[dvips]{epsfig}
\usepackage{graphicx}
\usepackage{color}
\usepackage{graphicx}
\usepackage{subcaption}
\newtheorem{theorem}{Theorem}[section]

\newtheorem{remark}[theorem]{Remark}

\begin{document}
\title[Numerical study of the 2D KBK system]
{Numerical study of the 2D Kaup-Broer-Kupershmidt Boussinesq system}

\author[T.~Gaudry]{Th\'eo Gaudry}
\address[T.~Gaudry]
{Université Bourgogne Europe, CNRS, IMB UMR 5584, 21000 Dijon, France}
\email{theo.gaudry@ube.fr}

\author[C. Klein]{Christian Klein}
\address[C.~Klein]
{Université Bourgogne Europe, CNRS, IMB UMR 5584, 21000 Dijon, 
France\\
Institut Universitaire de France}
\email{christian.klein@u-bourgogne.fr}

\author[J.-C. Saut]{Jean-Claude Saut}
\address[J.-C. Saut]{Laboratoire de Math\' ematiques, UMR 8628,\\
Universit\' e Paris-Saclay et CNRS\\ 91405 Orsay, France}
\email{jean-claude.saut@universite-paris-saclay.fr}

\author[N.~Stoilov]{Nikola Stoilov}
\address[N.~Stoilov]
{Université Bourgogne Europe, CNRS, IMB UMR 5584, 21000 Dijon, France}
\email{nikola.stoilov@ube.fr}

\begin{abstract}
In this work we consider the well posed version of the 
Kaup-Broer-Kupershmidt system in two dimensions 
irrotational case. We numerically construct soliton type solutions and show that they are unstable both against dispersion and singularity formation. Further, we study line solitons and their stability, as well as generally localised initial data.   
\end{abstract}

\date{\today}


\thanks{This work was supported by the ANR project 
ANR-17-EURE-0002 EIPHI and by the ANR project 
ISAAC-ANR-23-CE40-0015-01. We thank the anonymous referees for 
helpful suggestions.  }
\maketitle

\section{Introduction}
This note is concerned with a particular case of the Boussinesq 
`abcd' systems derived in \cite{BCS1, BCS2}, see also \cite{CFOPS}, namely the two-dimensional 
Kaup-Broer-Kupershmidt (KBK)  in the 
irrotational case. The KBK system has two versions, 
the  linearly ill-posed `bad' one, and the `good' one, that is well posed. We concentrate on the latter:
\begin{equation}
    \label{Kaup-true-2D}
    \left\lbrace
    \begin{array}{l}
    \eta_t+\nabla \cdot {\bf v}+\nabla\cdot(\eta {\bf v})-\Delta \nabla\cdot {\bf v}=0,\\
    {\bf v}_t+\nabla \eta+ \frac{1}{2}\nabla |{\bf v}|^2=0,
\end{array}\right.
    \end{equation}
    
    The terminology "KBK system" arises from the one-dimensional 
	version of \eqref{Kaup-true-2D} introduced and studied  in \cite{Broer, BGT, Kaup2, Ka, Kup} as a completely integrable system; we refer to the survey article \cite{KS-KBK} for more details.
    
    In \cite{BCS1} the KBK system in both  spatial dimensions is 
	labeled  as a `Schr\"{o}dinger type Boussinesq system' since the eigenvalues of the dispersion matrix behave as $\pm|\xi|^2$ as $|\xi|\to +\infty.$
	
	The KBK system in one-dimension  has a bi-Hamiltonian structure, see \cite{Kup}. In  two dimensions,  with the standard Poisson bracket, the Hamiltonian is 
 
\begin{equation*}\label{energy}
	H = 
	\int_{\mathbf{R}^{2}}^{}\left[\eta^{2}+(1+\eta)|\mathbf{v}|^{2}+|\nabla 
	\mathbf{v}|^{2}\right].
\end{equation*}

Conservation of energy, however, is not sufficient to ensure the global existence of solutions of the Cauchy problem which remains an open problem.

Another important issue, necessary for the rigorous justification of the KBK system as an asymptotic model on relevant time scales, is the long time existence for the system written as 


\begin{equation*}
    \label{Kaup-true-2D-bis}
    \left\lbrace
    \begin{array}{l}
    \eta_t+\nabla \cdot {\bf v}+\epsilon \nabla\cdot(\eta {\bf v})-\epsilon \Delta \nabla\cdot {\bf v}=0,\\
    {\bf v}_t+\nabla \eta+ \frac{\epsilon}{2}\nabla |{\bf v}|^2=0,
\end{array}\right.
    \end{equation*}
    where the small parameter $\epsilon$ measures the comparable nonlinear and dispersive effects.
    
    The long time  existence on time scales of order $O(1/\epsilon)$ for the KBK system has been established in \cite{SWX} and recently improved in \cite{Tes-Sel}  to $O(1/\epsilon^{3/2})$ in the two-dimensional case.

    Although being completely integrable in 1D, (see the references 
	in \cite{KS-KBK}), a description of  the qualitative properties 
	of solutions of the 1D KBK system is  far from being complete, 
	for instance the possible decomposition of compact initial data 
	into solitons and radiation over sufficiently long times (that 
	is, a version of the soliton resolution conjecture) is unknown.
    
     The generalisation of this equation to two dimensions is not known to be 
	 integrable. This article numerically studies the properties of 
	 solutions to this two-dimensional version of the KBK system, and 
	 in particular looks for stable structures. We examine 2D 
	 solitons and their stability as well as line solitons, that is 
	 solutions of the 1D equation, trivially extended in the 
	 transversal direction, and their stability. Finally we consider compact data in general position and study their evolution. 
	 
	 The paper is organised as follows: In section 2, we collect some 
	 known facts on the KBK system in 1D. In section 3, we write the 
	 KBK system in 2D in  the form of a system of fractional NLS 
	 equation, discuss the possibility of a self-similar blow-up and 
	 of static solutions localised in 2D. In section 4, we 
	 numerically construct such a static solution. The numerical time 
	 evolution of 2D KBK solutions is discussed in section 5. 
	 Perturbations of the static solutions are numerically studied in 
	 section 6. In section 7 we numerically study the transverse 
	 stability of line solitons. In section 8 we study Gaussian 
	 initial data. We add some concluding remarks in section 9.

\section{KBK in 1D}
In this section we briefly collect some known facts on the 1D KBK system, more details can be found in the survey article \cite{KS-KBK}. The KBK system possesses solitary wave solutions that can be written  for the velocity 
$C\in\mathbb{R}$, $|C|<1$  in the form
\begin{align}
	v& =\frac{2(1-C^2)}{\cosh(\sqrt{1-C^2}(x-Ct-x_{0}))-C},
	\nonumber\\
	\eta & = Cv-\frac{1}{2}v^2 
	\label{soliton}.
\end{align}
Note that $\eta$ in (\ref{soliton}) can violate the non-cavitation 
condition $\eta+1>0$ for $C>-0.5$ which is important for the regularity of 
solutions of other Boussinesq systems, 
but does not affect the regularity of solutions of 1D KBK. 
Whereas in applications $\eta+1$ should never be 
negative, we will in this paper also consider cases where this 
happens.

Actually, $v$ is the unique solution (up to translations) of the ODE:

\begin{equation*}
v''+\frac{1}{2}v^3-\frac{3}{2}Cv^2-(1-C^2)v=0.
\end{equation*}

Angulo in \cite{Ang} used the following  higher order conserved quantity, 
\begin{equation*}
\begin{split}
I_3(\eta,v)=&\frac{1}{8}\int [4(v_{xx}^2)+8(v_x^2)+4v^2+4(\eta_x^2)+4\eta^2+6v^2(v_x)^2 \\
&-16\eta vv_{xx}-4\eta(v_x)^2+10\eta v^2+2\eta^3+v^4+6\eta^2v^2+\eta v^4]dx,
\end{split}\label{Ang}
\end{equation*}

\noindent consequence of the complete integrability of the KBK system, to prove that the Cauchy problem is globally well-posed in $H^s(\mathbb{R})\times H^{s-1}(\mathbb{R}), s \geq 2$ for initial data close to a translation of the solitary wave. He furthermore proved the orbital stability of the solitary wave. The asymptotic stability of the solitary wave for velocities $c\in (-1,-c^*)$ was recently  proven in \cite{ZGX, Mau-Mun}.

\begin{remark}
The explicit form of multi-soliton solutions and the proof of their stability appear to be open problems as well as the global existence of solutions corresponding to arbitrary initial data.
\end{remark}




\section{KBK in two dimensions}
The theory of the two-dimensional KBK system is much less understood, 
and no global existence of solutions is known, even for small initial data and the possible finite time blow-up is not excluded. Note, however, that linear dispersive estimates such as Strichartz and local Kato smoothing estimates are available in dimensions one and two as a consequence of the general study in \cite{Mel}.

In this section we outline how the system (\ref{Kaup-true-2D}) can be 
written as a system of fractional 2D nonlinear Schr\"odinger (NLS) 
equations, and how static solutions to said system can be obtained. 
Further, we discuss the possibility of a self-similar blow-up of solutions to the 2D KBK system. 

We use the 
following convention for the Fourier transform for sufficiently regular and localized functions $g$:
\begin{equation*}
\widehat{g}(k_{x},k_{y})  :=\int_{\mathbb{R}^{2}}    
  e^{-i(k_{x}x+k_{y}y)}\,g(x,y) \, dxdy
\end{equation*}
\begin{equation*}
    g(x,y)  =\frac{1}{(2\pi)^{2}}\int_{\mathbb{R}^{2}}
    e^{i(k_{x}x+k_{y}y)}\,\widehat{g}(k_{x},k_{y}) \,dk_{x}dk_{y};
\end{equation*}
here  $(x,y)\in\mathbb{R}^{2}$ are the variables in the physical 
space, $(k_{x},k_{y})\in\mathbb{R}^{2}$ the variables in Fourier space. 

\subsection{Fractional NLS form}

In the irrotational case we consider in this paper, it is possible to introduce a velocity potential $V$ via 
$\mathbf{v}=\nabla V$ in the system (\ref{Kaup-true-2D}). This leads  to
\begin{equation}
    \label{Kaup-true-2DV}
    \left\lbrace
    \begin{array}{l}
    \eta_t+\Delta V+\nabla\cdot(\eta \nabla V)-\Delta^{2} V=0,\\
    V_t+ \eta+ \frac{1}{2} |\nabla V|^2=0,
\end{array}\right.
   \end{equation}
where some vanishing condition at infinity for $V$ has been applied in 
integrating the second equation for localized data. Line solitons are solutions to the 2D KBK system, that do not depend on $y$, and are soliton solutions of the 1D KBK equation in $x$. For these solutions the potential $V$ does not vanish as $y$ tends to infinity.  

In Fourier space, the system (\ref{Kaup-true-2DV}) reads
\begin{equation}
    \label{Kaup-true-2DF}
    \left\lbrace
    \begin{array}{l}
    \hat{\eta}_t -|\mathbf{k}|^{2}(1+|\mathbf{k}|^{2}) \hat{V}+i\mathbf{k}\cdot\widehat{(\eta \nabla V)}=0,\\
    \hat{V}_t+ \hat{\eta}+ \frac{1}{2} \widehat{|\nabla V|^2}=0.
\end{array}\right.
   \end{equation}
   
To handle the potential difficulties at infinity
we introduce the fractional Laplacian of $V$: $W = i\sqrt{-\Delta}V$ and its Fourier image $\hat{W}= i|\mathbf{k}|\hat{V}$. This allows to write the system 
(\ref{Kaup-true-2DF})  for the quantities $\hat{u}_{\pm}$ 
given by
\begin{equation*}
	\hat{u}_{\pm}:=\hat{W}\pm\frac{\hat{\eta}}{\sqrt{1+|\mathbf{k}|^{2}}}
	\label{upm2}.
\end{equation*}
in the form
\begin{equation}
	(\hat{u}_{\pm})_{t}=\mp i|\mathbf{k}|\sqrt{1+|\mathbf{k}|^{2}}\hat{u}_{\pm}-\frac{i}{2}|\mathbf{k}|\widehat{|\nabla V|^{2}}
	\mp\frac{i\mathbf{k}\widehat{\eta\nabla 
	V}}{\sqrt{1+|\mathbf{k}|^{2}}}.
	\label{upmt2}
\end{equation}
This can be seen as a system of fractional NLS equations where the 
dispersion is proportional to $|\mathbf{k}|^{2}$ for large 
$|\mathbf{k|}$ (the dispersion relation for the linear part being
$\omega=\pm|\mathbf{k}|\sqrt{1+|\mathbf{k}|^{2}}$). The system (\ref{upmt2}) appears to be ideally suited for the use of \emph{exponential 
integrators} in the time integration, see below. 

\subsection{Static solutions to the 2D KBK system}
Static, i.e., time independent solutions to the system 
(\ref{Kaup-true-2DV}) solve the equations 
\begin{align}
	\eta+ \frac{1}{2} |\nabla 
V|^2&=0,\nonumber\\
	\nabla\cdot\left((- \frac{1}{2} |\nabla 
V|^2\nabla V)+\nabla(V-\Delta V)\right) &= 0
	\label{Vstatic}.
\end{align}
We do not explore the general solution of these 
equations, but concentrate on the localised solutions satisfying
\begin{equation}
	(- \frac{1}{2} |\nabla 
	V|^2\nabla V)+\nabla(V-\Delta V) = 0.
	\label{Vst}
\end{equation}
This corresponds to 
ground states of a vector NLS equation discussed in \cite{CW}. There 
are radially symmetric solutions to the equation (\ref{Vst}). In standard polar coordinates $x=r\cos\phi$, $y= 
r\sin\phi$, $r>0$ and $0\leq \phi<2\pi$, one gets from 
(\ref{Vst}) the equation
\begin{equation}
	V_{r}-\left(V_{rr}+\frac{1}{r}V_{r}\right)_{r}-\frac{1}{2}V_{r}^{3}=0
	\label{Vr}.
\end{equation}

\subsection{Dynamical rescaling}
An interesting question in the context of KBK is whether there can be 
a blow-up, and if yes, whether the type of the blow-up can be 
determined. The NLS form of the 2D KBK system (\ref{upmt2}) suggests 
that a  self-similar blow-up is possible. In this case, a possible 
approach to treat the blow-up  is given by a 
\emph{dynamical rescaling}, see for instance \cite{MPSS,Sul}: 
\begin{equation}
	X = \frac{x-x_{s}}{L},\quad Y = \frac{y-y_{s}}{L},\quad \eta = 
	U/L^{2}, 
	\label{dyn}\quad \frac{d\tau}{dt}=\frac{1}{L^{2}}
\end{equation}
where $L=L(\tau)$ vanishes for $\tau\to\infty$ and  $x_{s}$, $y_{s}$ give the  location 
of the blow-up. Note that the function $V$ is not changed by this 
rescaling, and that the 
$L^{2}$ norm of $\eta$ is not invariant. In this sense the equation 
is $L^{2}$ supercritical, and as shown below, the behavior  of the 
solutions is reminiscent of solutions to supercritical NLS equations. The transformation (\ref{dyn}) leads to the dynamically rescaled form of the 
system (\ref{Kaup-true-2DV})
\begin{equation}
    \label{KBKresc}
    \left\lbrace
    \begin{array}{l}
    U_\tau-a(2U+XU_{X}+YU_{Y})+L^{2}\Delta V+\nabla\cdot(U \nabla V)-\Delta^{2} V=0,\\
    V_\tau-a(XV_{X}+YV_{Y})+ U+ \frac{1}{2} |\nabla V|^2=0,
\end{array}\right.
   \end{equation}
where 
\begin{equation*}
	a=L_{\tau}/L
	\label{a}.
\end{equation*}
In the case of a self-similar blow-up, we expect the 
functions $U$, $V$ to become constant as well as $a$. Denoting these 
constants in $\tau$ with a superscript $\infty$,  we get in the limit $\tau\to\infty$ for (\ref{KBKresc})
\begin{equation*}
    \label{KBKrinf}
    \left\lbrace
    \begin{array}{l}
    -a^{\infty}(2U^{\infty}+XU^{\infty}_{X}+YU^{\infty}_{Y})+\nabla\cdot(U^{\infty} \nabla V^{\infty})-\Delta^{2} V^{\infty}=0,\\
    -a^{\infty}(XV^{\infty}_{X}+YV^{\infty}_{Y})+U^{\infty}+ \frac{1}{2} |\nabla V^{\infty}|^2=0,
\end{array}\right.
   \end{equation*}

In an $L^{2}$ critical case, one expects an algebraic decrease of $L$ 
with $\tau$,   $L\propto 1/\tau^{\gamma}$, $\gamma\geq1$. With 
(\ref{dyn}) this leads to $L\propto 
(t^{*}-t)^{\gamma/(2\gamma-1)}$ where $t^{*}$ is the blow-up time 
(this follows from $\frac{d\tau}{dt}=1/L^{2}$ by 
integration). In 
the case of an exponential dependence of $L$ on $\tau$ expected in an 
$L^{2}$ supercritical case as here, 
$L\propto \exp(a^{\infty}\tau)$, we find $L\propto (t^{*}-t)^{1/2}$. 
Note that $a^{\infty}=0$ in the former case.
For the $L^{\infty}$ norm of $\eta$, we have a behavior 
proportional to $1/L^{2}$, 
and the same for the square of the 
$L^{2}$ norm of $\eta$. Obviously the rates coincide in the limit 
$\gamma\to\infty$, thus it will numerically only be possible to 
distinguish the case $\gamma\sim 1$ from an exponential dependence. 
Note that with our numerical approach, it is impossible to discover 
loglog corrections to the blow-up rate as known from the 
$L^{2}$-critical blow-up in NLS. Our goal is to distinguish an 
algebraic dependence of the scaling function $L$ from an exponential 
one. 

\section{Numerical construction of localised static solutions to the 
2D KBK system}\label{secstatic}
To numerically construct a localised solution to equation (\ref{Vr}), one has 
to find a non-trivial solution to a nonlinear ODE that is vanishing for 
$r\to\infty$ and regular for $r=0$, i.e., bounded with vanishing 
derivative $V_{r}$ for $r\to0$. The standard approach is to look for the solutions with an iterative scheme. The problem is that $V=0$ is a 
trivial solution to equation (\ref{Vr}), and that iterations will 
converge to this solution or not converge at all for general initial 
iterates. Therefore we first construct a solution in polar 
coordinates, which will then be interpolated to a 2D Fourier grid. 

\subsection{Polar coordinates}
We apply  the approach of 
\cite{CKS}: we introduce the variable $s=r^{2}$ in which equation 
(\ref{Vr}) becomes less singular,
\begin{equation}
	(4sV_{ss}+4V_{s}-V)_{s}+2sV_{s}^{3}=0
	\label{Vs}.
\end{equation}
In addition we work on the finite interval $s\in[0,s_{0}]$, where we
choose $s_{0}\gg 1$ 
such that the numerical solution for $V$ is vanishing for 
$s=s_{0}$ with numerical precision (we work in this paper always with 
double precision, i.e., with an accuracy of roughly $10^{-16}$). On the interval $[0,s_{0}]$ we apply the 
transform $s=s_{0}(1+l)/2$ with $l\in[-1,1]$ and for $l$ the standard 
Chebyshev discretisation $l_{j}=\cos(j\pi/N)$, $j=0,\ldots,N$ with $N$ 
a positive integer. In this setting derivatives are approximated via 
Chebyshev differentiation matrices \cite{trefethen,weireddy}. The 
vanishing condition for $s=s_{0}$ is implemented with a 
$\tau$-method, see \cite{trefethen}. 

\begin{figure}[htb!]
\includegraphics[width=0.49\textwidth]{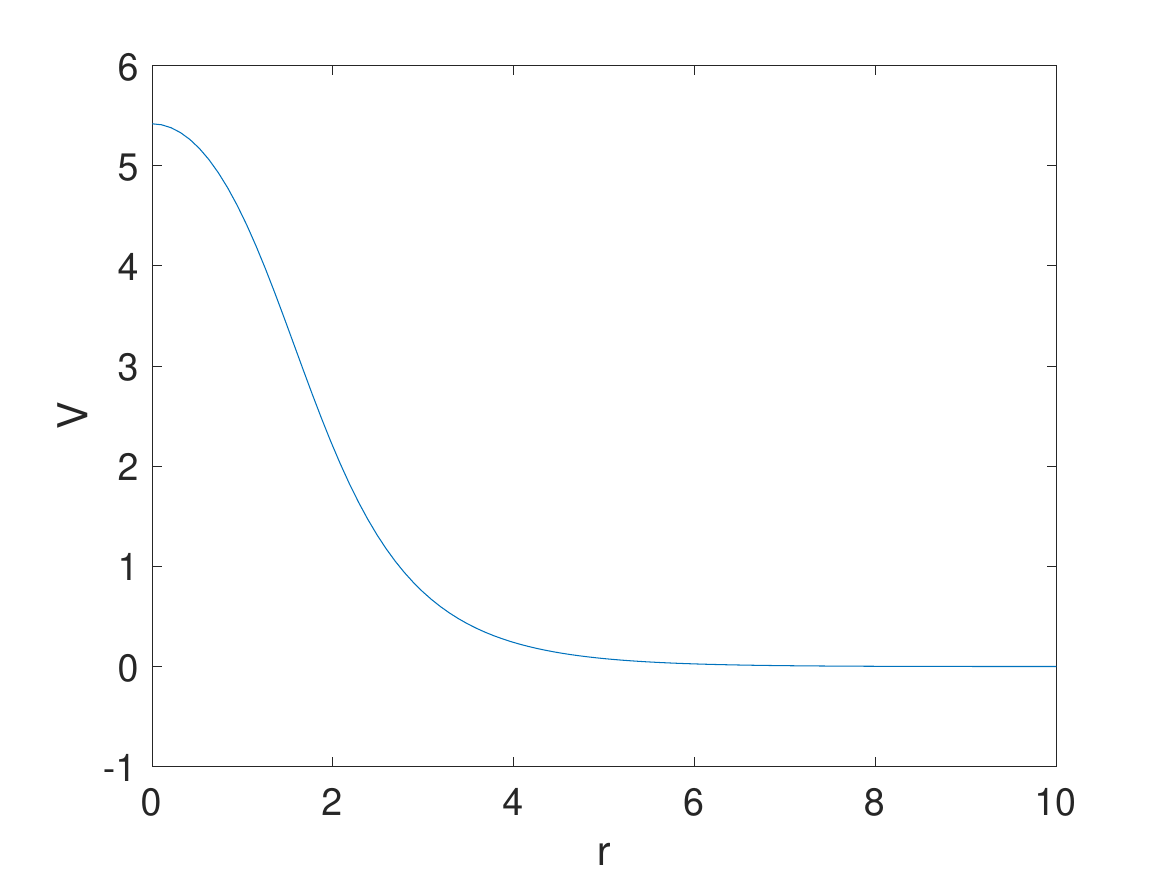}
\includegraphics[width=0.49\textwidth]{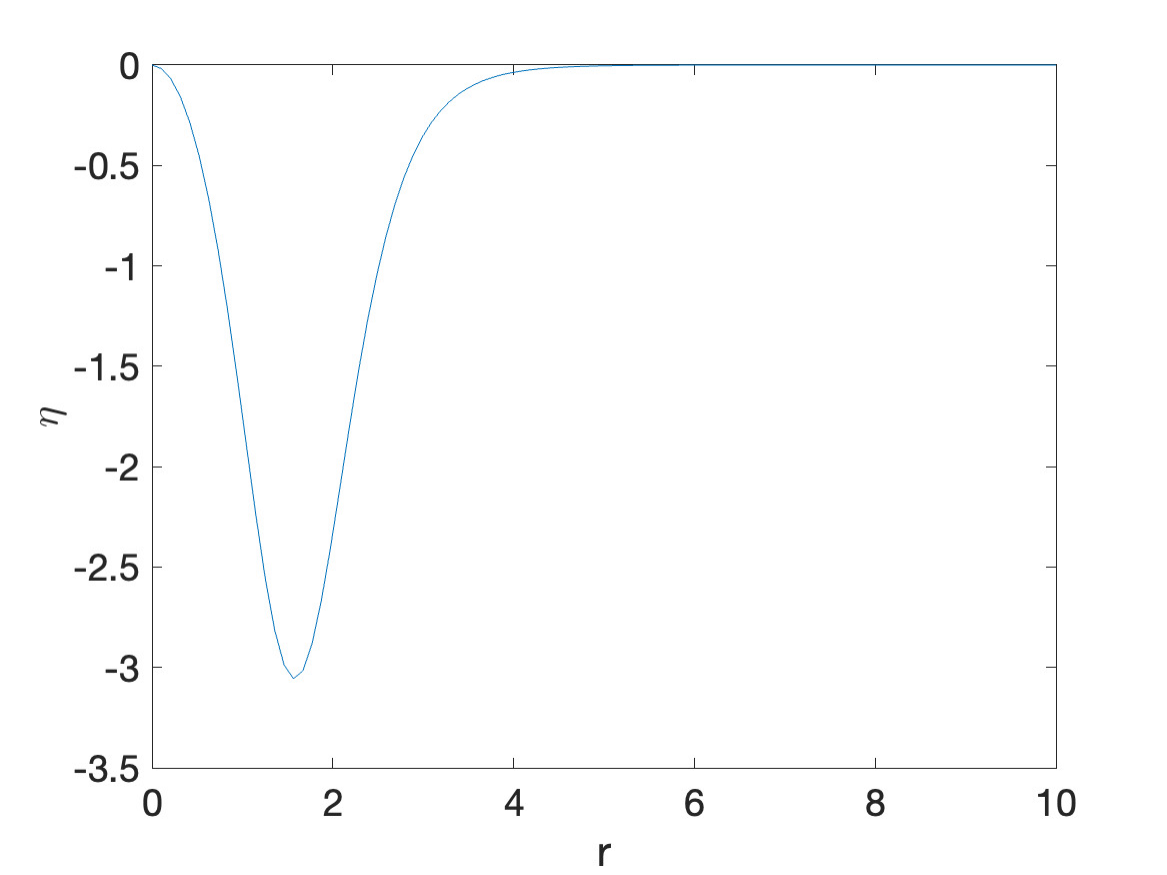}
\caption{Solution to equation (\ref{Vr}) in dependence of $r$ on the 
left and the corresponding solution for $\eta$ on the right. }
\label{Vradial}
\end{figure}

This 
leads for (\ref{Vs}) to a system of $N+1$ nonlinear equations that are 
solved iteratively with a Newton iteration. We use  $V=4\exp(-s)$ as 
an initial iterate. The iteration will not converge with this initial
iterate to a nontrivial solution, which is why we apply a relaxation: 
instead of accepting the $n+1$st 
iterate $V^{(n+1)}$ of the Newton iteration, we use $\mu 
V^{(n+1)}+(1-\mu)V^{(n)}$ for the next step, where $0<\mu<1$ 
(here we apply $\mu=0.1$). With this choice of $\mu$, the iteration 
with $s_{0}=100$ and $N=150$ converges (the iteration is stopped once 
the residual of (\ref{Vs}) drops below $10^{-9}$) to the solution shown in 
Fig.~\ref{Vradial} on the left. The corresponding solution $\eta$ is 
shown on the right of the same figure. 

Note that we were unable to construct non-trivial solutions vanishing 
at the origin, however, there may be a choice for the initial 
iterates for which this is possible. There may be 
also solutions with zeros for $V$, see the discussion in \cite{IW},
but to find these is beyond the 
scope of the current paper. 

\subsection{Fourier grid in 2D}\label{FFT2D}
In order to study the time evolution of solutions to the 2D KBK system 
and perturbations of exact solutions such as the static solution 
constructed above, we will not limit ourselves to radially symmetric 
situations. To treat the general case, we consider solutions on a 
2-torus $\mathbb{T}_{2\pi L_{x}}\times\mathbb{T}_{2\pi L_{y}}$ with 
$L_{x},L_{y}>0$. This means we work in a doubly periodic setting 
where we approximate rapidly decreasing functions on large enough 
tori as essentially periodic within the finite numerical precision. 
We introduce the standard discretisation of the discrete Fourier 
transform (DFT) on this 2-torus in each direction, and compute the 
DFT with a fast Fourier transform (FFT). 

\begin{figure}[htb!]
\includegraphics[width=0.49\textwidth]{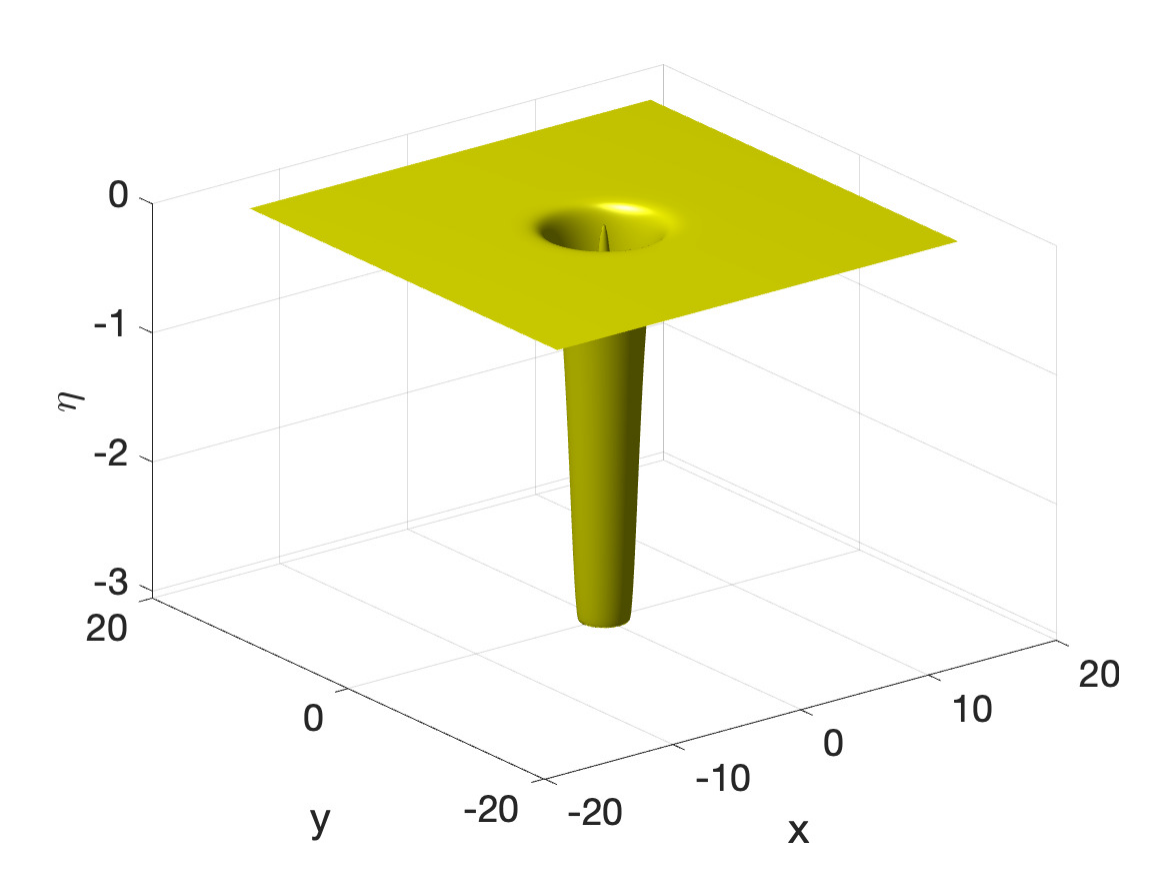}
\includegraphics[width=0.49\textwidth]{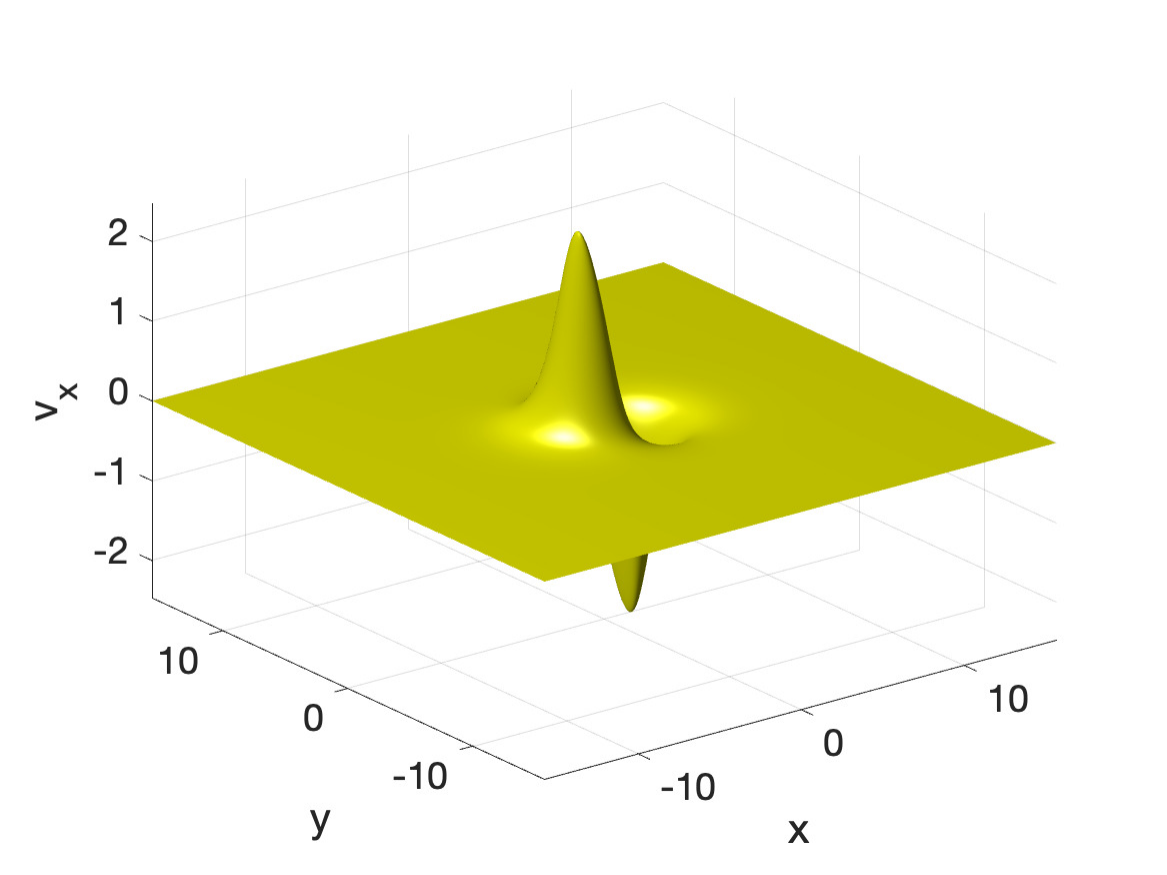}
\caption{Static solution $\eta$ to the system KBK (\ref{Kaup-true-2D}) on the 
left and $v_{x}$ on the right. }
\label{figstatic}
\end{figure}

To treat the static solution on this Fourier grid, 
we interpolate the solution in polar coordinates to the latter and 
solve as in \cite{SGN2d} the 
resulting system of nonlinear equations for (\ref{Vst}) with a Newton-Krylov method. 
This leads to the solutions shown  in Fig.~\ref{figstatic}. The 
solutions $V$ and $\eta$ are simply the rotated 
versions of the radially symmetric solutions in 
Fig.~\ref{Vradial}. The function $v_{x}$ ($v_{y}$ is just $v_{x}$ rotated by 
90 degrees for symmetry reasons) is just the radial derivative of 
$V$. We note that the $L^{2}$ norm of  $\eta$ is $||\eta||_{2}^{2}\sim  96.5967$.

Note that the iteration with the same initial iterate as in the 
radially symmetric case will converge (even with relaxation) to a 
solution of the scalar fractional NLS equation. 

\section{Numerical approach for the time evolution}
In this section we present the approach for the numerical time 
evolution for the KBK solutions, a 2D variant of the approach in 
\cite{KS-KBK}.

We apply the same discretisation of the spatial dependence as in 
section \ref{FFT2D}. The Fourier transform in the system 
(\ref{upmt2}) is approximated by a DFT which leads to a $2N_{x}N_{y}$ 
dimensional system of ordinary differential equations in $t$. In an 
abuse of notation we denote the DFT of a function $u$ as the Fourier 
transform by $\hat{u}$. Thus the system (\ref{upmt2}) is of the 
form 
\begin{equation*}
	(\hat{u}_{\pm})_{t}=\mathcal{L}_{\pm}\hat{u}_{\pm}+\mathcal{N}_{\pm}
	\label{LU},
\end{equation*}
where $\mathcal{L}_{\pm}=\mp i|\mathbf{k}|\sqrt{1+|\mathbf{k}|^{2}}$ and where 
$\mathcal{N}_{\pm}$ corresponds to the nonlinear dependence in (\ref{upmt2}). 
Since we will need large values for $N_{x}$ and $N_{y}$, this is a
classical example of a stiff system which means that explicit time 
integration schemes will be inefficient due to stability conditions,
see for instance the discussion in \cite{HO} and references therein. 
As in \cite{KS-KBK} we use an \emph{exponential time differencing 
scheme} (ETD) of classical order 4 by Cox and Matthews \cite{CM}, see 
\cite{HO} and \cite{KR} for a detailed discussion. 
 Since the method is explicit (only information of the 
functions at the previous time step is needed), it is not important that 
the nonlinear part is not diagonal in (\ref{upmt2}) as the linear part. 
As discussed in \cite{etna,KR}, the exactly conserved energy can be 
used to control the numerical error in the time integration since 
this is not guaranteed by the time integration scheme.  The relative energy
\begin{equation*}
	\Delta E :=|E(t)/E(0)-1|
	\label{rel}
\end{equation*}
typically overestimates the numerical error by 1-2 orders of 
magnitude. We will always aim for a $\Delta E$ considerably smaller than 
$10^{-3}$.

We first test the code for the example of the static solution. To this end 
we use $L_{x}=L_{y}=5$, $N_{x}=N_{y}=2^{9}$ with $N_{t}=10^{3}$ time 
steps for $t\leq10$. The difference between the numerical solution for 
$t=10$ and the initial data is of the order of $10^{-12}$. Note that 
an exact static solution would not test the time evolution, but since 
the solution is numerically constructed, this tests whether the 
static solution is well enough approximated that no time dependence 
of the solution is observed. Thus we can also conclude that the numerical 
errors with which the static solution is obtained must be smaller 
than $10^{-12}$. 

\begin{figure}[htb!]
\includegraphics[width=0.49\textwidth]{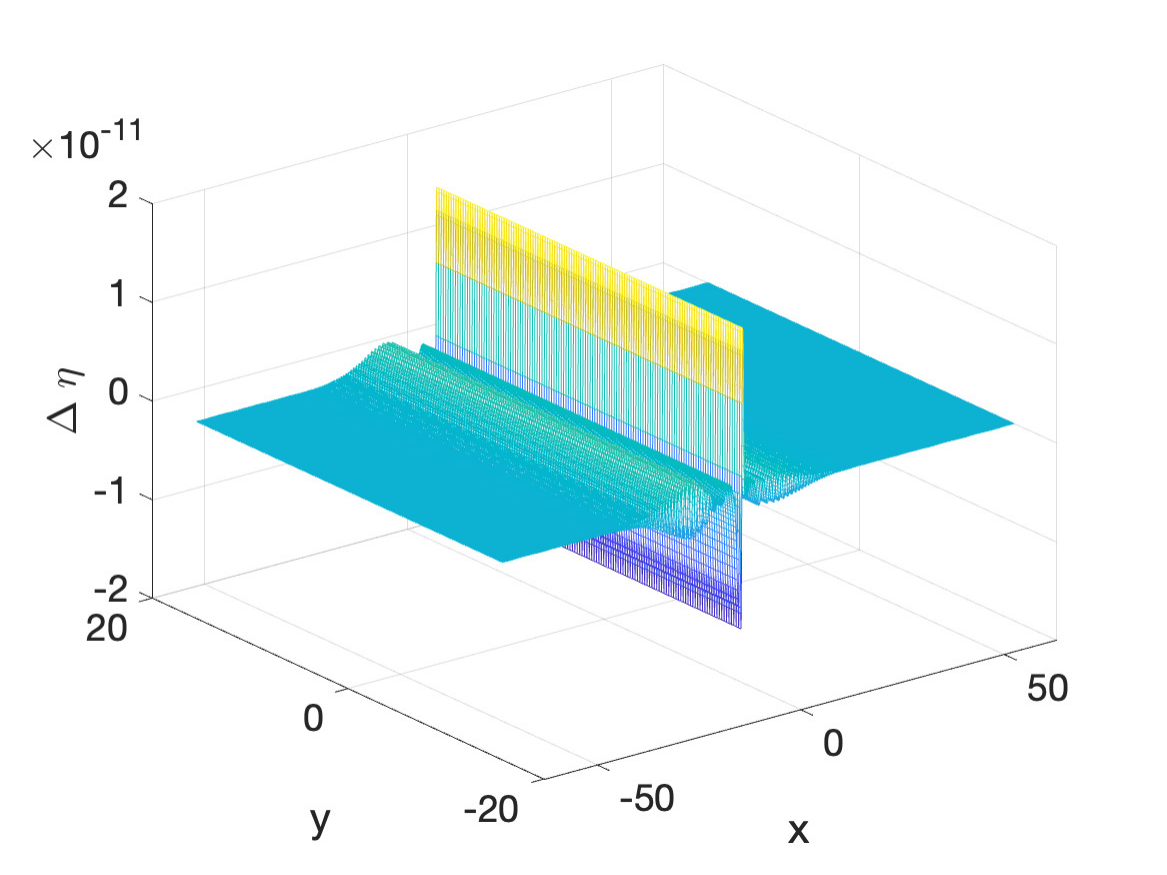}
\includegraphics[width=0.49\textwidth]{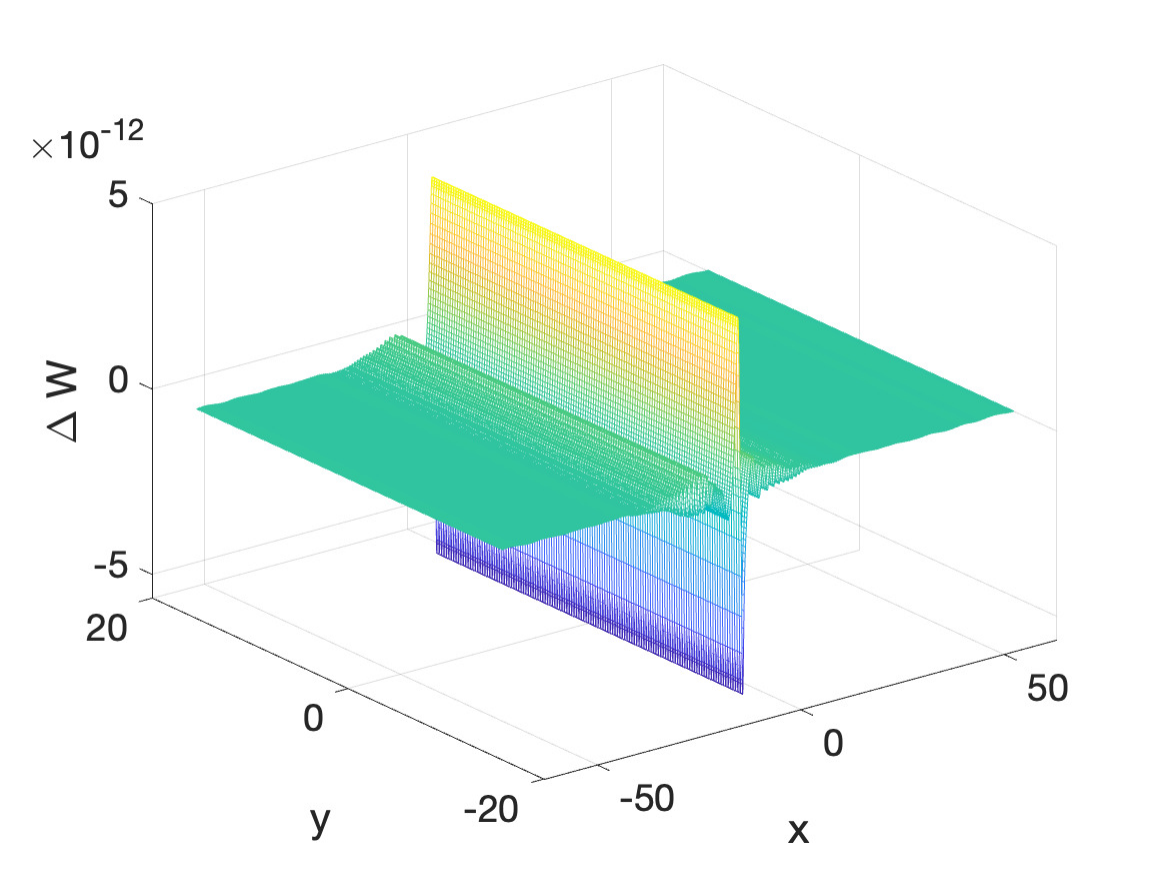}
\caption{Difference of the KBK solution for line solitary wave initial 
data (\ref{soliton}) for $C=0.8$ and the exact solution for $t=1$, on the 
left $\eta$, on the right $W$. }
\label{linetest}
\end{figure}

To test a time dependent exact solution, we consider the exactly 
known line solitary wave. This does not test the $y$-dependence since 
the solution is not dependent on this variable, however this was already achieved 
with the test of the static solution. We use $L_{x}=20$, $L_{y}=5$, 
$N_{x}=2^{11}$, $N_{y}=2^{7}$ and $N_{t}=2000$ time steps for $t\leq 
1$. As an example we consider the line solitary wave (\ref{soliton}) 
for $C=0.8$ as exact initial data. We show the difference between the 
numerical and the initial data shifted by a factor $C$ in $x$ in 
Fig.~\ref{linetest}. It can be seen that the difference is of order 
$10^{-11}$. This is a  non-trivial result because, as it will be 
shown in section \ref{line}, the line solitary wave for this value of 
$C$ is strongly 
unstable. The numerically computed energy is conserved to better than 
$10^{-12}$ during the whole computation.

These tests show that the code is able to solve the KBK system to 
essentially machine precision even for initial data which can produce a blow-up under small  perturbation. 

\section{Perturbations of the static solution}\label{secpertstatic}
In this section we study perturbations of the static solution 
constructed in section \ref{secstatic}. First, we consider a perturbation of the form 
\begin{equation}
	\eta(x,y,0) = \lambda\eta_{static}(x,y),\quad V(x,y,0) = V_{static}
	\label{inistatic},
\end{equation}
where $\lambda\sim1$. 

\begin{figure}[htb!]
\includegraphics[width=0.49\textwidth]{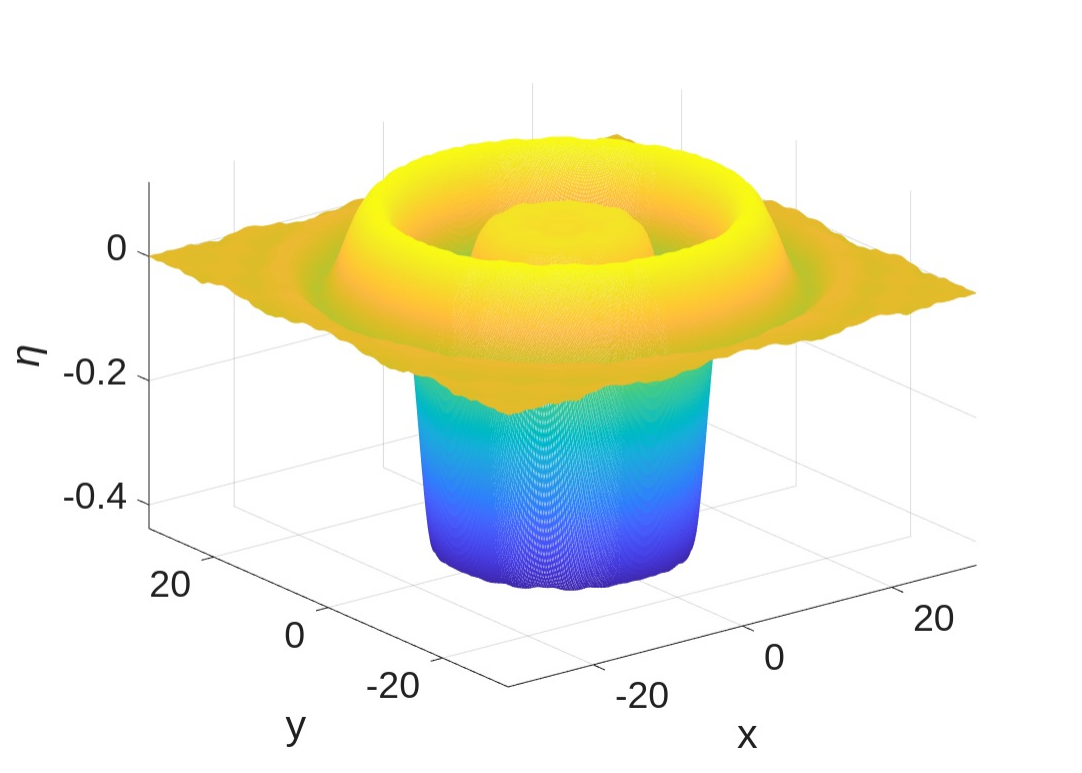}
\includegraphics[width=0.49\textwidth]{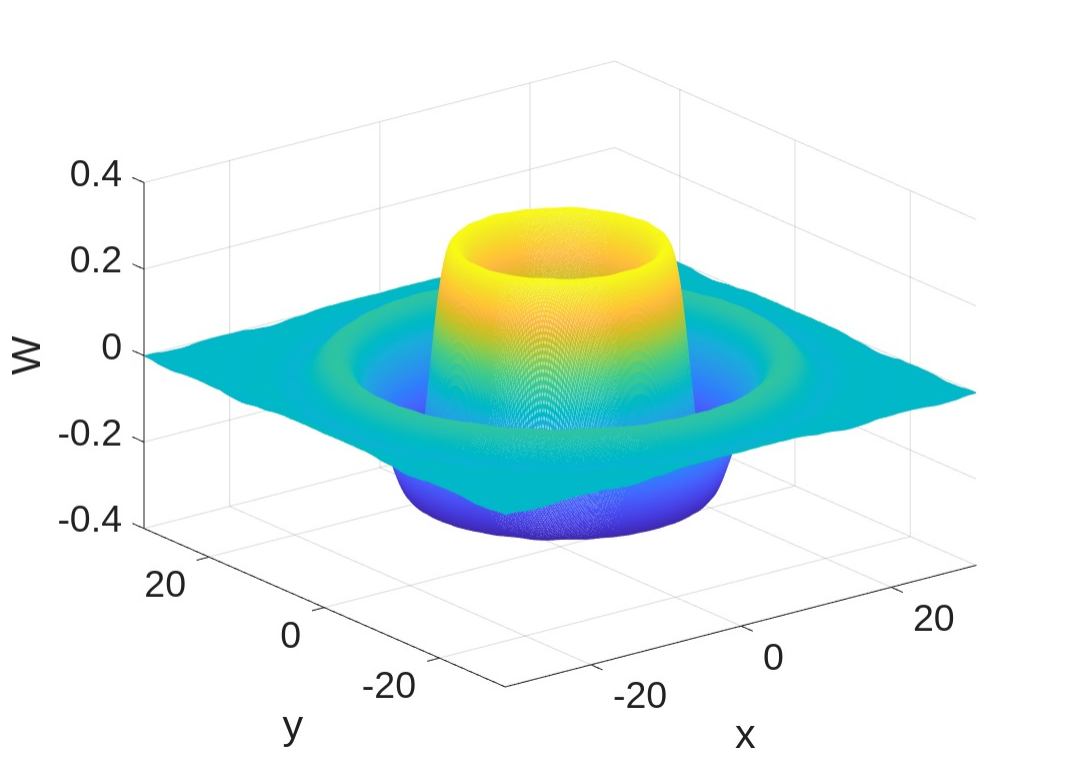}
\caption{Solution to the system KBK (\ref{Kaup-true-2D}) for the 
initial data (\ref{inistatic}) with $\lambda=0.99$ for $t=20$}, on the 
left $\eta$, on the right $W$. 
\label{static099eta}
\end{figure}

For $\lambda=0.99$, a perturbation with a smaller $L^{2}$ norm of 
$\eta$ than the static solution, we apply  $L_{x}=L_{y}=10$, 
$N_{x}=N_{y}=2^{9}$ with $N_{t}=10^{4}$ time 
steps for $t\leq20$. We find that the solution is simply dispersed as 
shown by the solution at the final time in Fig.~\ref{static099eta}.

\begin{figure}[htb!]
\includegraphics[width=0.49\textwidth]{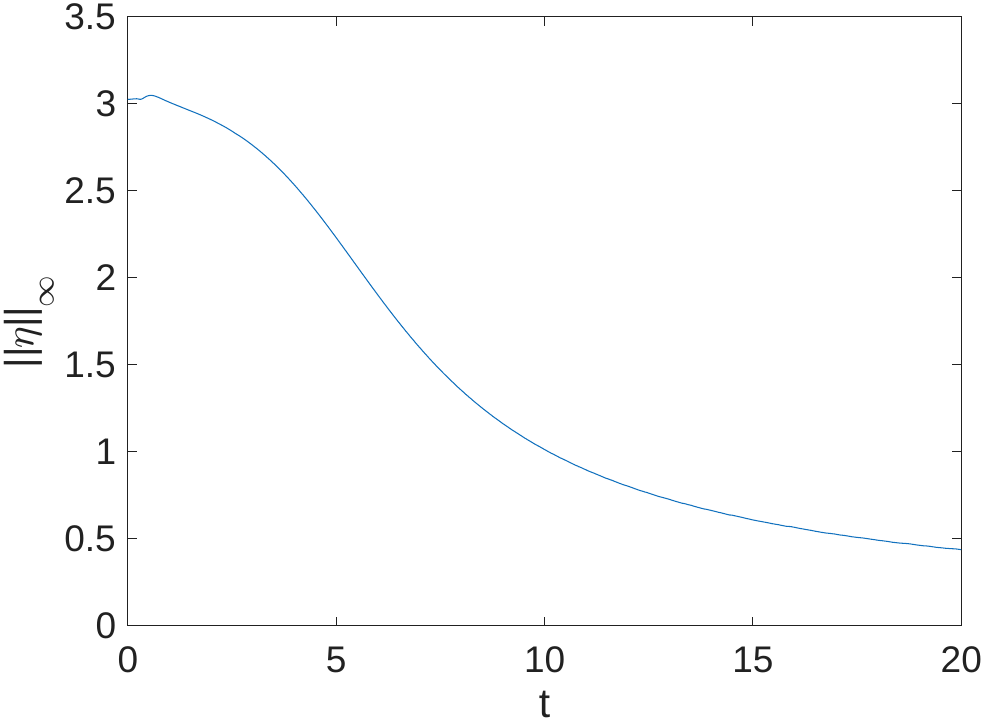}
\includegraphics[width=0.49\textwidth]{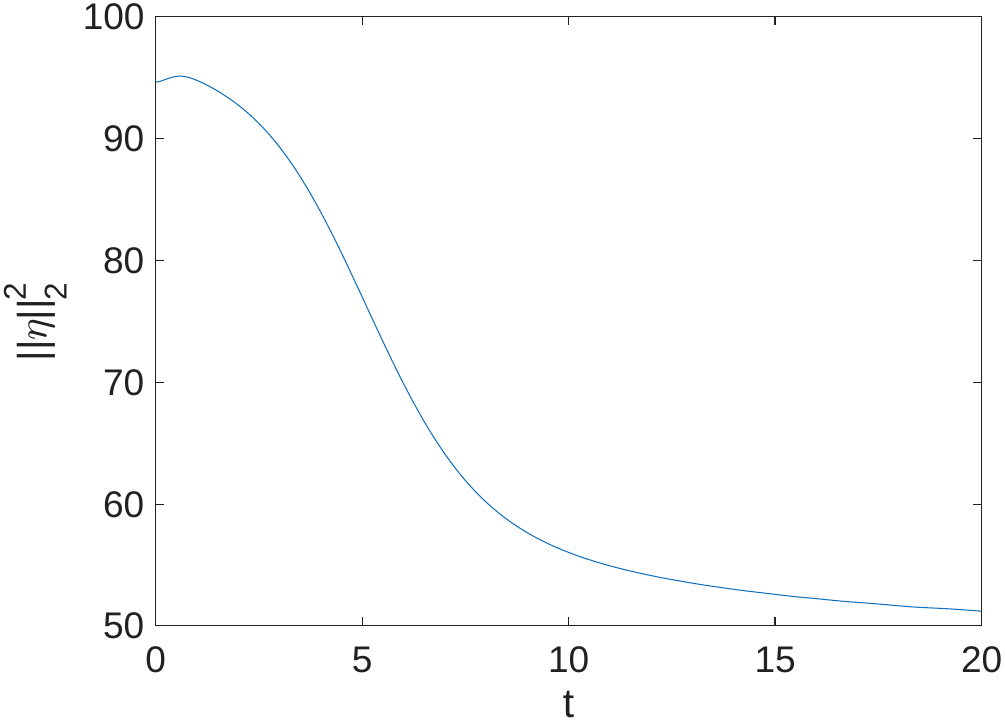}
\caption{Norms of the solution to the system KBK (\ref{Kaup-true-2D}) for the 
initial data (\ref{inistatic}) with $\lambda=0.99$ in dependence of 
time, on the 
left $||\eta||_{\infty}$, on the right $||\eta||_{2}^{2}$. }
\label{static099etainf}
\end{figure}

The $L^{\infty}$ norm and the $L^{2}$ norm of $\eta$ shown in 
Fig.~\ref{static099etainf} confirm that no stable structure appears 
in the time evolution. Apparently the initial data are simply 
dispersed to infinity. 

We now consider the case $\lambda=1.01$, a perturbation with a larger 
$L^{2}$ norm than the static solution. 
We apply  $L_{x}=L_{y}=5$, $N_{x}=N_{y}=2^{11}$ with $N_{t}=2 
\times10^{4}$ time steps for $t\leq5.2$. The numerical simulation breaks down at $t \approx 5.1740$.

\begin{figure}[htb!]
\includegraphics[width=0.49\textwidth]{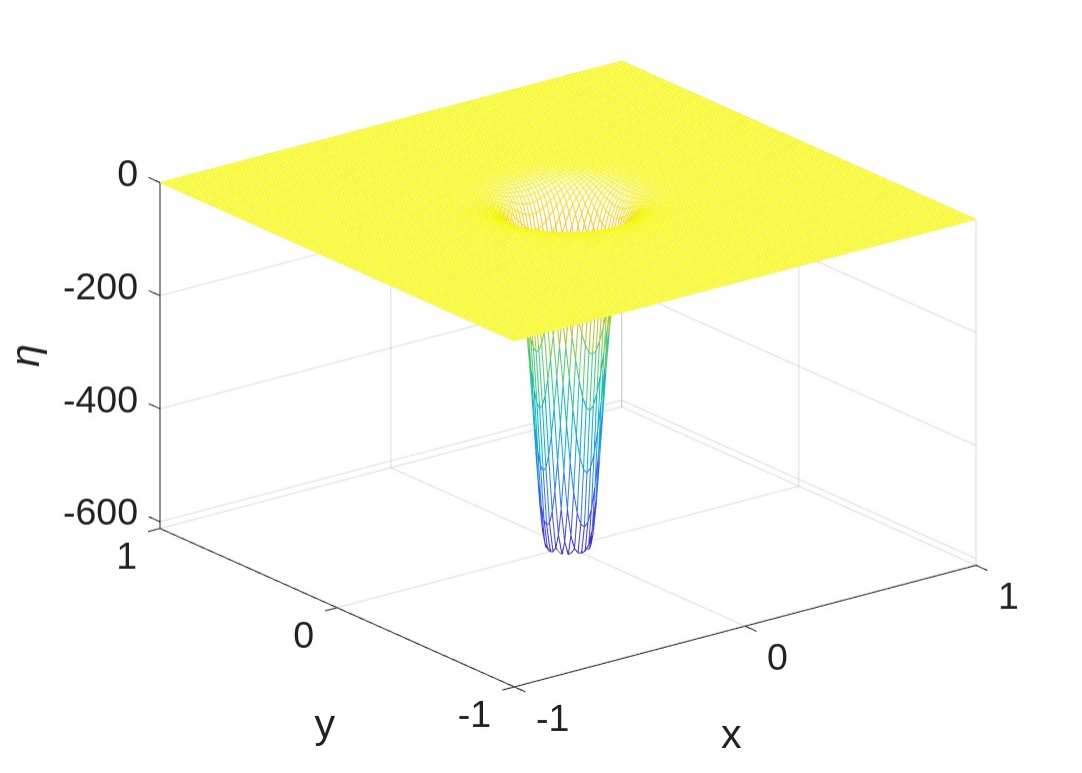}
\includegraphics[width=0.49\textwidth]{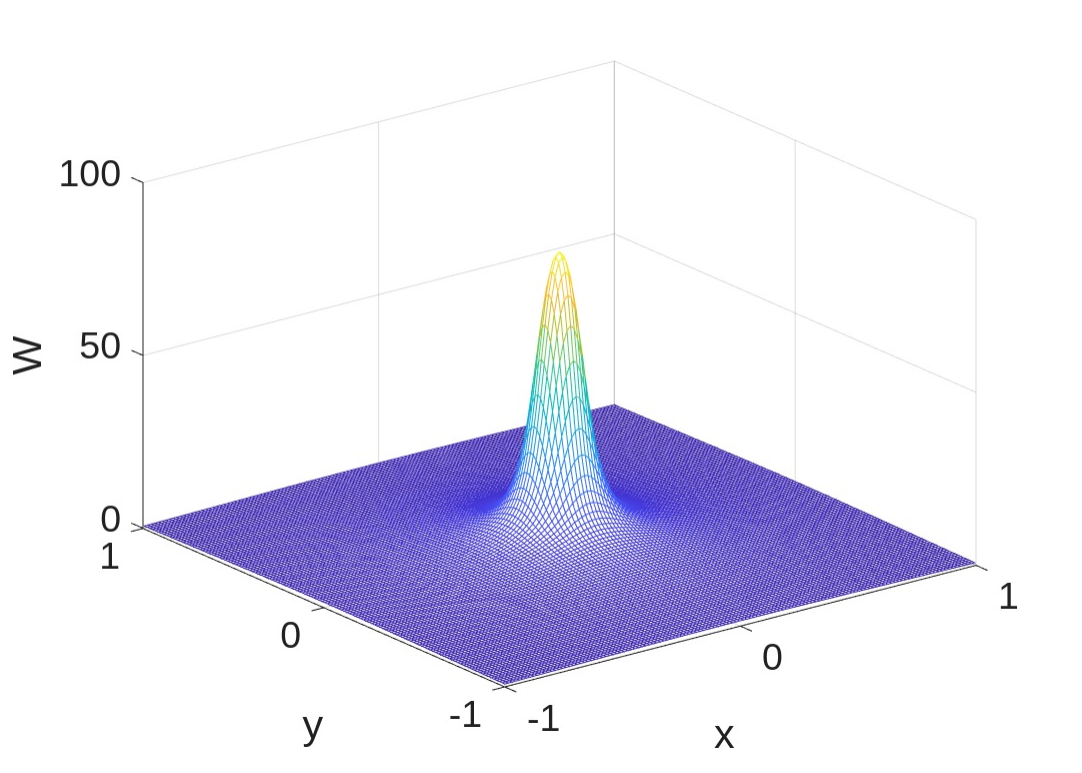}
\caption{Solution to the system KBK (\ref{Kaup-true-2D}) for the 
initial data (\ref{inistatic}) with $\lambda=1.01$ for $t=5.1740$ in 
a close-up, on the 
left $\eta$, on the right $W$. }
\label{static101eta}
\end{figure}

In contrast to the dispersive regime observed for $\lambda = 0.99$, 
the perturbation now undergoes a rapid self-focusing process, as can 
be seen in Fig.~\ref{static101eta}.  The amplitude of $\eta$ becomes 
increasingly sharper while its spatial support shrinks. Both the 
$L^\infty$ and $L^2$ norm exhibit strong growth as the solution approaches the blow-up time.

We numerically observe that near blow-up, the $L^\infty$ and $L^2$ norms behave as
\[
\|\eta\|_\infty \sim (t^* - t)^\alpha, \qquad 
\|\eta\|_2 ^2 \sim (t^* - t)^\alpha,
\]
with $\alpha < 0$. We use the Matlab algorithm \emph{fminsearch} to 
fit $\ln \|\eta\|_2 ^2$ and $\ln \|\eta\|_\infty$ to $\alpha \ln(t^* 
- t) + \beta$ over the last $1000$ time steps and obtain:
\begin{itemize}
    \item for $\|\eta\|_2 ^2$: $\alpha = -0.8809$, $\beta = 4.4756$, $t^* = 5.1768$,
    \item for $\|\eta\|_\infty$: $\alpha = -0.9177$, $\beta = 1.1564$, $t^* = 5.1772$.
\end{itemize}

\begin{figure}[htb!]
    \centering
    \begin{subfigure}[b]{0.49\textwidth}
        \centering
        \includegraphics[width=\textwidth]{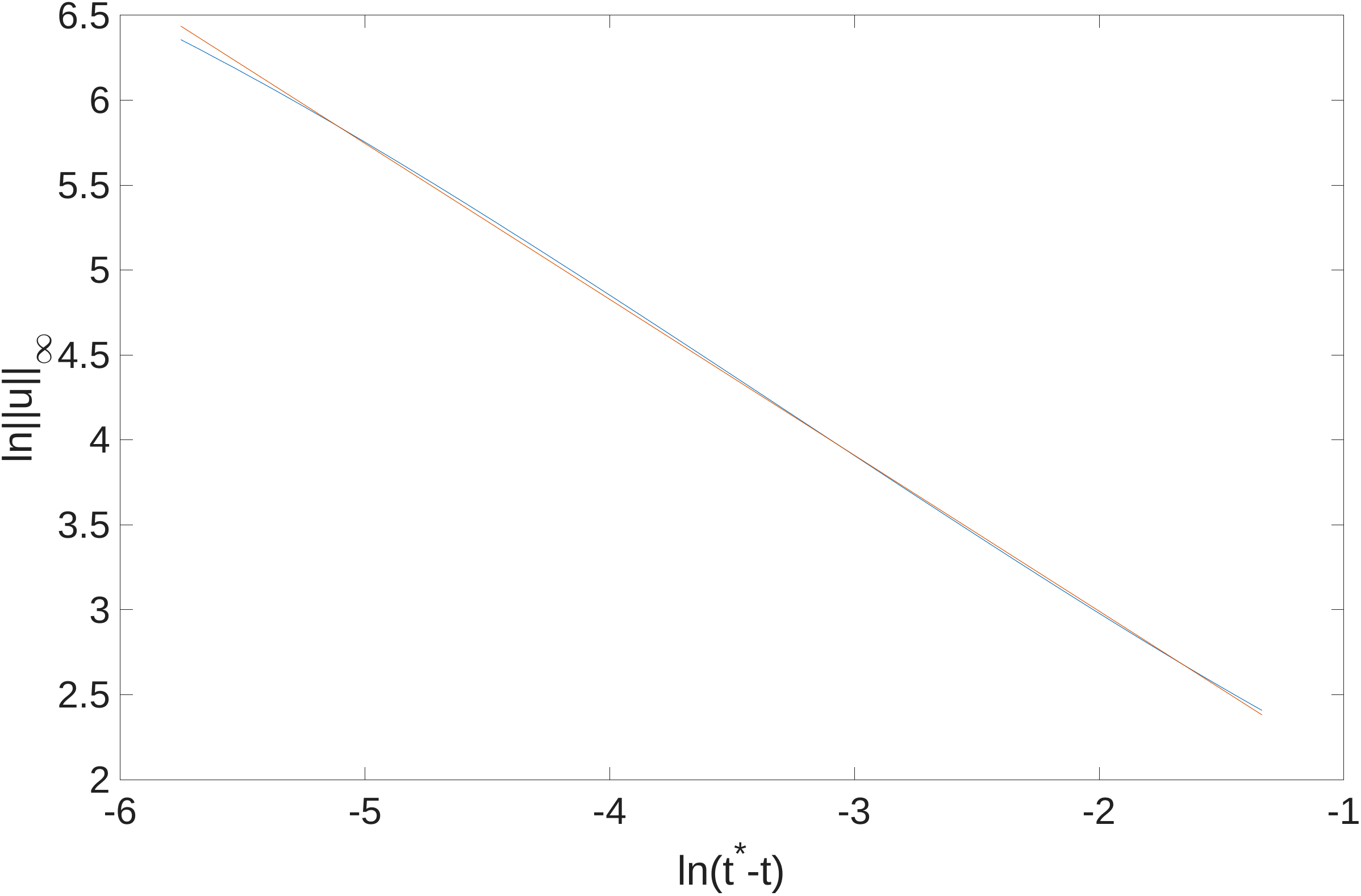}
    \end{subfigure}
    \hfill
    \begin{subfigure}[b]{0.5\textwidth}
        \centering
        \includegraphics[width=\textwidth]{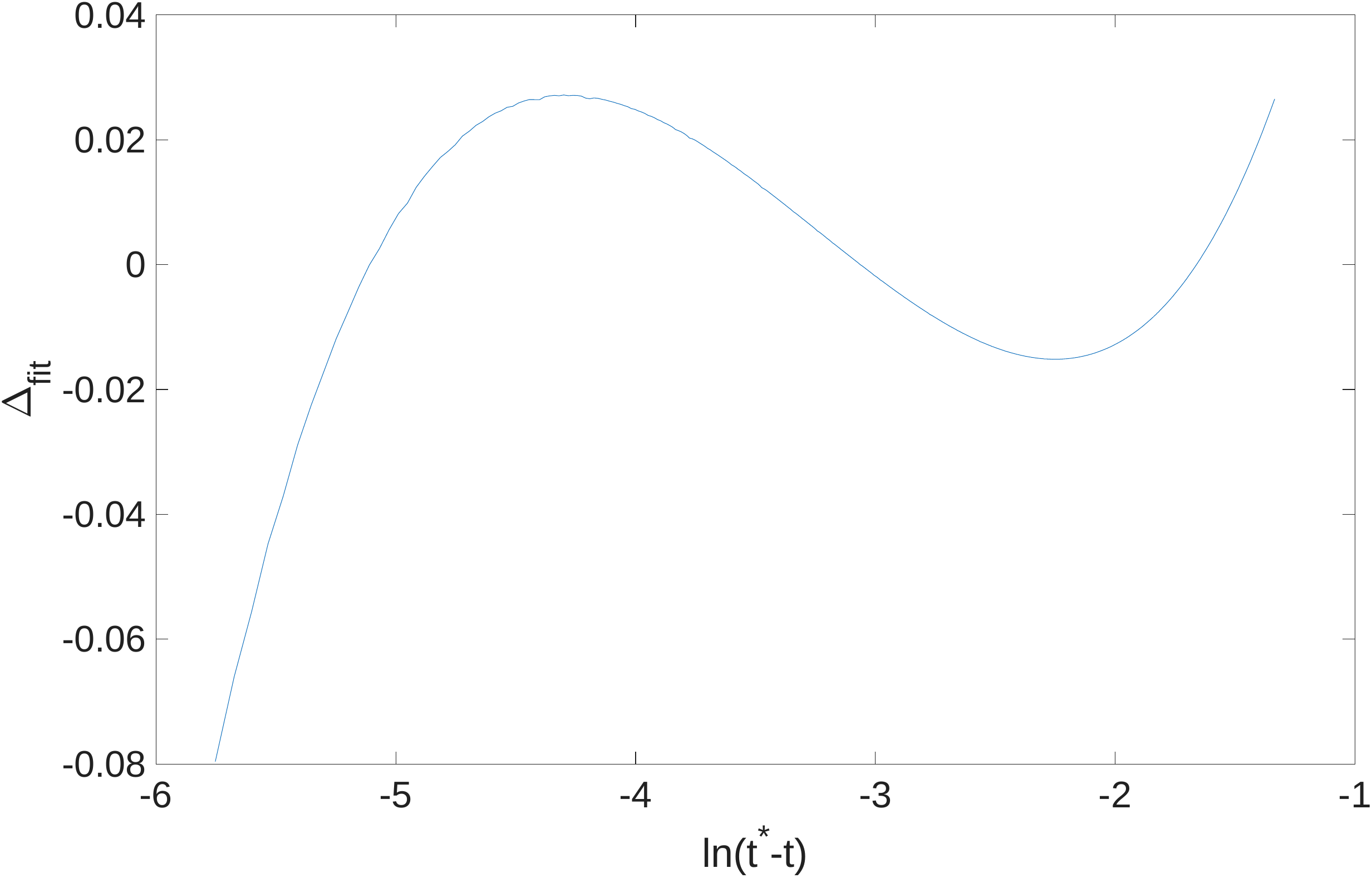}
    \end{subfigure}

    \begin{subfigure}[b]{0.49\textwidth}
        \centering
        \includegraphics[width=\textwidth]{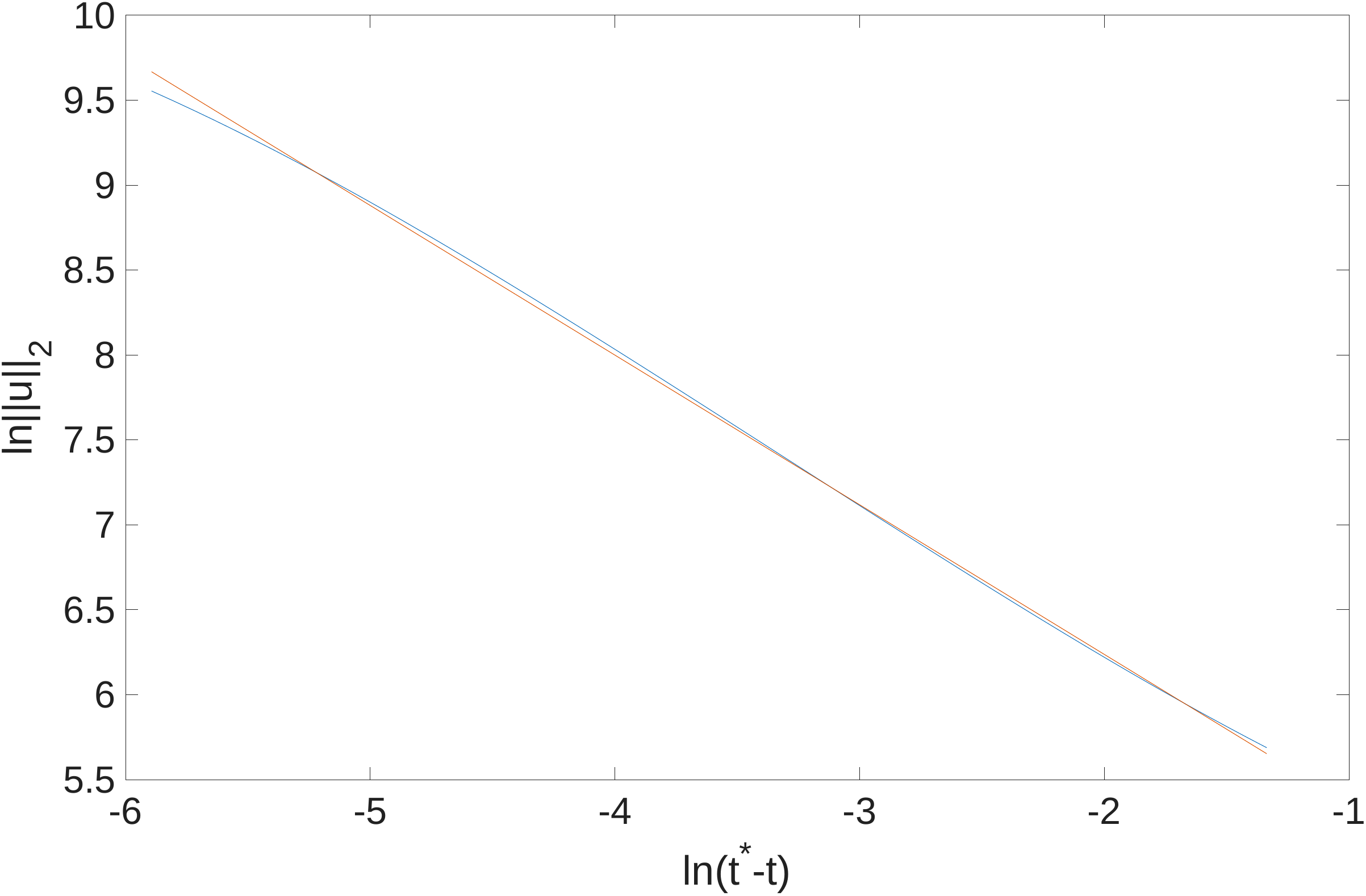}
    \end{subfigure}
    \hfill
    \begin{subfigure}[b]{0.5\textwidth}
        \centering
        \includegraphics[width=\textwidth]{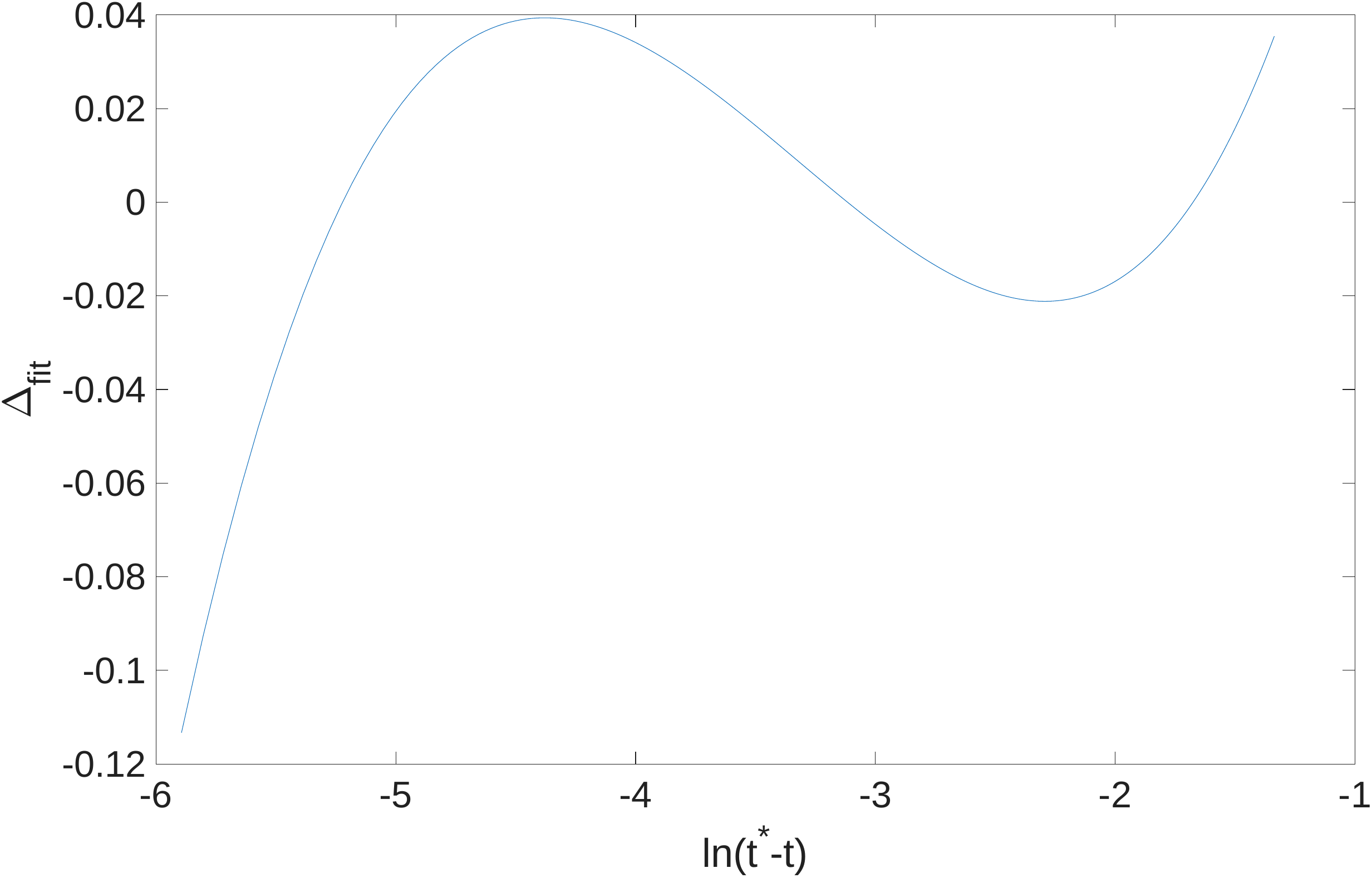}
    \end{subfigure}

    \caption{On top the logarithm of the $L^\infty$ norm of the last 
	$1000$ time steps of the solution to the system KBK 
	(\ref{Kaup-true-2D}) for the 
initial data (\ref{inistatic}) with $\lambda=1.01$ and the fitted 
curve $\alpha \ln(t^* - t) + \beta$ on the left, and the difference 
$\Delta_{fit}$ (`fitting error') between both 
	curves on the right. Analogous plots on the bottom for the 
    $L^2$ norm.}
    \label{fig:combined_norms_stat}
\end{figure}

The fitting error $\Delta_{fit}$, i.e., the relative difference between the 
logarithm of the considered norm and both fitted models is  
about $1\%$ (see Fig.~\ref{fig:combined_norms_stat}), and the fitted 
parameters remain stable when varying the fitting window or the 
initial guess. The agreement between the fitted values for the 
$L^{\infty}$ norm, which is locally determined, and the $L^{2}$ norm 
being computed on the whole considered interval indicates that the 
results are stable. However, it is numerically difficult to decide 
whether one is close enough to the actual blow-up that the computed 
norms already catch the asymptotic behavior. With this caveat, the 
obtained values of $\alpha$ are compatible with the theoretically 
predicted exponent $-1$ for an exponential dependence of the scaling 
factor $L$ on $\tau$.

\begin{figure}[htb!]
\includegraphics[width=0.49\textwidth]{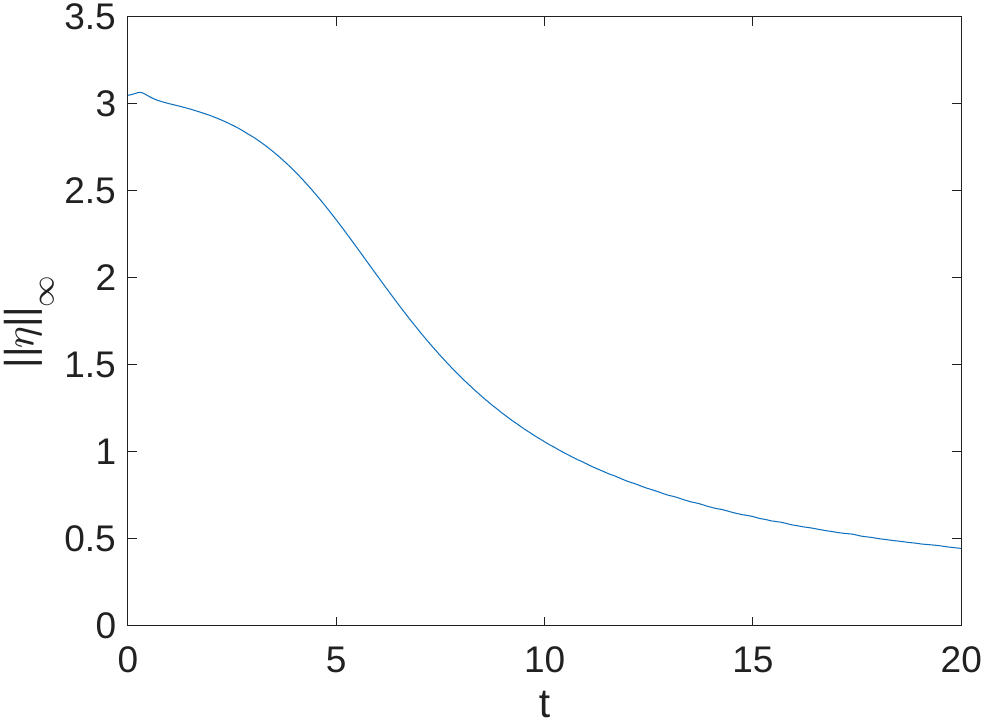}
\includegraphics[width=0.49\textwidth]{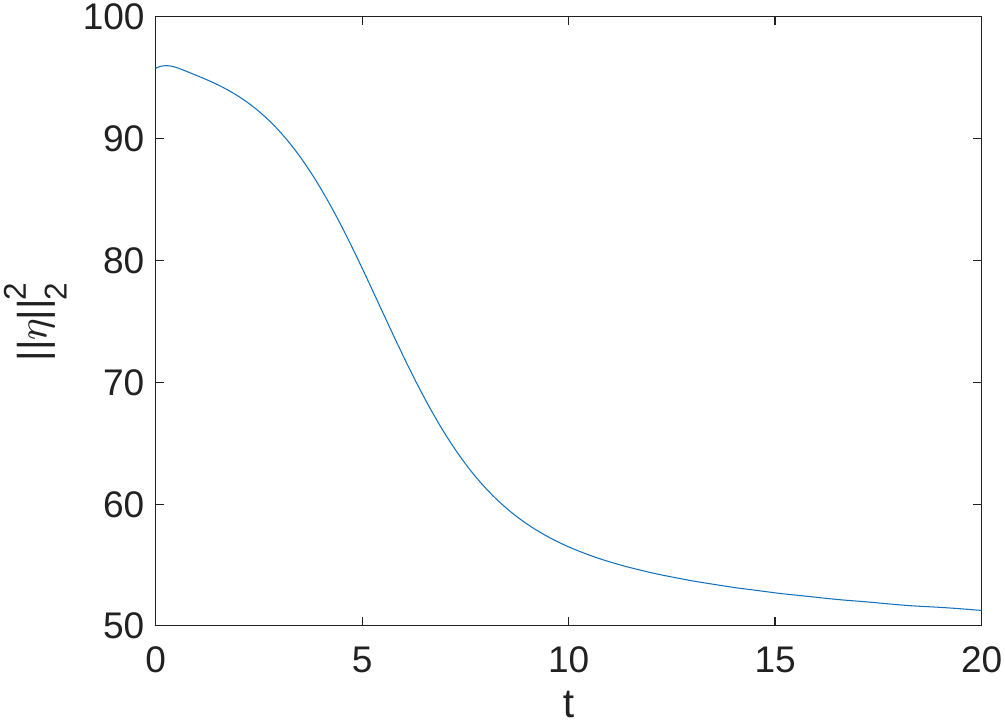}
\caption{Norms of the solution to the system KBK (\ref{Kaup-true-2D}) for the 
initial data (\ref{inistatic2}) with $\mu=+0.1$ in dependence of 
time, on the 
left $||\eta||_{\infty}$, on the right $||\eta||_{2}^{2}$. }
\label{normstatpg}
\end{figure}

We introduce perturbations around the static solution of the form:
\begin{equation}
	\eta(x,y,0) = \eta_{static}(x,y)\pm \mu \exp(-x^{2}-y^{2}),\quad V(x,y,0) = V_{static}
	\label{inistatic2}.
\end{equation}
 We use $L_x = L_y = 10$, $N_x = N_y = 2^{9}$, and $N_t = 1 \times 10^4$ time steps up to $t \leq 20$. For $\mu = +0.1$, corresponding to a perturbation with smaller $L^2$ norm of $\eta$ than the static solution, the solution is once again dispersed to infinity, as indicated by the decay of both the $L^\infty$ and $L^2$ norms of $\eta$ (see Fig.~\ref{normstatpg}).

We now consider the case where $\mu =-0.1$, a perturbation with 
bigger $L^2$ norm of $\eta$ than the static solution. We use $L_x = 
L_y = 5$, $N_x = N_y = 2^{11}$, and $N_t = 2 \times 10^4$ time steps 
up to $t \leq 6$. The computation breaks down at $t\approx 5.46$. 
Once again $\eta$ shows focusing behavior and both $L^2$ and $L^\infty$ norms grow strongly as the blow-up time is approached.
Note that the blow-up occurs for negative values of 
$1+\eta$ which should not appear in the context of water waves.

The fitting of $\ln \|\eta\|_2$ and $\ln \|\eta\|_\infty$ to $\alpha 
\ln(t^* - t) + \beta$ over the last $1000$ time steps, gives:
\begin{itemize}
    \item for $\|\eta\|_2 ^2$: $\alpha = -0.7839$, $\beta = 4.6157$, $t^* = 5.4477$,
    \item for $\|\eta\|_\infty$: $\alpha = -0.8419$, $\beta = 1.2640$, $t^* = 5.4499$.
\end{itemize}

\begin{figure}[htb!]
    \centering
    \begin{subfigure}[b]{0.49\textwidth}
        \centering
        \includegraphics[width=\textwidth]{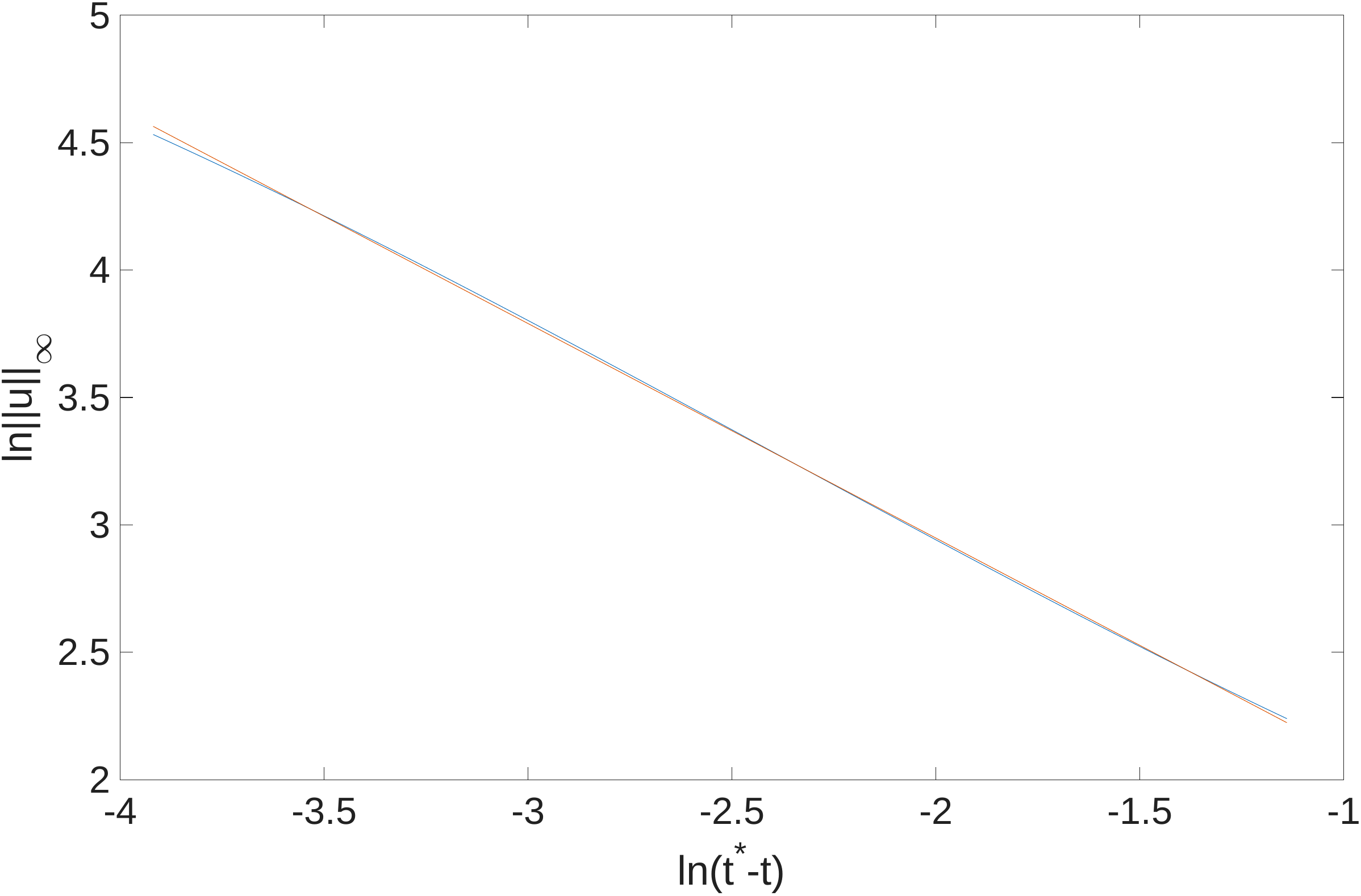}
    \end{subfigure}
    \hfill
    \begin{subfigure}[b]{0.5\textwidth}
        \centering
        \includegraphics[width=\textwidth]{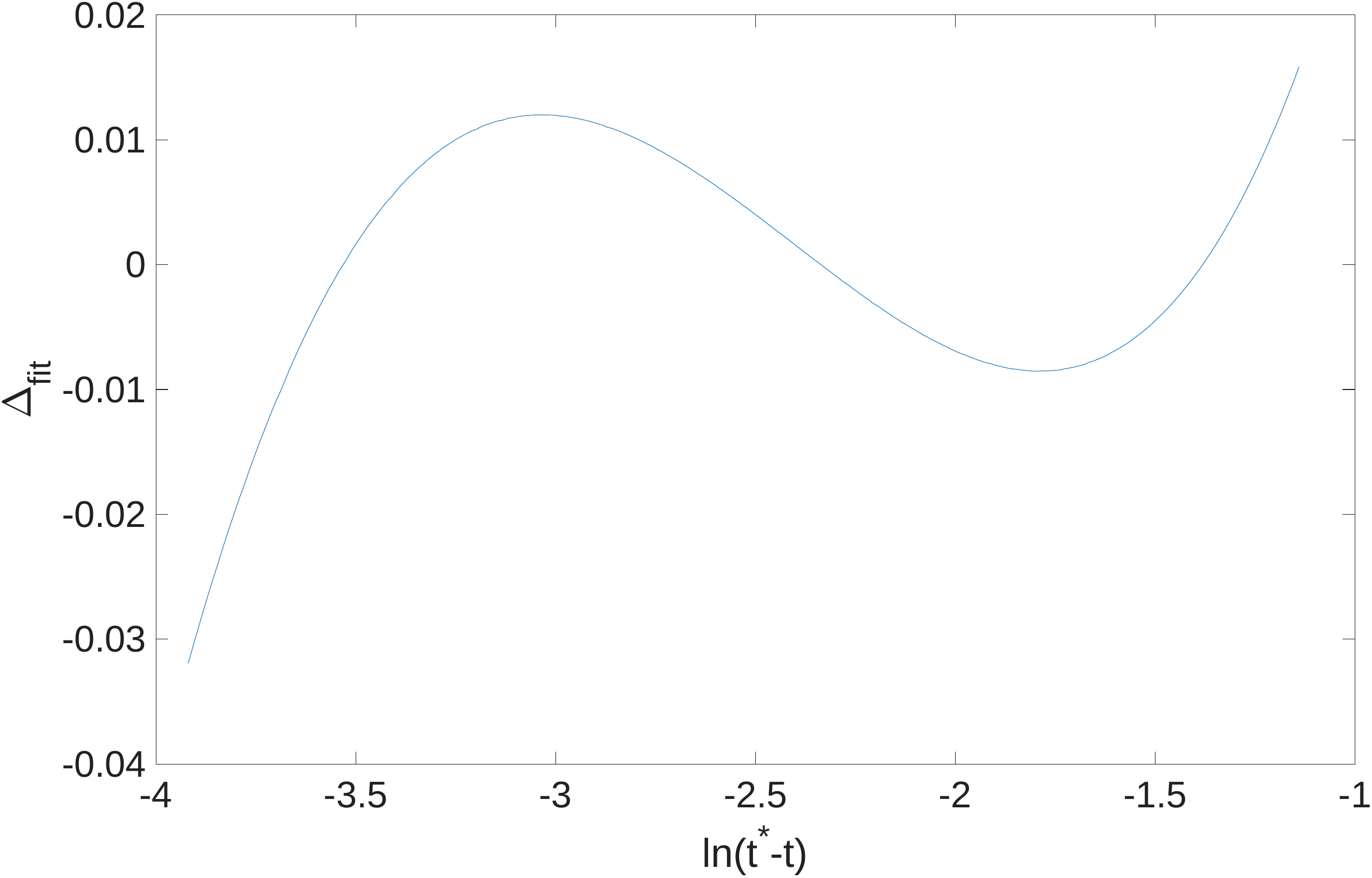}
    \end{subfigure}

    \begin{subfigure}[b]{0.49\textwidth}
        \centering
        \includegraphics[width=\textwidth]{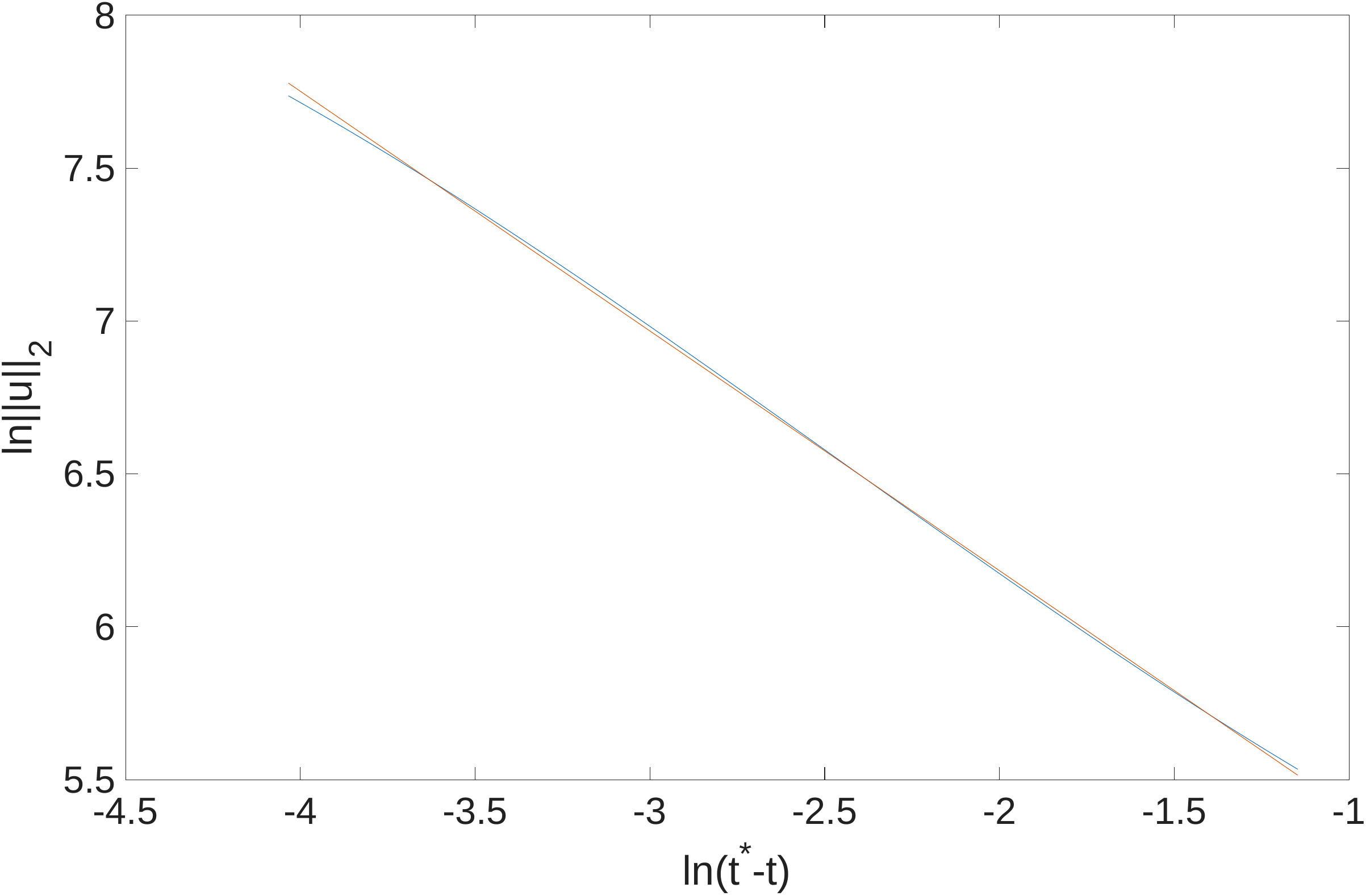}
    \end{subfigure}
    \hfill
    \begin{subfigure}[b]{0.5\textwidth}
        \centering
        \includegraphics[width=\textwidth]{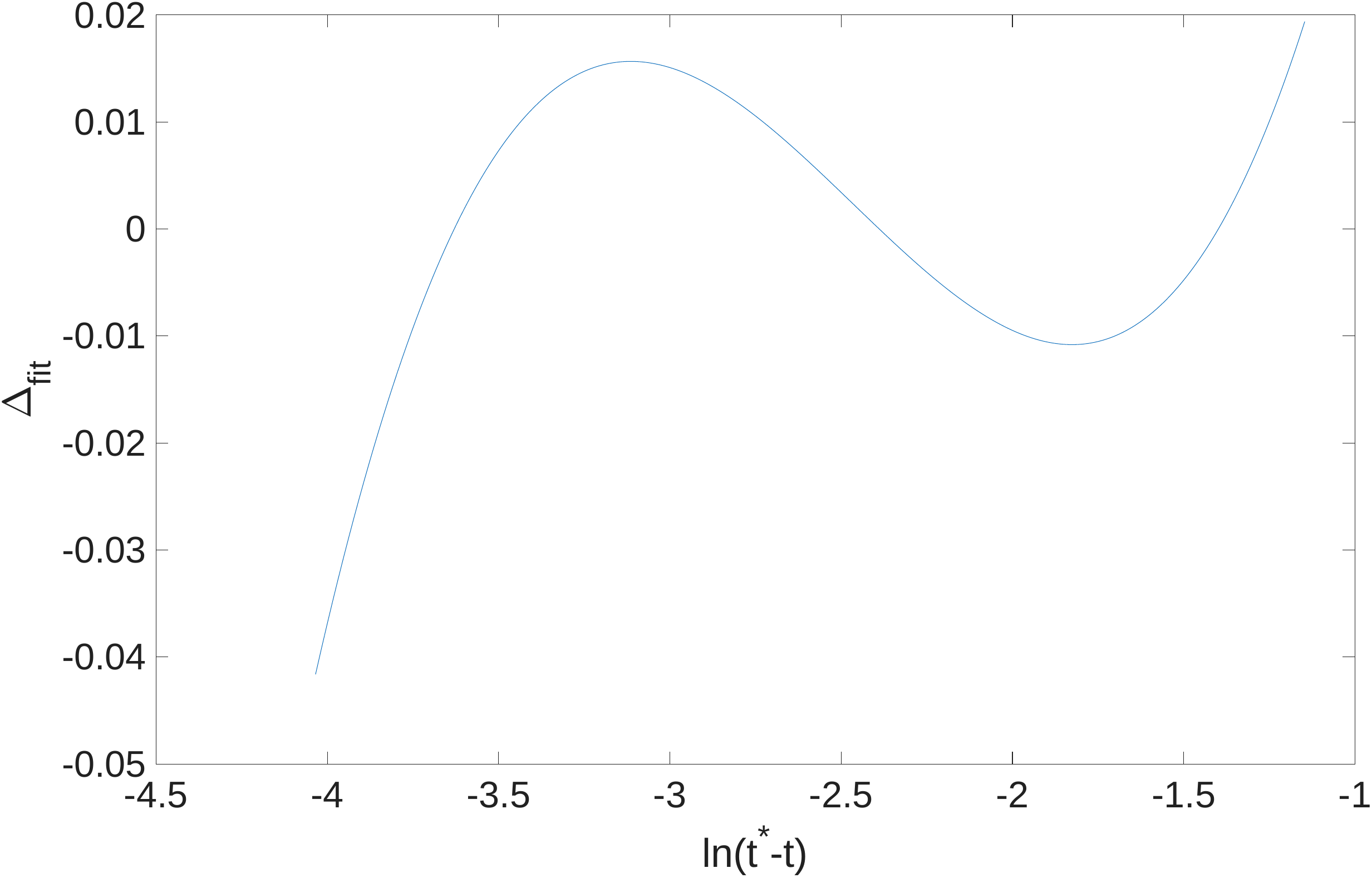}
    \end{subfigure}

    \caption{On top the logarithm of the $L^\infty$ norm of the last 
	$1000$ time steps of the solution to the system KBK 
	(\ref{Kaup-true-2D}) for the 
initial data (\ref{inistatic2}) with $\mu=-0.1$ and the fitted curve 
$\alpha \ln(t^* - t) + \beta$ on the left, and the difference $\Delta_{fit}$  (`fitting error') between both 
	curves on the right. Analogous plots on the bottom for the 
    $L^2$ norm}
    \label{fig:combined_norms_stat_mg}
\end{figure}

The relative fitting error stays below $1\%$ (see Fig.~\ref{fig:combined_norms_stat_mg}), and the fitted parameters are insensitive to changes in the fitting window and initial guess. The agreement between the exponents obtained from the $L^\infty$ and $L^2$ norms further indicates robustness. Overall, the measured values of $\alpha$ are consistent with the predicted exponent $-1$, corresponding to an exponential dependence of the scaling factor $L$ on $\tau$.

This behavior is reminiscent of the focusing regime in the 
supercritical NLS equation, where solutions above a critical threshold can undergo finite-time blow-up due to nonlinear self-focusing.

\section{Transverse stability of the line solitary wave}\label{line}
In this section we study the transverse stability of the 1D solitons 
(\ref{soliton}) that are exactly $y$-independent (and thus in 
$y$-direction infinitely extended) \emph{line solitary waves}. We 
consider initial data of the form
\begin{equation}
	\eta(x,y,0) = \eta_{C}(x)\pm \mu \exp(-x^{2}-y^{2}),\quad v_{x}= 
	v_{x,C}(x).
	\label{initrans}
\end{equation}
We work with $L_x = 20$, $L_y = 5$, $N_x = 2^{12}$, and $N_y = 
2^{10}$. 

\subsection{Line solitary waves not satisfying the non-cavitation condition}

First we consider the case $C = 0.8$ and $\mu = -0.1$. Time integration is performed using $N_t = 2 \times 10^4$ time steps over the interval $t \in [0, 18.4]$, and $N_t = 2 \times 10^3$ time steps over $t \in (18.4, 18.5]$.
The energy conservation used to check numerical accuracy drops below 
$10^{-3}$ at $t \approx 18.4695$. Thus we terminate the computation 
at this point due to a loss of numerical accuracy.

\begin{figure}[htb!]
    \centering
    \begin{subfigure}[b]{0.49\textwidth}
        \centering
        \includegraphics[width=\textwidth]{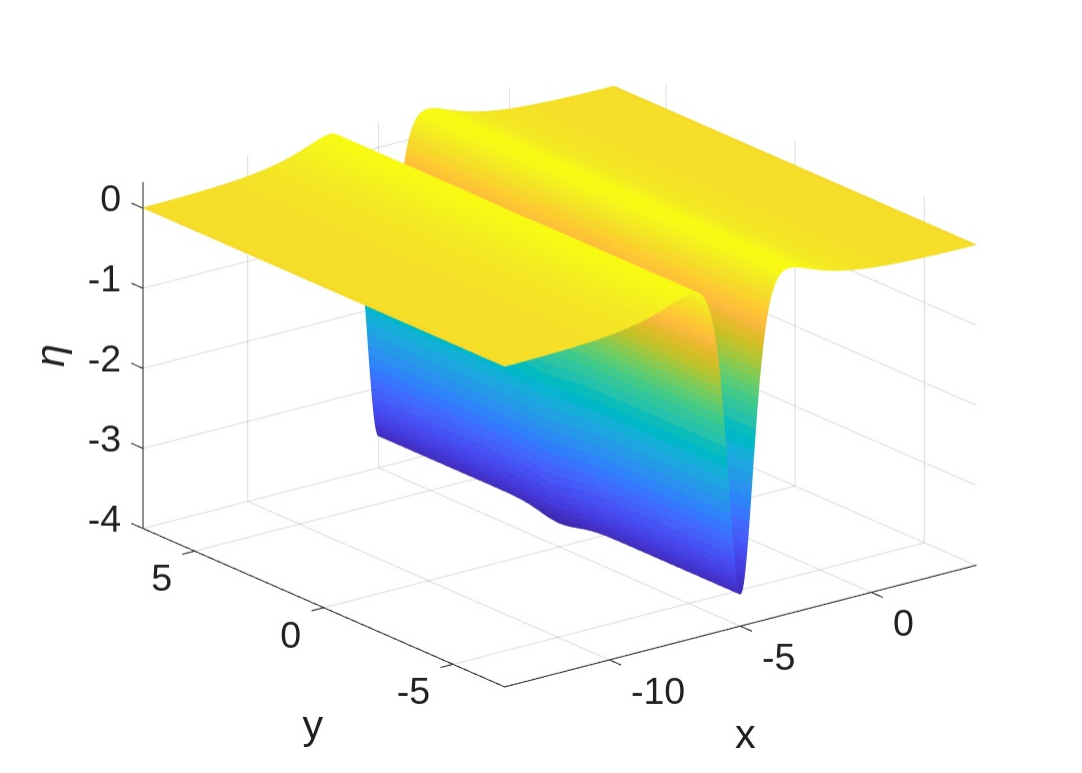}
    \end{subfigure}
    \hfill
    \begin{subfigure}[b]{0.49\textwidth}
        \centering
        \includegraphics[width=\textwidth]{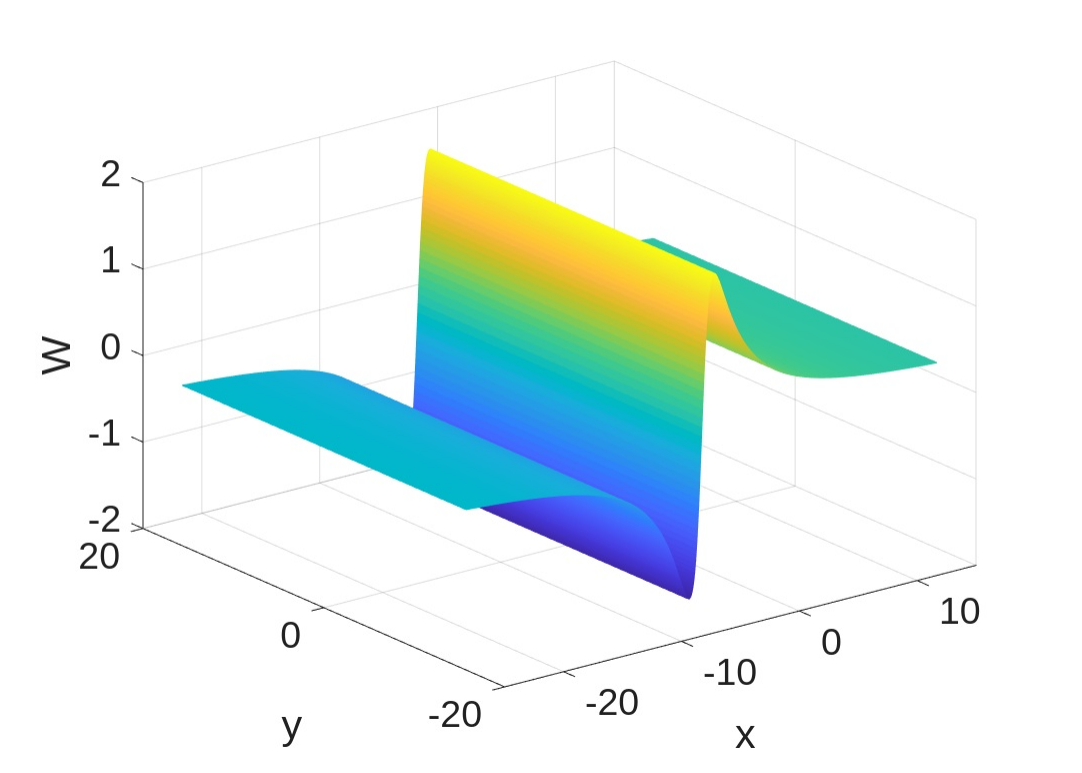}
    \end{subfigure}

    \begin{subfigure}[b]{0.49\textwidth}
        \centering
        \includegraphics[width=\textwidth]{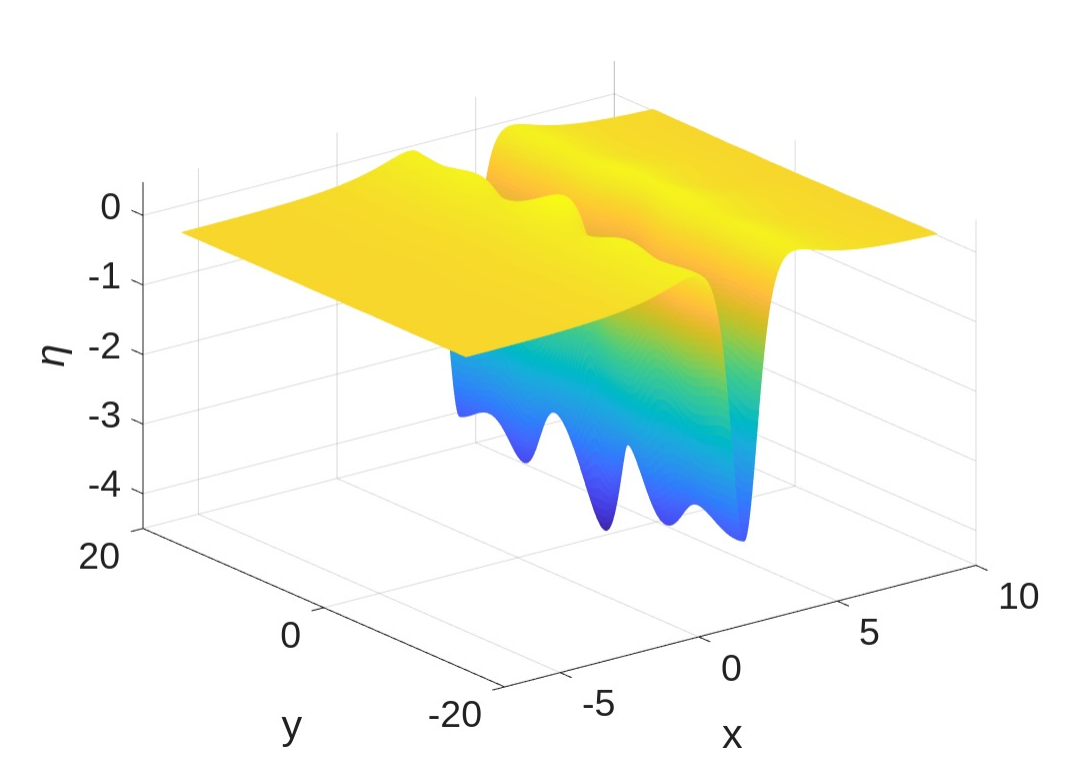}
    \end{subfigure}
    \hfill
    \begin{subfigure}[b]{0.49\textwidth}
        \centering
        \includegraphics[width=\textwidth]{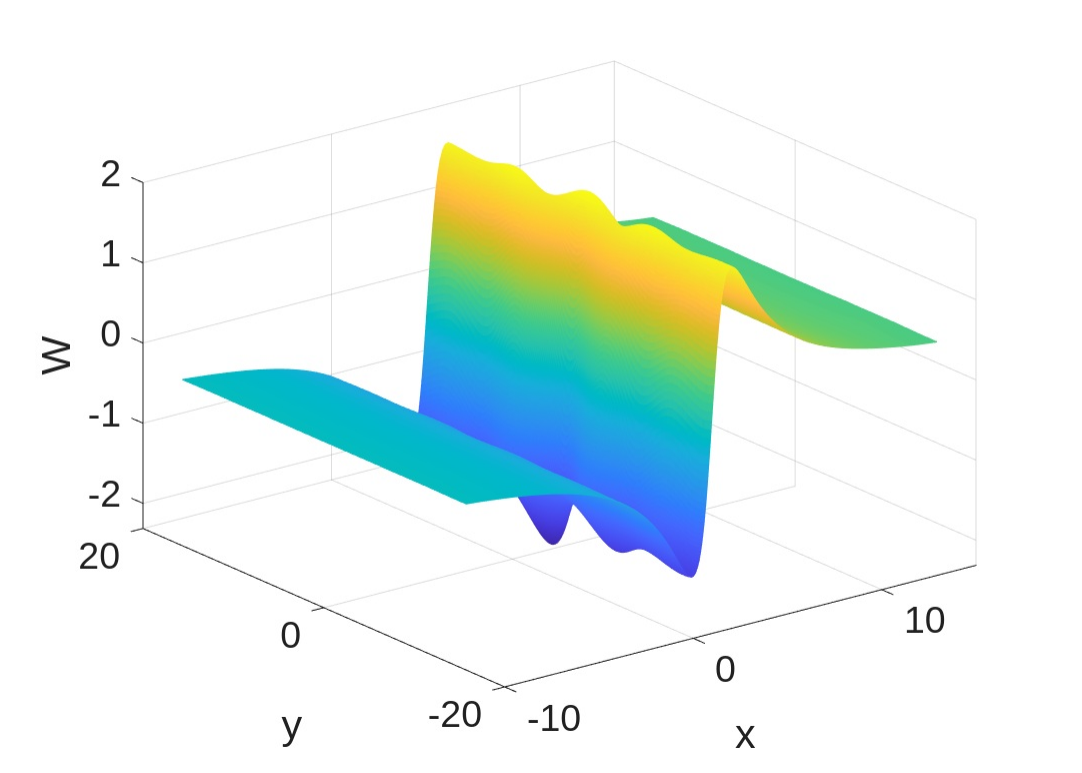}
    \end{subfigure}

    \begin{subfigure}[b]{0.49\textwidth}
        \centering
        \includegraphics[width=\textwidth]{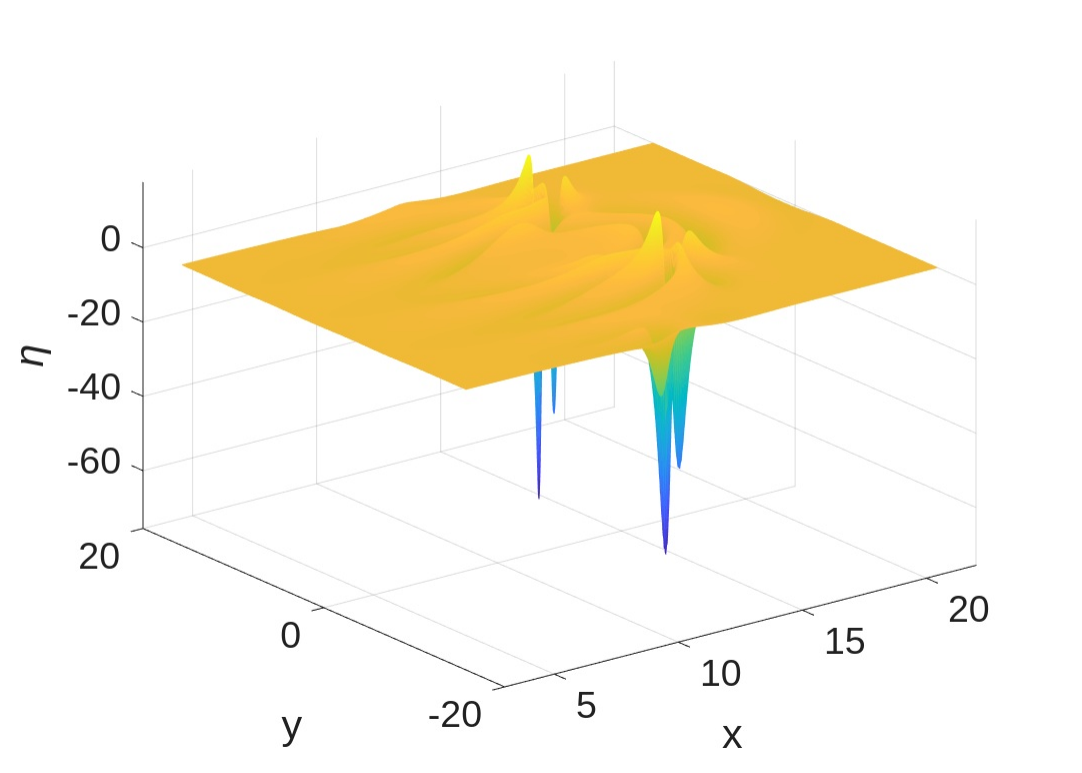}
    \end{subfigure}
    \hfill
    \begin{subfigure}[b]{0.49\textwidth}
        \centering
        \includegraphics[width=\textwidth]{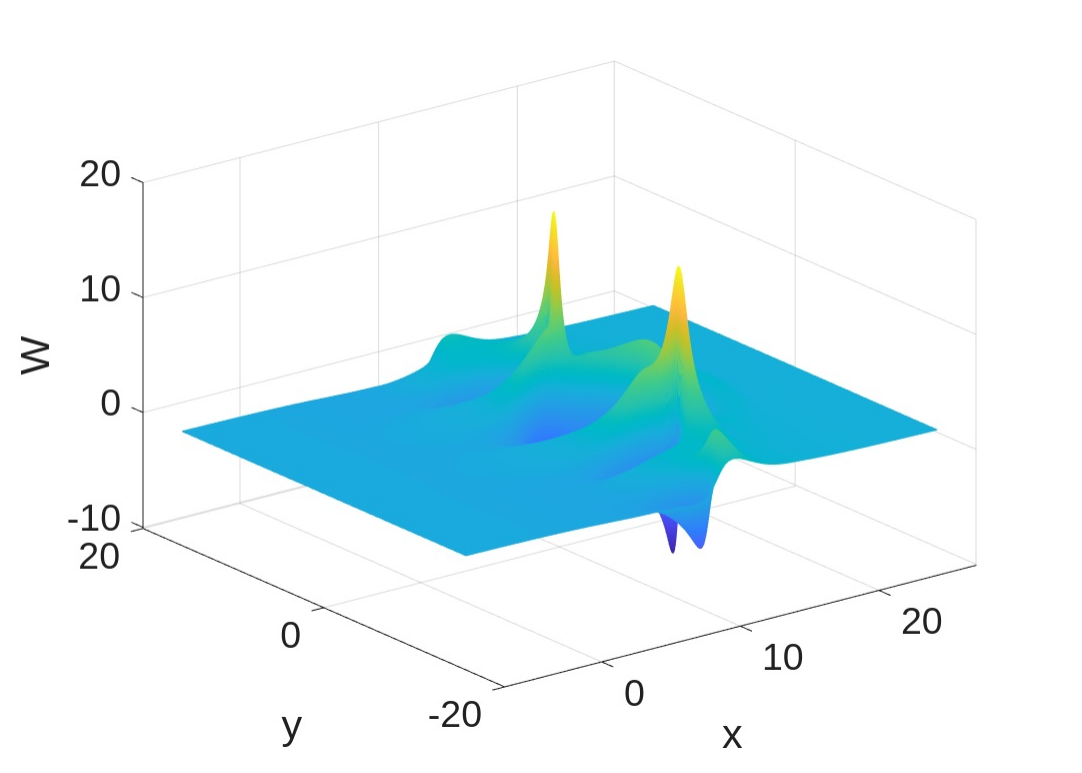}
    \end{subfigure}

\caption{Solution to the KBK system (\ref{Kaup-true-2D}) with initial data (\ref{initrans}) for $\mu = -0.1$ and $C=0.8$. Rows correspond to times $t=0$ (top), $t=10.028$ (middle), and $t=18.4$ (bottom). In each row, $\eta$ is shown on the left and $W$ on the right.}
    \label{fig:combined_mg}
\end{figure}

As can be seen in Fig.~\ref{fig:combined_mg}, under this initial 
condition two peaks emerge, becoming increasingly tall and narrow as 
the solution approaches the blow-up time, suggesting the formation of 
an $L^{\infty}$ blow-up. Once more the blow-up occurs for negative values of $1+\eta$ which should not appear in the context of water waves.

Note that the statment about blow-up in section 3.3 
has to be understood as local, i.e., in the case of two well 
separated blow-up profiles, the dynamical rescaling is to be applied 
locally to to each of them.

We recall the previously observed scaling near blow-up,
\[
\|\eta\|_\infty \sim (t^* - t)^\alpha, \qquad 
\|\eta\|_2 ^2 \sim (t^* - t)^\alpha,
\]
with $\alpha < 0$, and we extract $\alpha$ by fitting $\ln \|\eta\|_2$ and $\ln \|\eta\|_\infty$ to $\alpha \ln(t^* - t) + \beta$ over the last $1000$ time steps:
\begin{itemize}
    \item for $\|\eta\|_2 ^2$: $\alpha = -0.8223$, $\beta = 6.2734$, $t^* = 18.4725$,
    \item for $\|\eta\|_\infty$: $\alpha = -0.8689$, $\beta = 2.0275$, $t^* = 18.4737$.
\end{itemize}

\begin{figure}[htb!]
    \centering
    \begin{subfigure}[b]{0.49\textwidth}
        \centering
        \includegraphics[width=\textwidth]{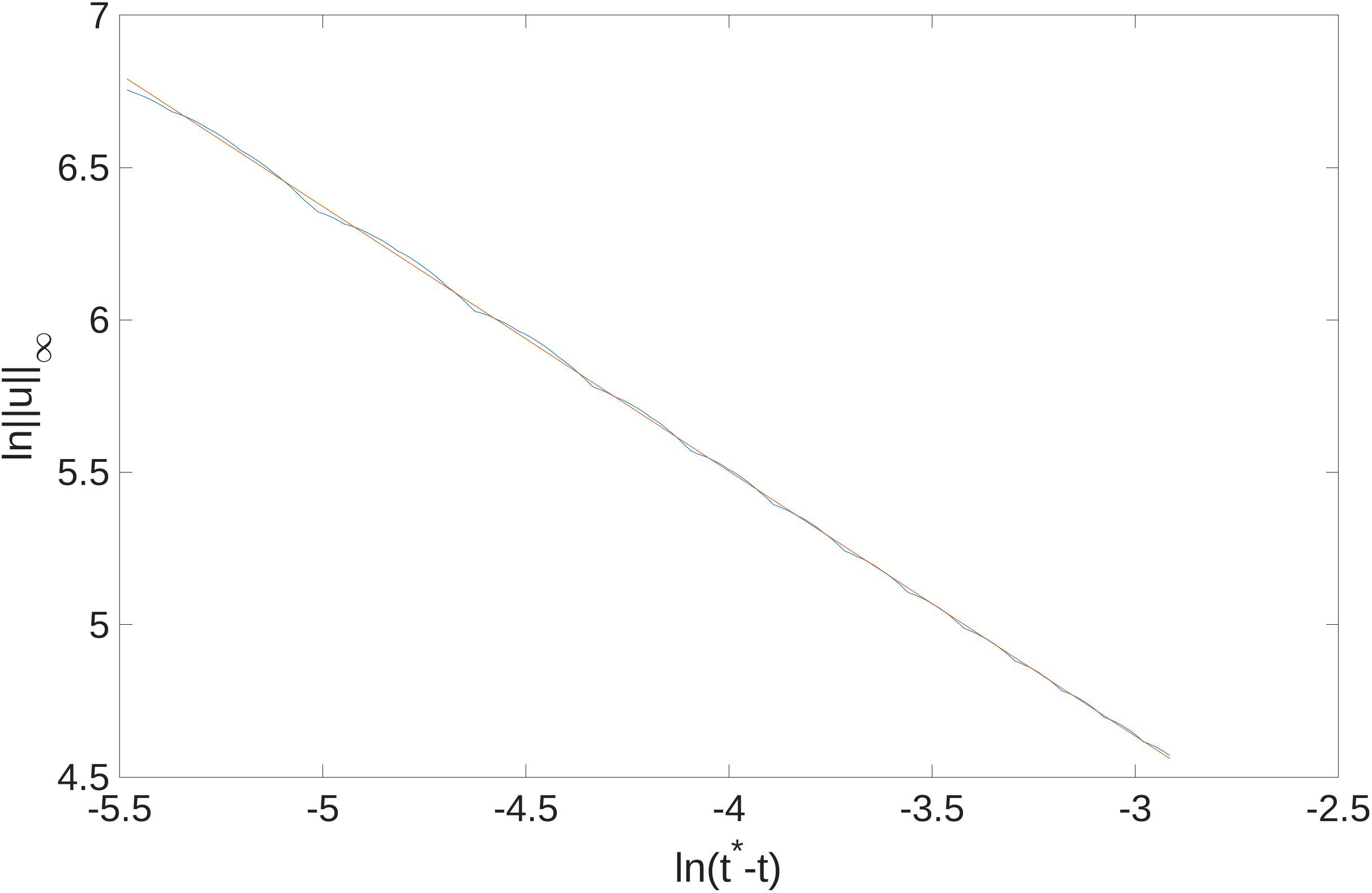}
    \end{subfigure}
    \hfill
    \begin{subfigure}[b]{0.5\textwidth}
        \centering
        \includegraphics[width=\textwidth]{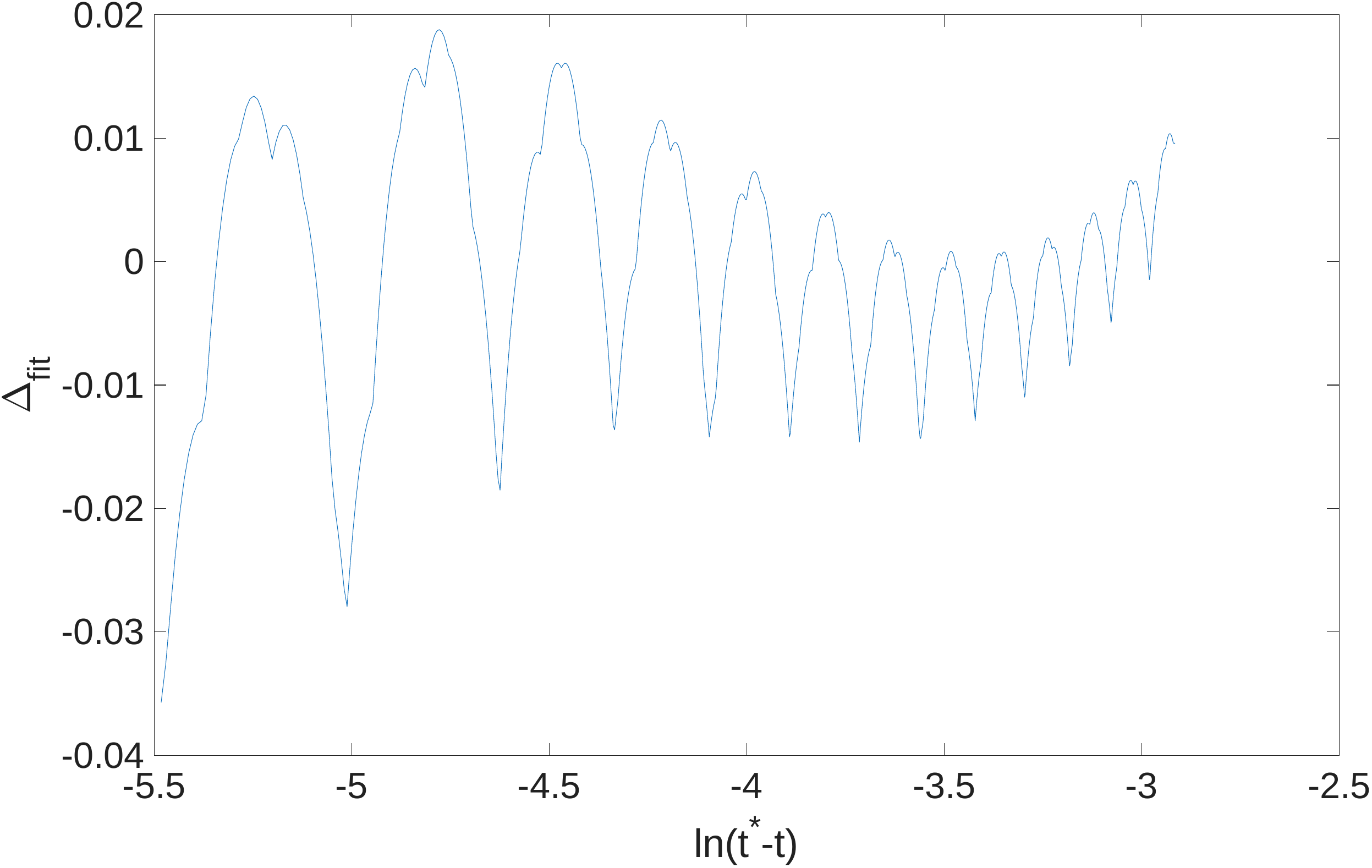}
    \end{subfigure}

    \begin{subfigure}[b]{0.49\textwidth}
        \centering
        \includegraphics[width=\textwidth]{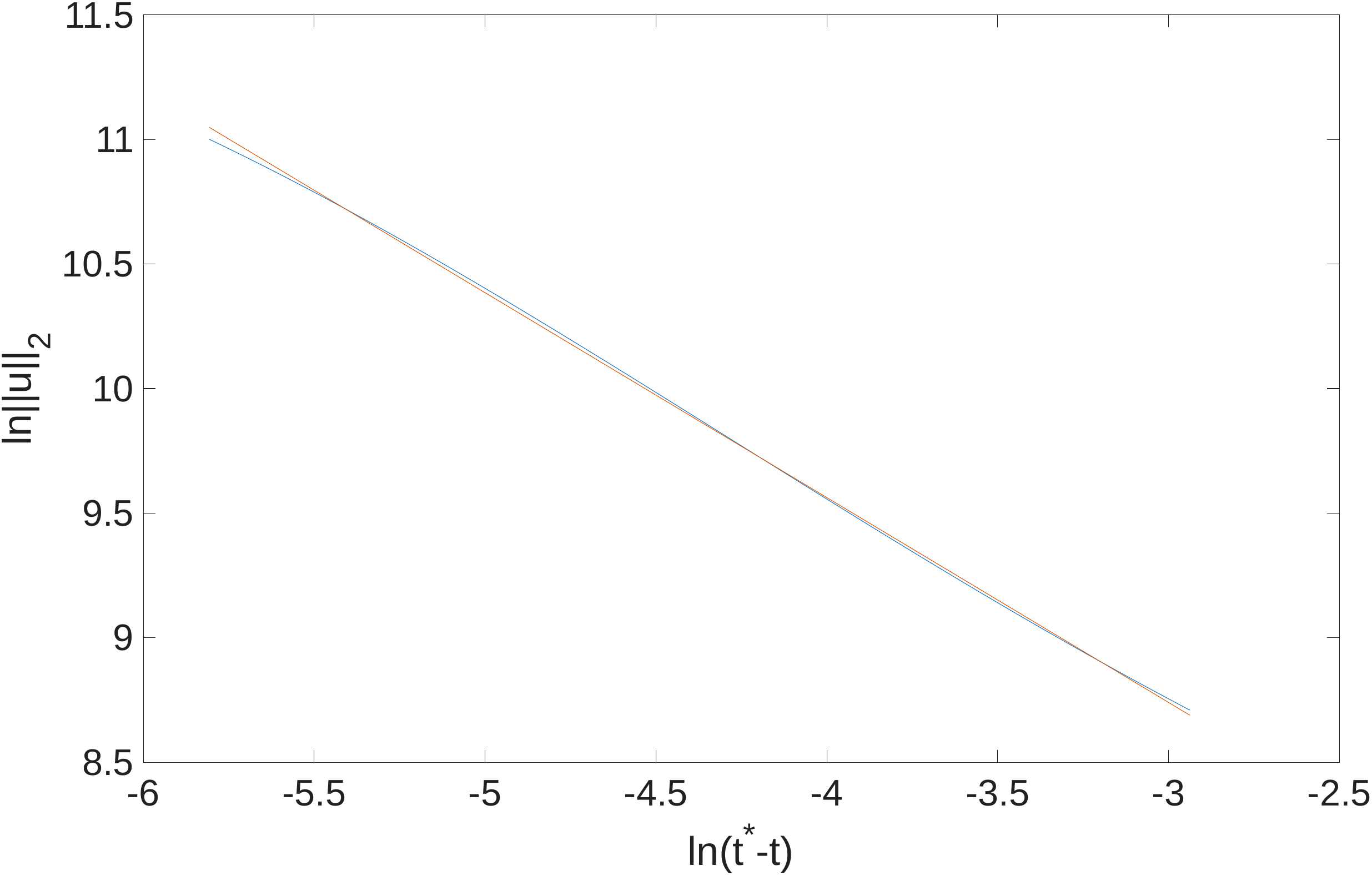}
    \end{subfigure}
    \hfill
    \begin{subfigure}[b]{0.5\textwidth}
        \centering
        \includegraphics[width=\textwidth]{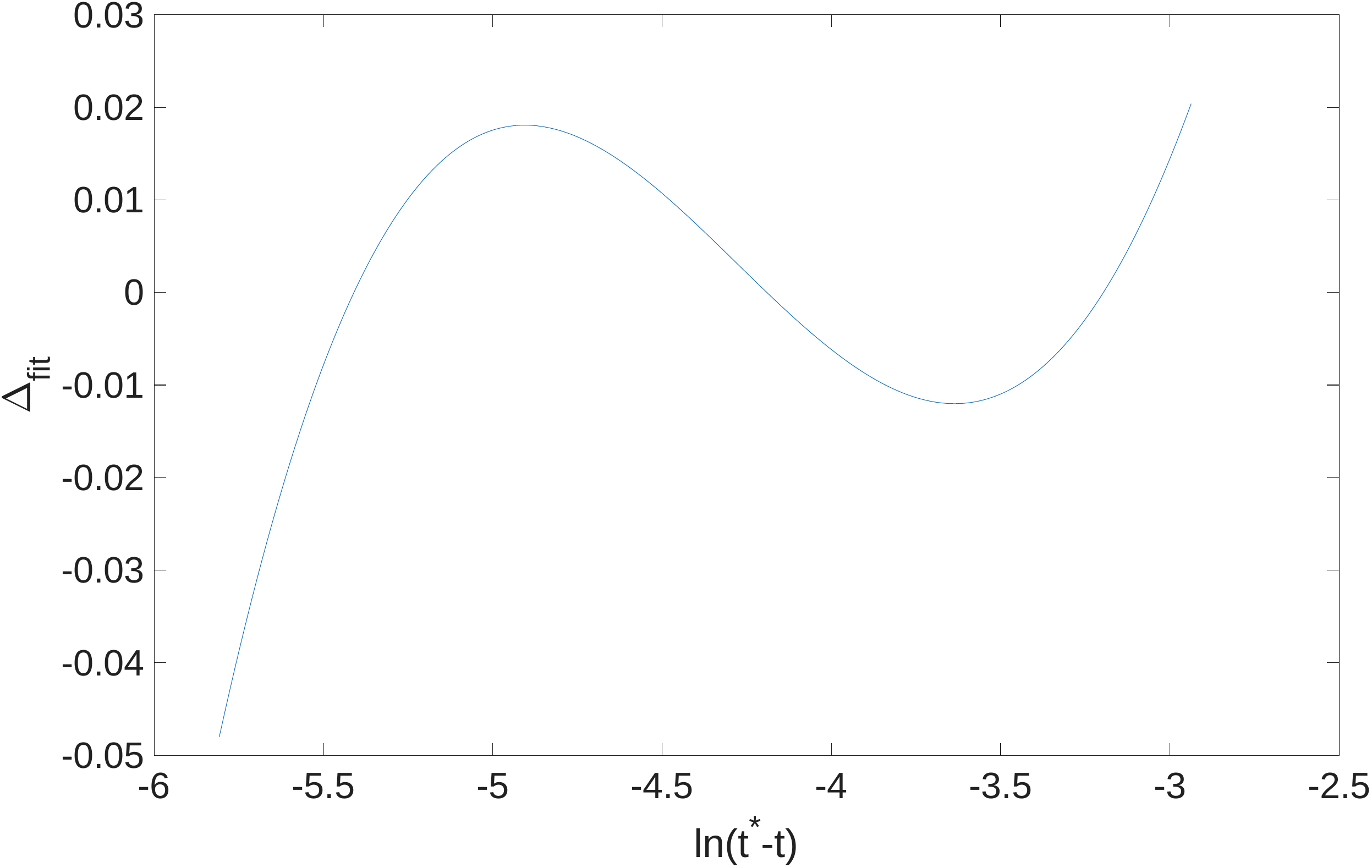}
    \end{subfigure}

    \caption{On top the logarithm of the $L^\infty$ norm of the last 
	$1000$ time steps of the solution to the system KBK 
	(\ref{Kaup-true-2D}) for the initial data(\ref{initrans}) with 
	$\mu = -0.1$, $C = 0.8$ and the fitted curve $\alpha \ln(t^* - t) + \beta$ 
	on the left, and the difference $\Delta_{fit}$ between both 
	curves on the right. Analogous plots on the bottom for the 
    $L^2$ norm}
    \label{fig:combined_norms}
\end{figure}

Relative fitting errors remain below $1\%$ (see Fig.~\ref{fig:combined_norms}), and the parameters are stable under changes in the fitting window and initial guess. The close agreement between the exponents from the $L^\infty$ and $L^2$ norms supports robustness, and the measured $\alpha$ is consistent with the predicted value $-1$.

\begin{figure}[htb!]
    \centering
    \begin{subfigure}[b]{0.49\textwidth}
        \centering
        \includegraphics[width=\textwidth]{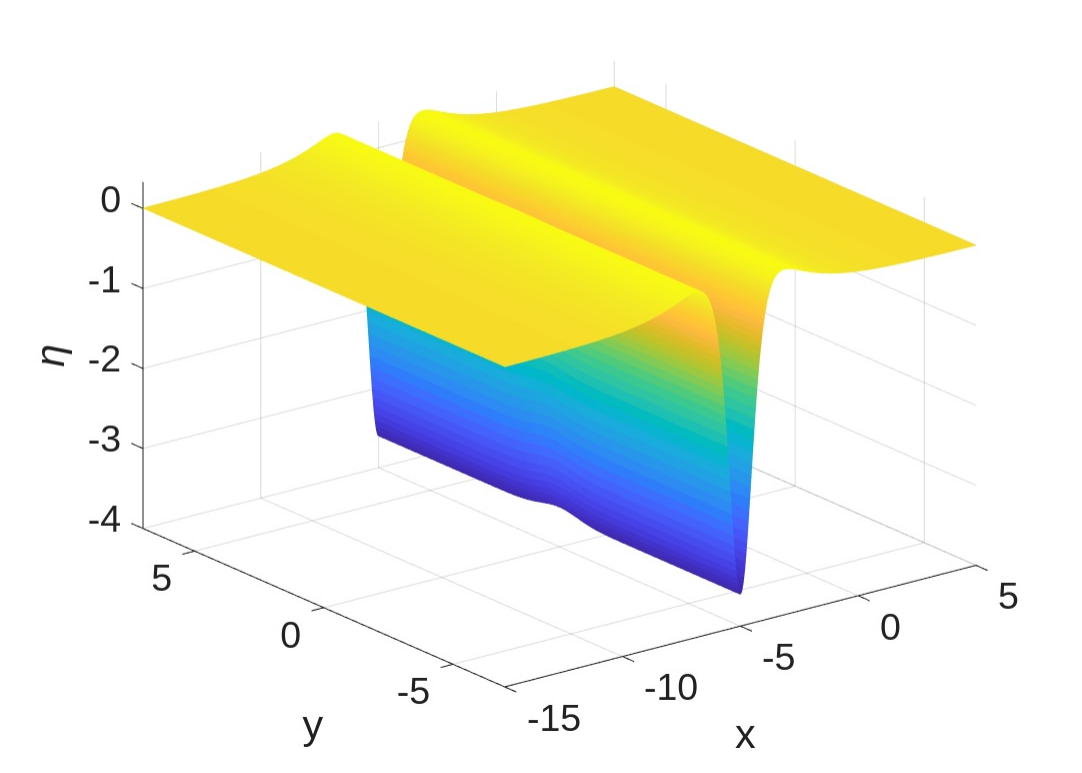}
    \end{subfigure}
    \hfill
    \begin{subfigure}[b]{0.49\textwidth}
        \centering
        \includegraphics[width=\textwidth]{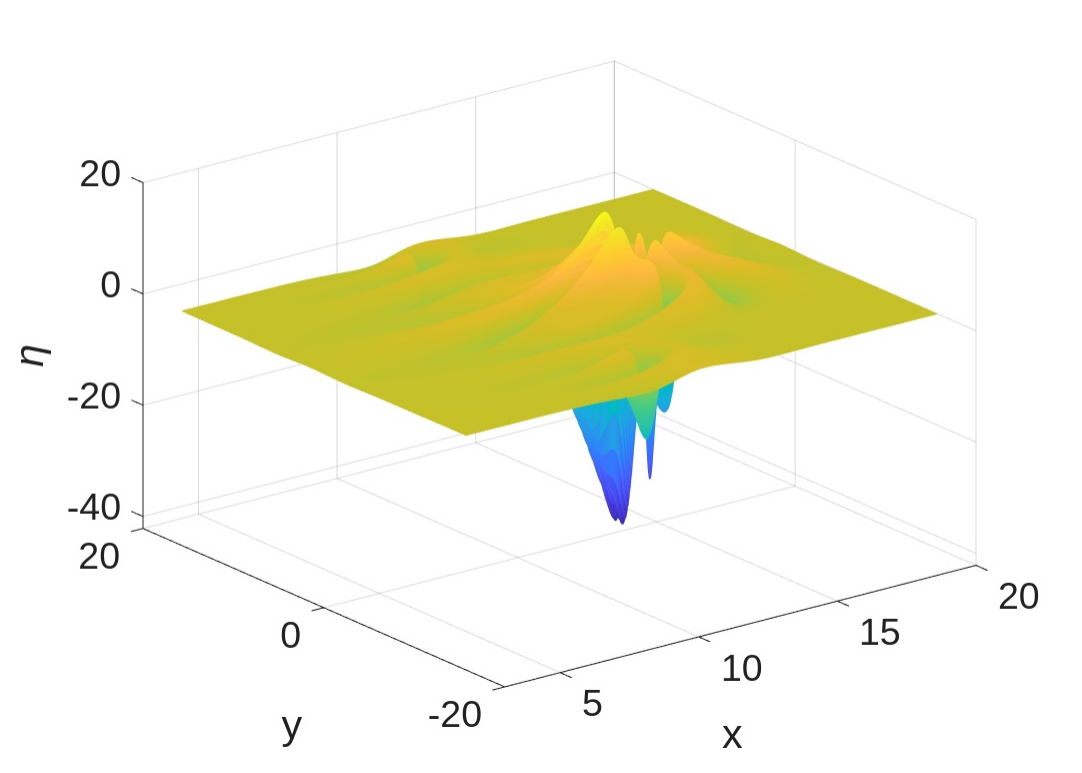}
    \end{subfigure}

    \begin{subfigure}[b]{0.49\textwidth}
        \centering
        \includegraphics[width=\textwidth]{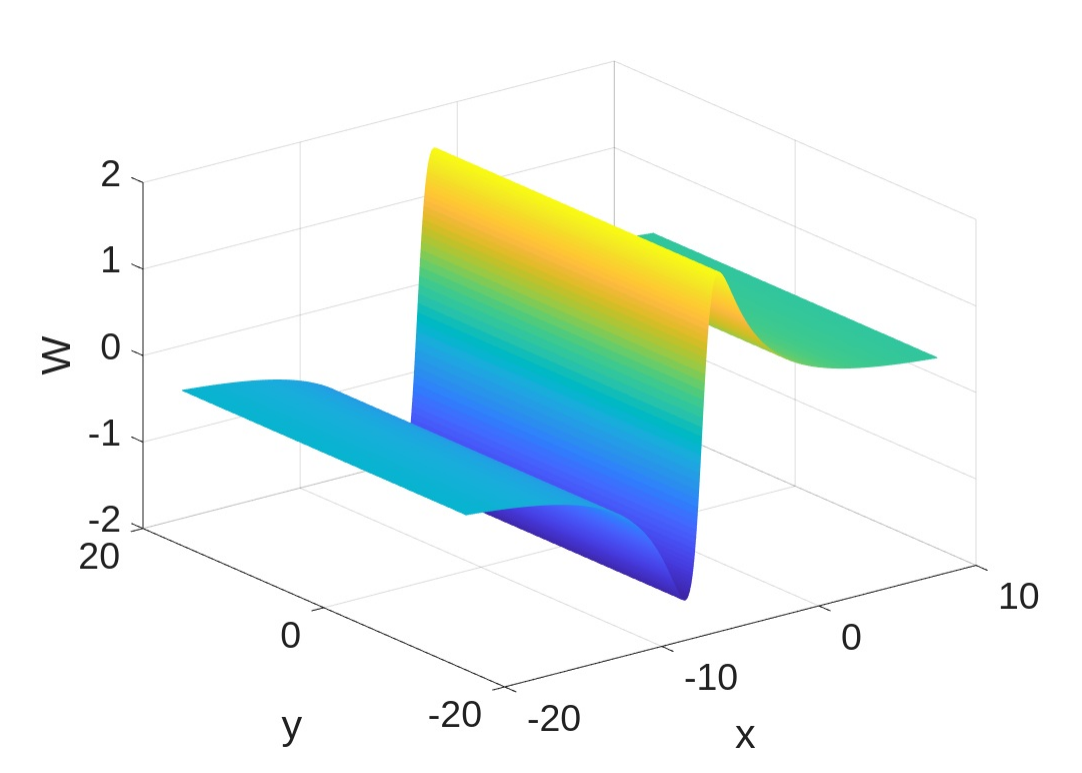}
    \end{subfigure}
    \hfill
    \begin{subfigure}[b]{0.49\textwidth}
        \centering
        \includegraphics[width=\textwidth]{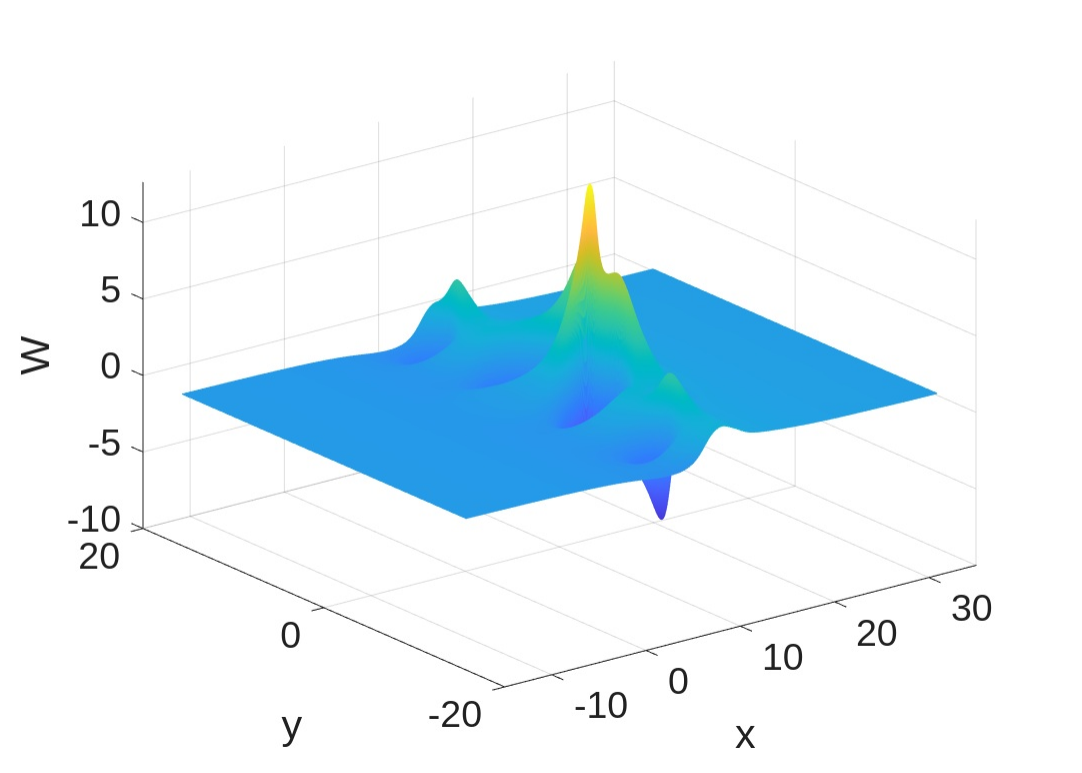}
    \end{subfigure}
\caption{Solution of the KBK system (\ref{Kaup-true-2D}) with initial data (\ref{initrans}), for $\mu = +0.1$ and $C = 0.8$. Top row: $\eta$ at $t=0$ (left) and $t=18.1$ (right). Bottom row: $W$ at $t=0$ (left) and $t=18.1$ (right).}
    \label{fig:combined_transp}
\end{figure}

We now examine the initial data (\ref{initrans}) for $C = 0.8$ and 
$\mu = +0.1$, for which a similar behavior is observed. Time 
integration is done over the interval $t \in [0, 18.1]$ using $N_t = 
2 \times 10^4$ time steps, and over $t \in (18.1, 18.2]$ using $N_t = 
2 \times 10^3$ time steps. The energy conservation falls below 
$10^{-3}$ at approximately $t \approx 18.1995$. Fig.~\ref{fig:combined_transp} shows two peaks that sharpen and grow as blow-up is approached, indicating an $L^\infty$ blow-up. 

The fitting of $\ln \|\eta\|_2$ and $\ln \|\eta\|_\infty$ to $\alpha 
\ln(t^* - t) + \beta$ over the last $1000$ time steps, gives:
\begin{itemize}
    \item for $\|\eta\|_2 ^2$: $\alpha = -0.5401$, $\beta = 6.9425$, $t^* = 18.2005$,
    \item for $\|\eta\|_\infty$: $\alpha = -0.9363$, $\beta = 1.5591$, $t^* = 18.2014$.
\end{itemize}

\begin{figure}[htb!]
    \centering
    \begin{subfigure}[b]{0.49\textwidth}
        \centering
        \includegraphics[width=\textwidth]{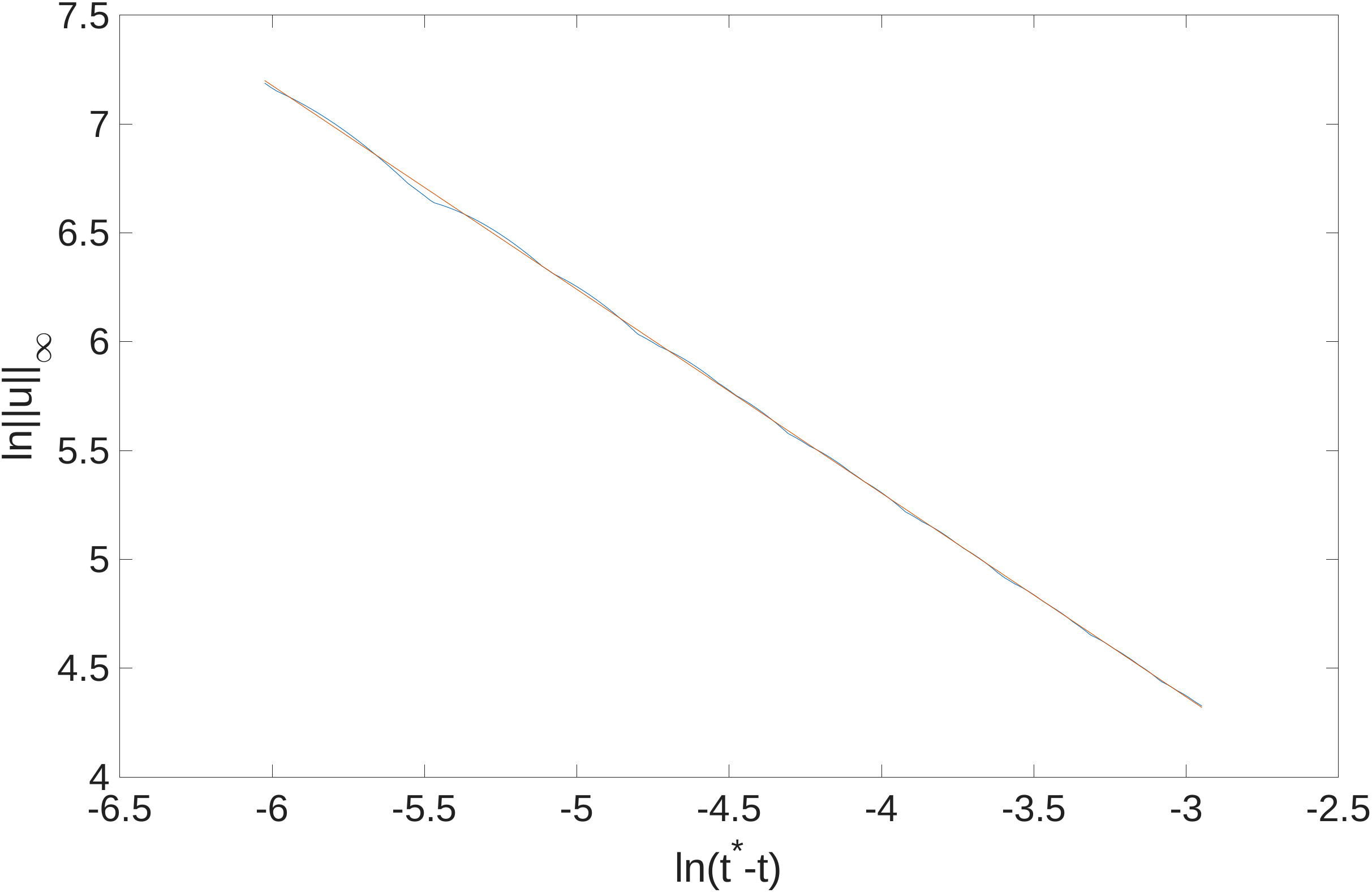}
    \end{subfigure}
    \hfill
    \begin{subfigure}[b]{0.5\textwidth}
        \centering
        \includegraphics[width=\textwidth]{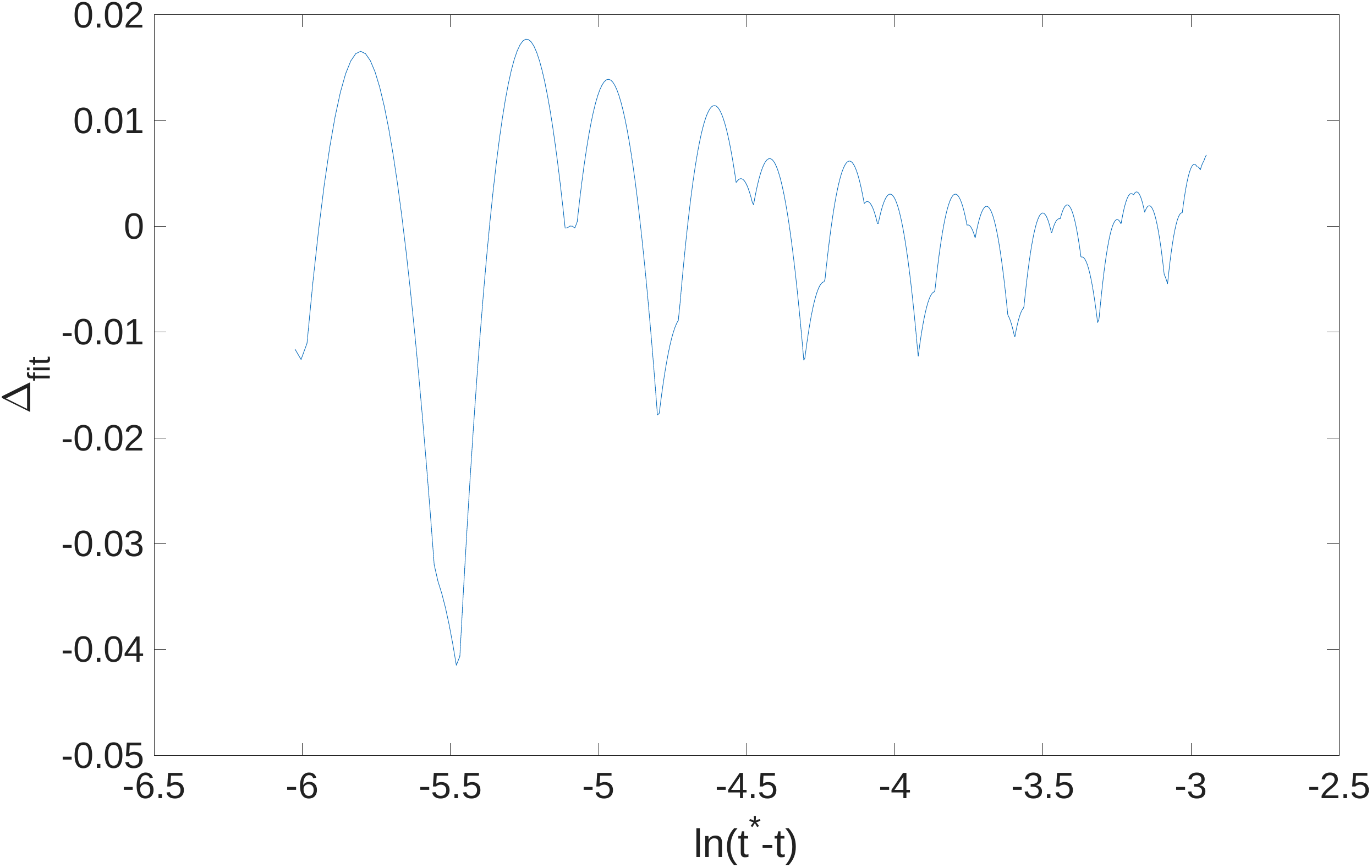}
    \end{subfigure}

    \begin{subfigure}[b]{0.49\textwidth}
        \centering
        \includegraphics[width=\textwidth]{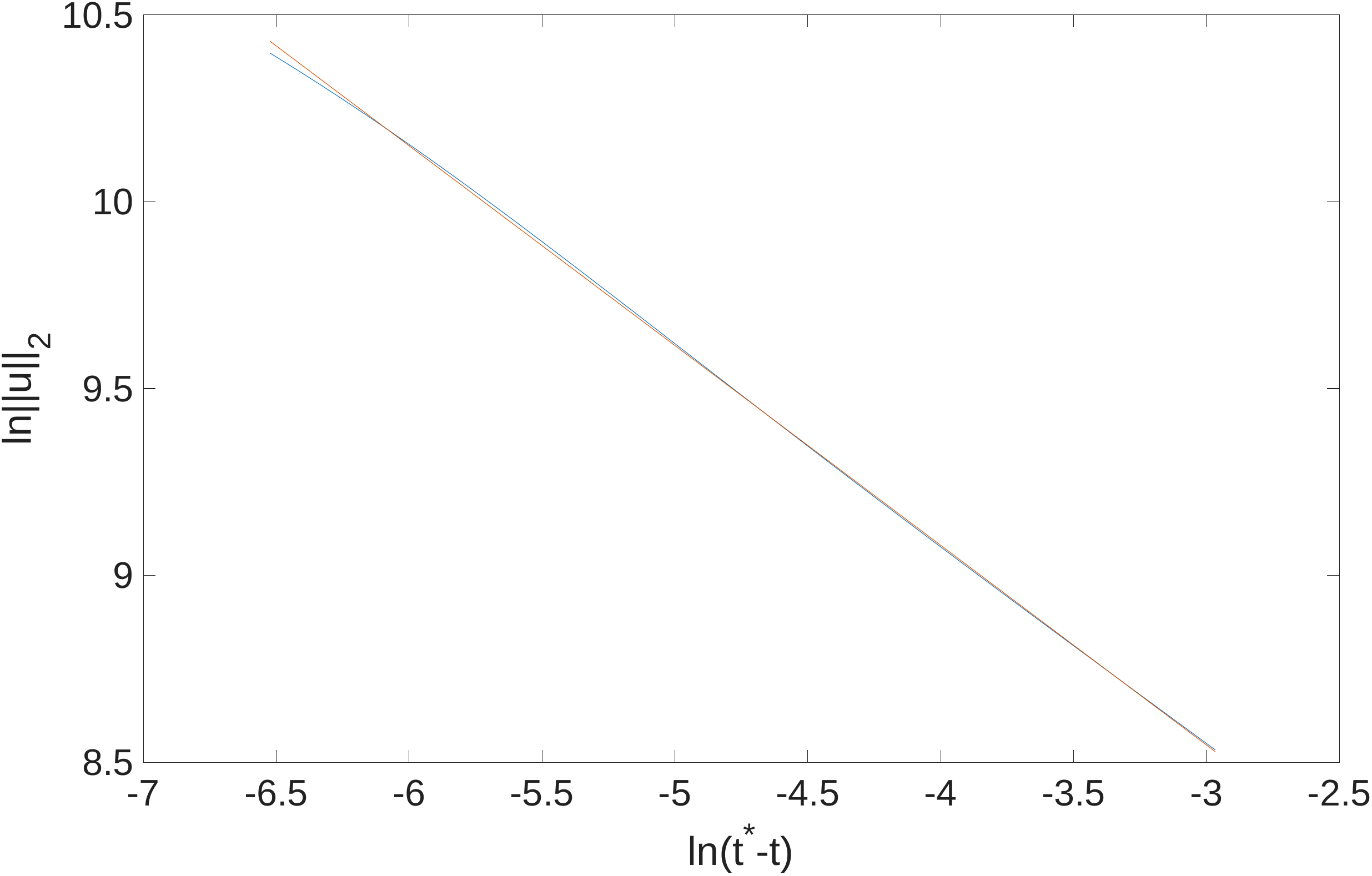}
    \end{subfigure}
    \hfill
    \begin{subfigure}[b]{0.5\textwidth}
        \centering
        \includegraphics[width=\textwidth]{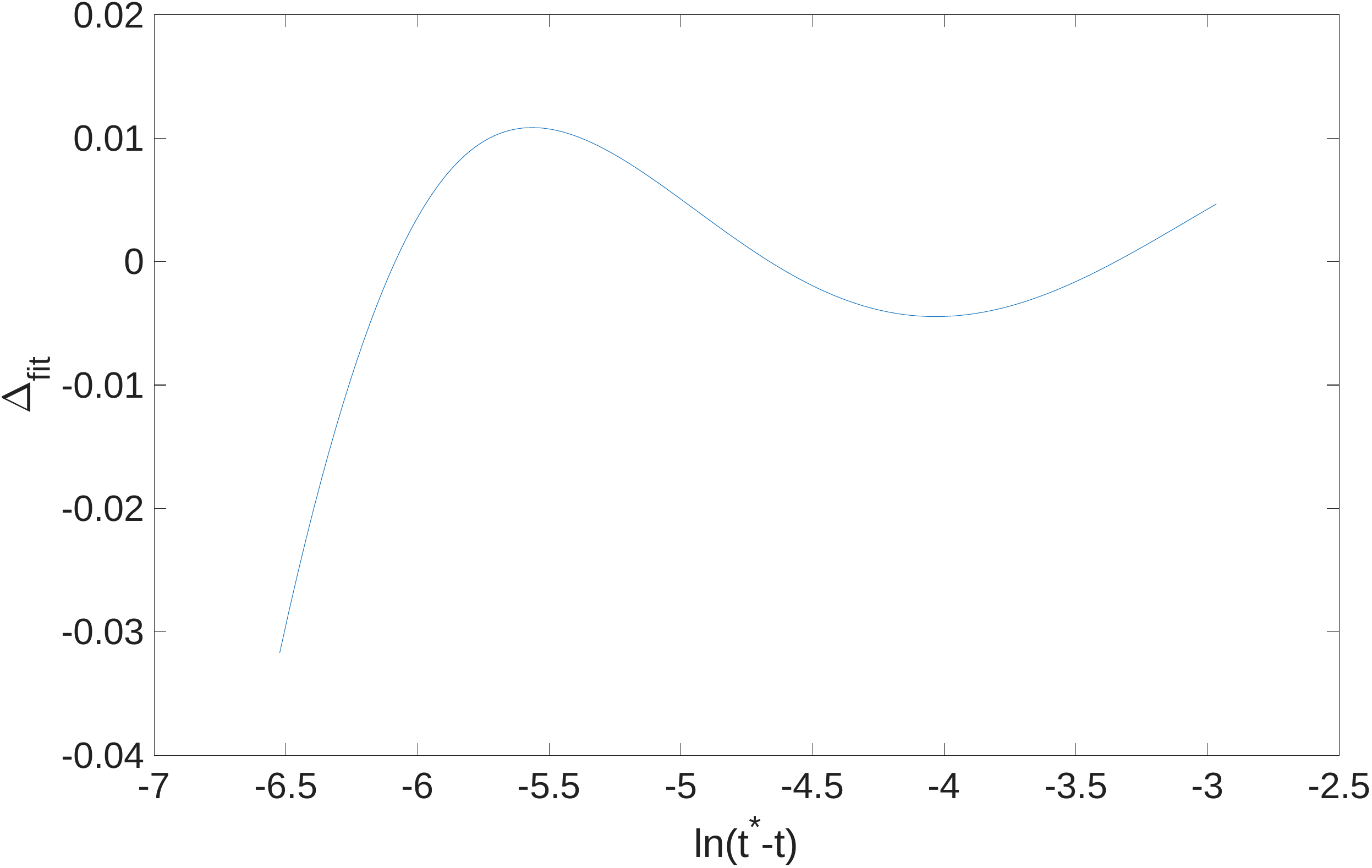}
    \end{subfigure}

    \caption{On top the logarithm of the $L^\infty$ norm of the last 
	$1000$ time steps of the solution to the system KBK 
	(\ref{Kaup-true-2D}) for the initial data(\ref{initrans}) with 
	$\mu = +0.1$, $C = 0.8$ and the fitted curve $\alpha \ln(t^* - t) + \beta$ 
	on the left, and the difference $\Delta_{fit}$ between both 
	curves on the right. Analogous plots on the bottom for the 
    $L^2$ norm}
    \label{fig:combined_norms2}
\end{figure}

For the $L^2$-based quantity, the fitted exponent $\alpha = -0.5401$ 
deviates significantly from the anticipated value $\alpha = -1$. 
Despite the discrepancy in the $L^2$ exponent, it is important to 
note that both fits are numerically stable, with relative errors 
smaller than $1\%$ (see Fig.~\ref{fig:combined_norms2}) and very 
close estimates of the blow-up time. Overall, if one were to infer 
the blow-up rate from these results, the evidence favors a scaling 
closer to the exponential dependence of the scaling factor $L$ on 
the rescaled time $\tau$ ($\alpha = -1$) rather than the alternative 
of an algebraic dependence  ($\alpha <-1$).

\begin{figure}[htb!]
    \centering

    \begin{subfigure}[b]{0.49\textwidth}
        \centering
        \includegraphics[width=\textwidth]{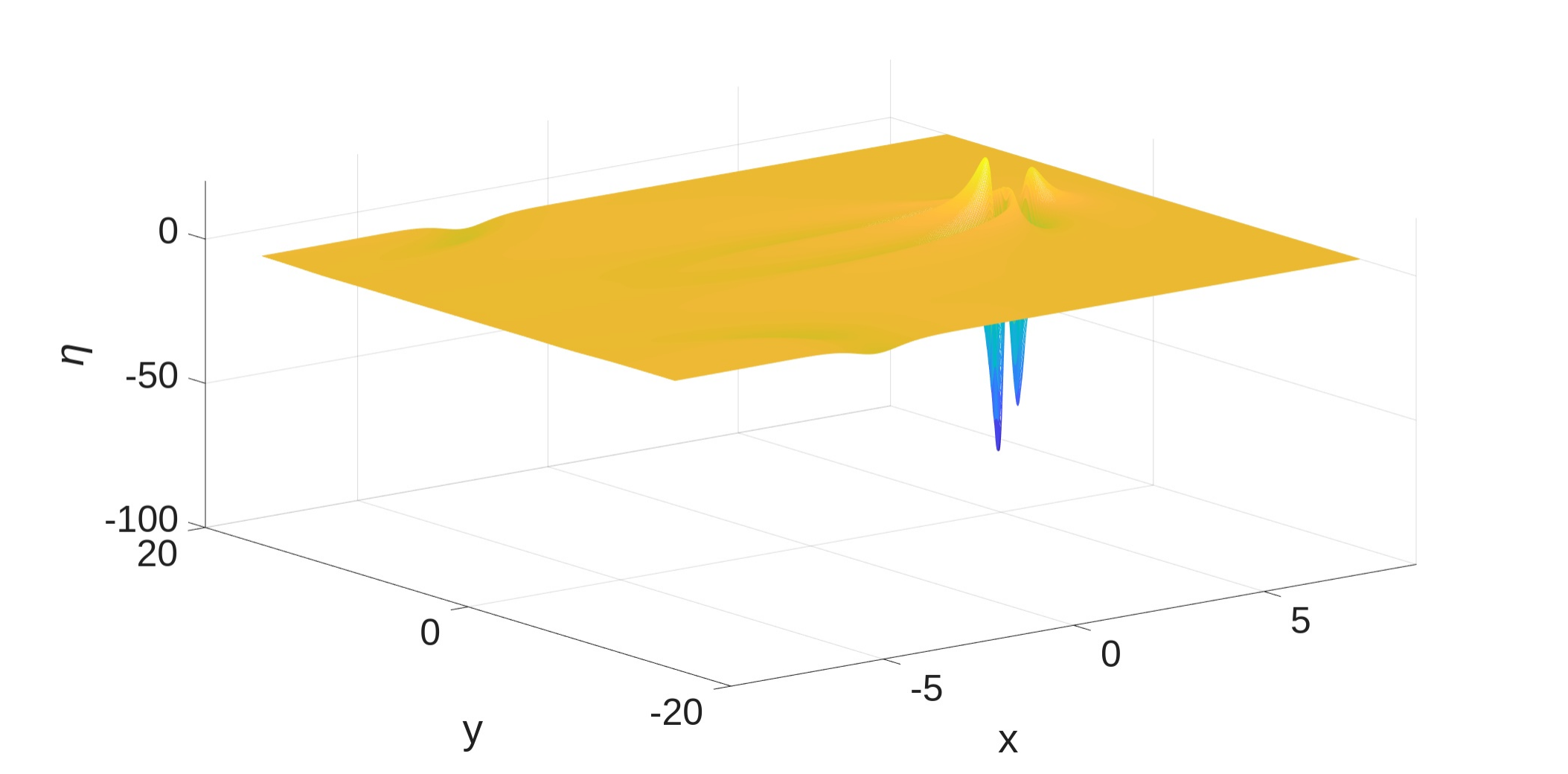}
    \end{subfigure}
    \hfill
    \begin{subfigure}[b]{0.49\textwidth}
        \centering
        \includegraphics[width=\textwidth]{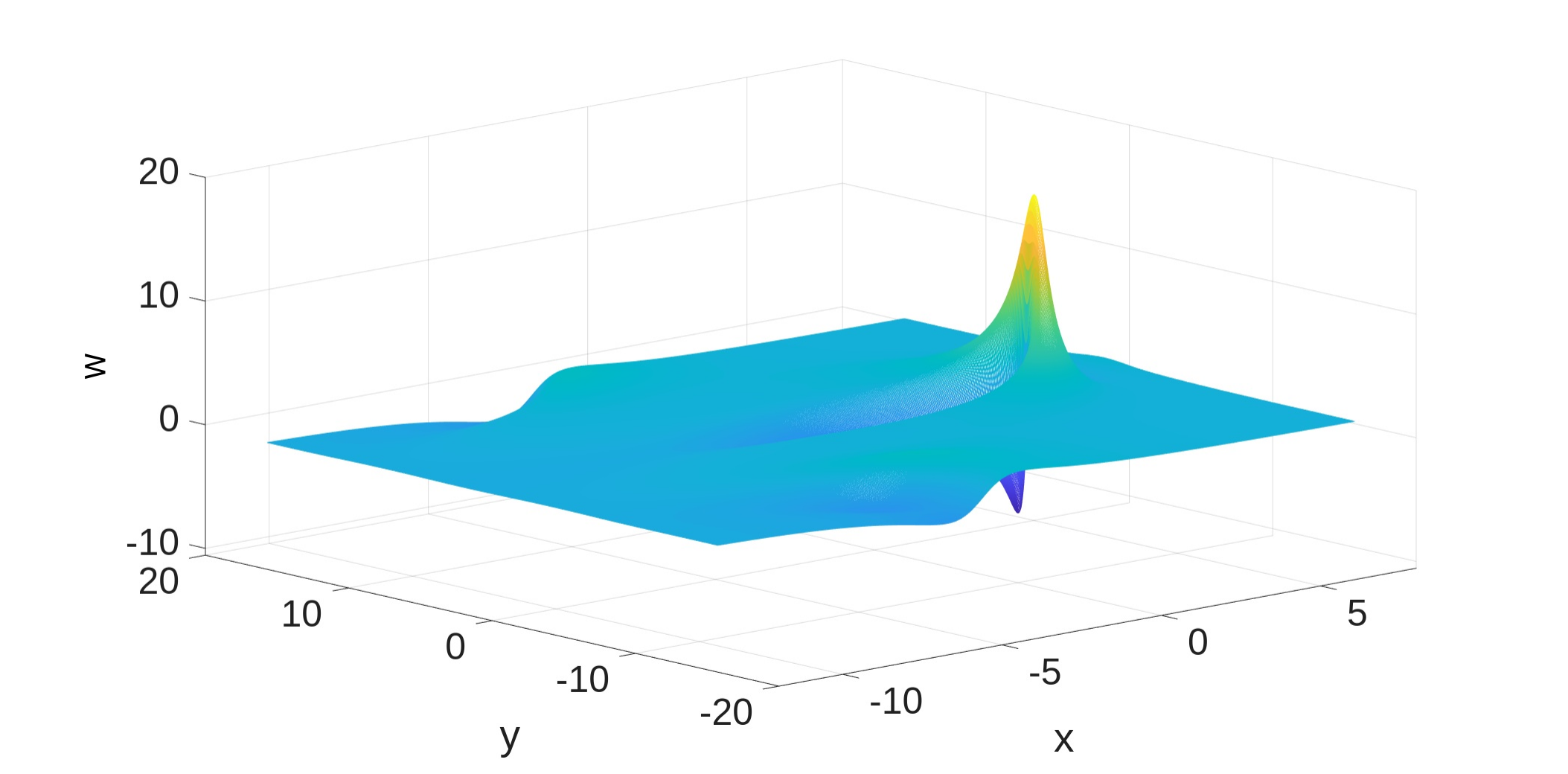}
    \end{subfigure}
\caption{Solution of the KBK system (\ref{Kaup-true-2D}) with initial data (\ref{initrans}), for $\mu = -0.1$ and $C = 0$ for $t=35.2$, on the left $\eta$, on the right $W$}.
    \label{fig:gaussc0}
\end{figure}

Note that in these cases $C>-0.5$, and thus, as 
pointed out in Section 2, the wave profile $\eta$ of the solitary 
wave solutions (\ref{soliton}) initially violates the non-cavitation 
condition $\eta+1>0$. Simulations with both larger and smaller 
perturbations, using different velocities $C>-0.5$, 
have also been conducted and exhibit the same type of behavior 
(an example with $C=0$ is given in 
Fig.~\ref{fig:gaussc0}). Notably, unlike the $L^{2}$-critical NLS 
equation, even when the $L^{2}$ norm of $\eta$ for the perturbed line 
soliton is smaller than that of the previously constructed stationary 
solution, a blow-up with a similar structure to the one described in 
this section still occurs. We do not show these 
examples here, as they are qualitatvely similar to the  previously 
presented examples.

\subsection{Line solitary waves satisfying the non-cavitation condition}

\begin{figure}[htb!]
    \centering
    \begin{subfigure}[b]{0.49\textwidth}
        \centering
        \includegraphics[width=\textwidth]{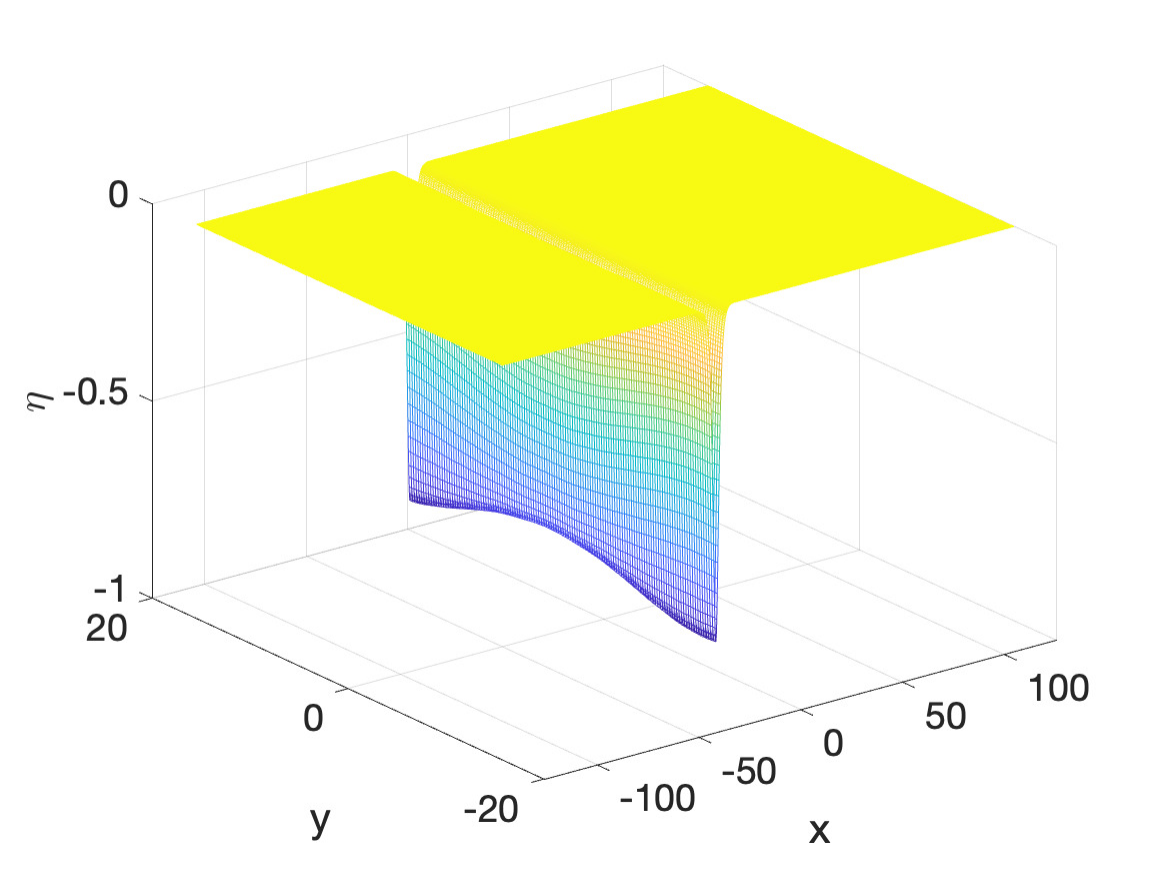}
    \end{subfigure}
    \hfill
    \begin{subfigure}[b]{0.49\textwidth}
        \centering
        \includegraphics[width=\textwidth]{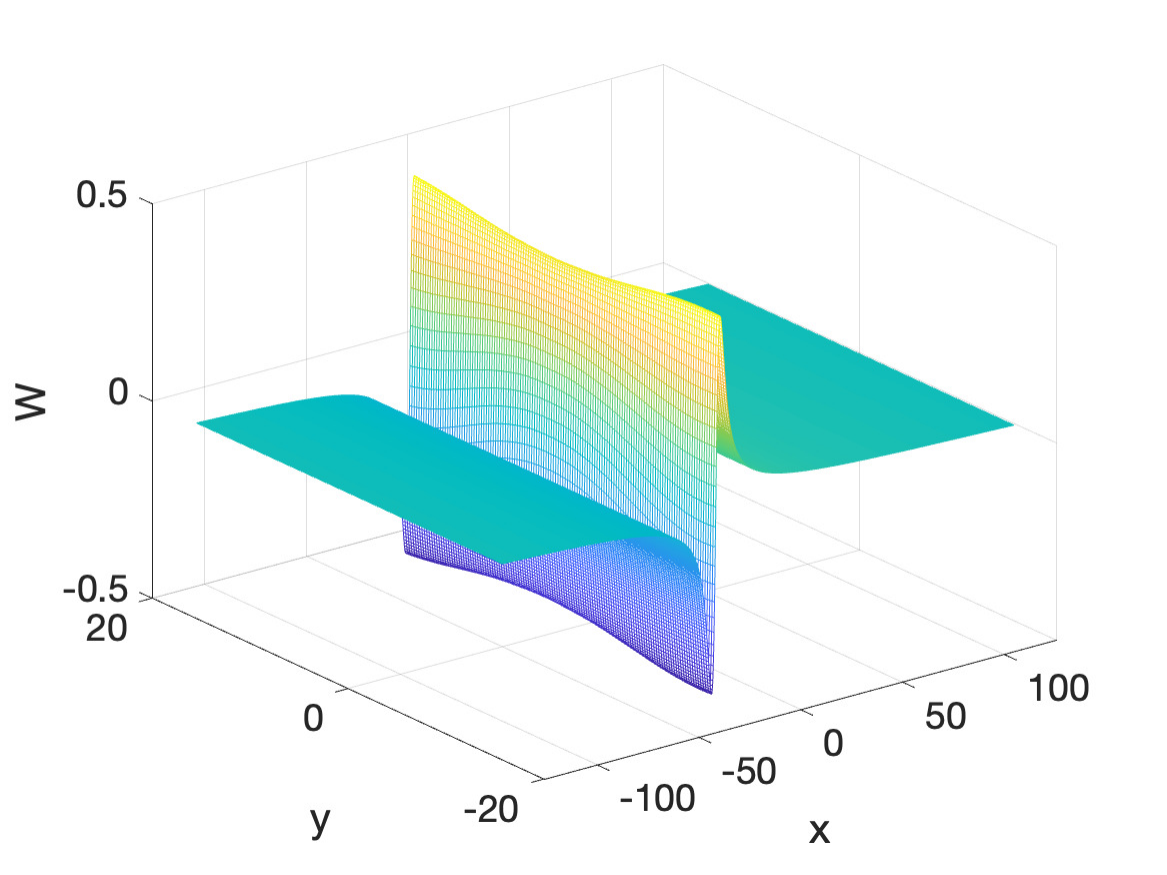}
    \end{subfigure}

    \begin{subfigure}[b]{0.49\textwidth}
        \centering
        \includegraphics[width=\textwidth]{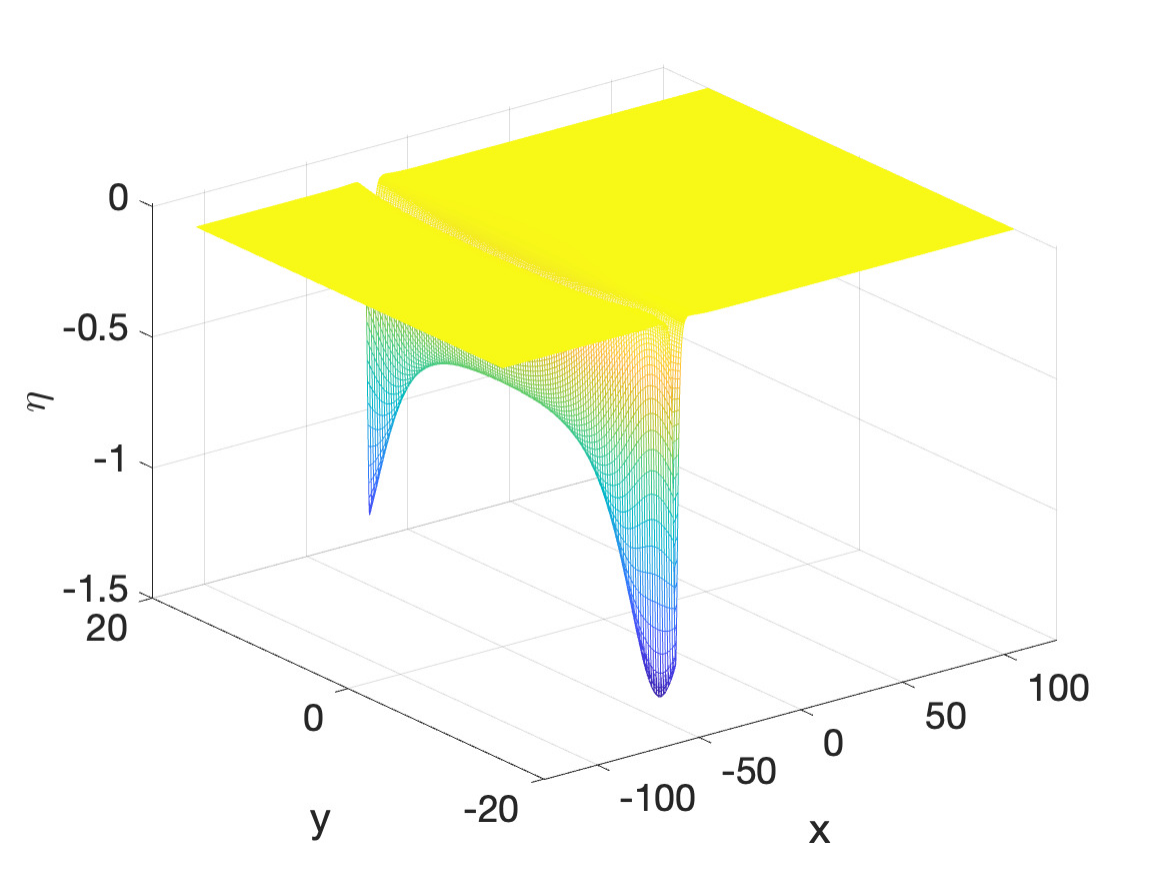}
    \end{subfigure}
    \hfill
    \begin{subfigure}[b]{0.49\textwidth}
        \centering
        \includegraphics[width=\textwidth]{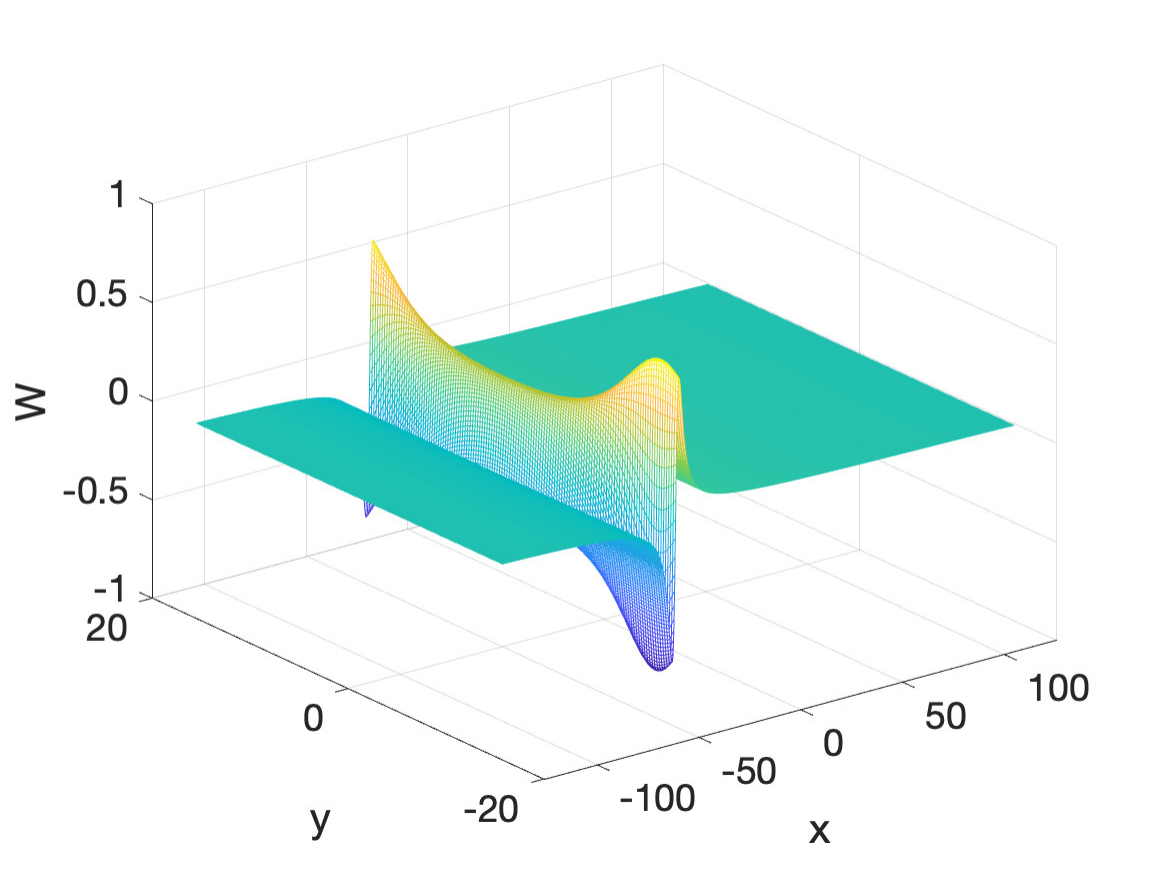}
    \end{subfigure}

    \begin{subfigure}[b]{0.49\textwidth}
        \centering
        \includegraphics[width=\textwidth]{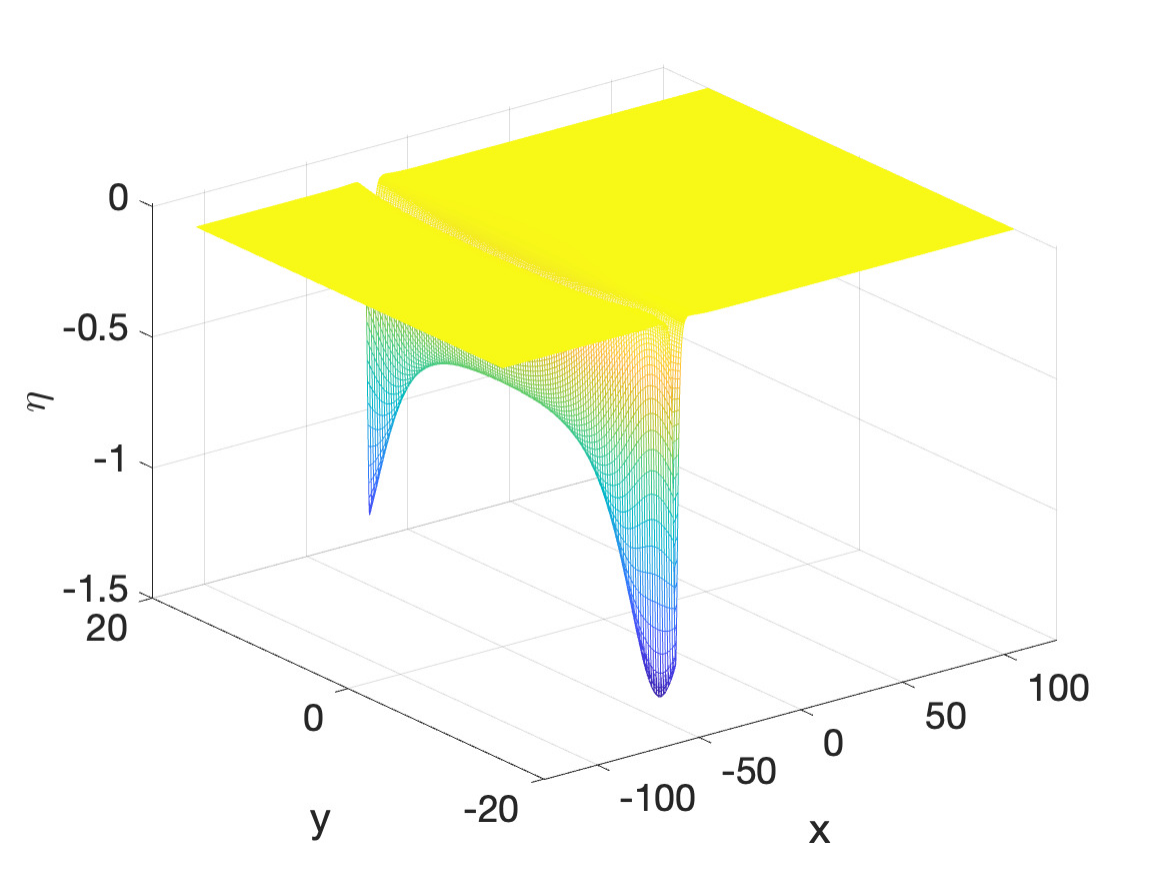}
    \end{subfigure}
    \hfill
    \begin{subfigure}[b]{0.49\textwidth}
        \centering
        \includegraphics[width=\textwidth]{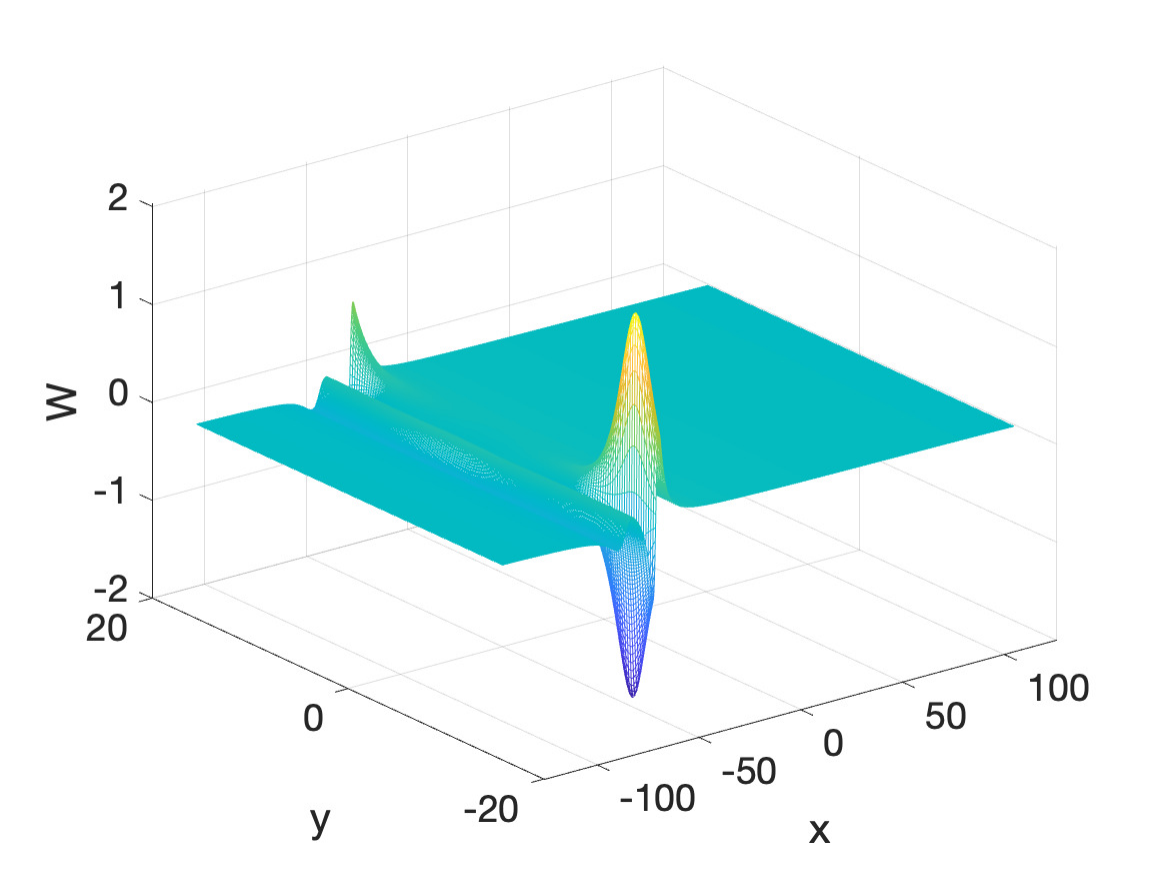}
    \end{subfigure}

\caption{Solution to the KBK system 
(\ref{Kaup-true-2D}) with initial data (\ref{initrans}) for $\mu = 
0.01$ and $C=-0.6$. Rows correspond to times $t=85$ (top), $t=120$ (middle), and $t=15$ (bottom). In each row, $\eta$ is shown on the left and $W$ on the right.}
    \label{fig:combined_06}
\end{figure}

\begin{figure}[htb!]
\includegraphics[width=0.49\textwidth]{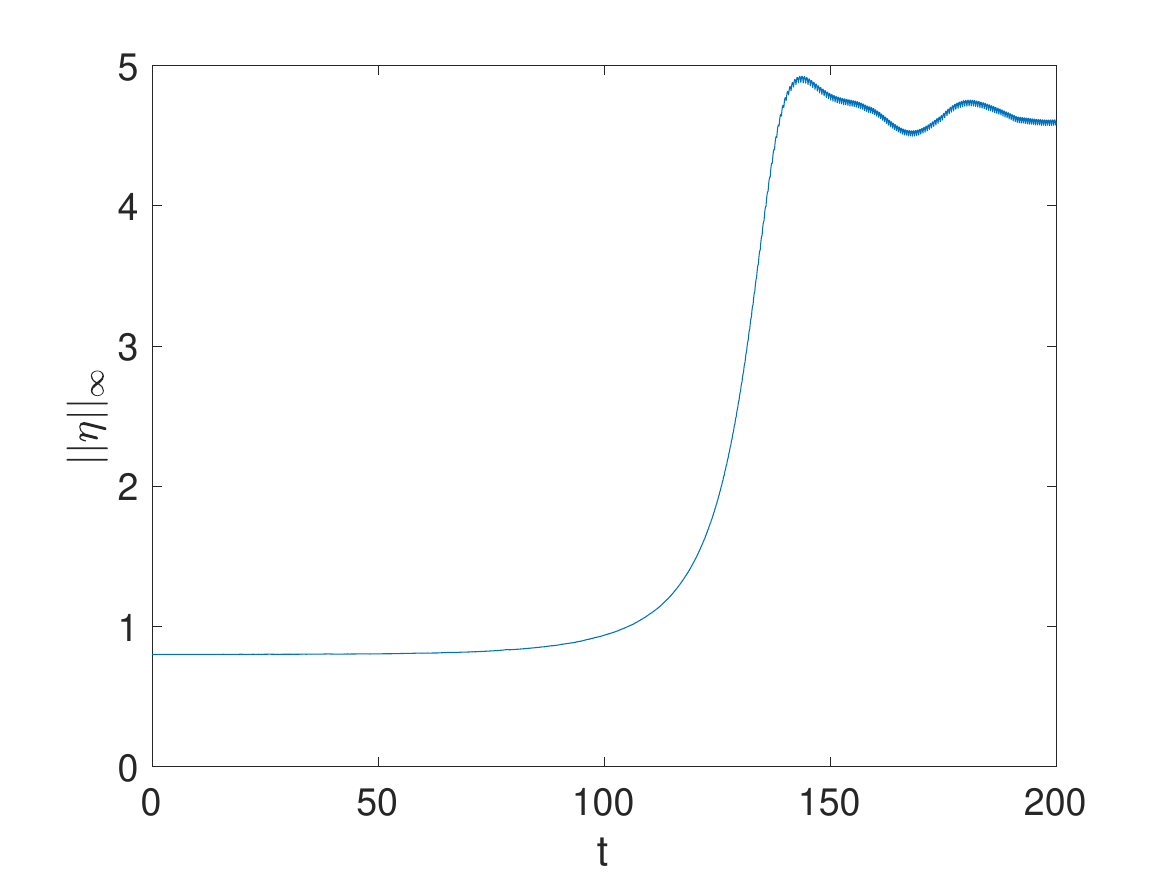}
\includegraphics[width=0.49\textwidth]{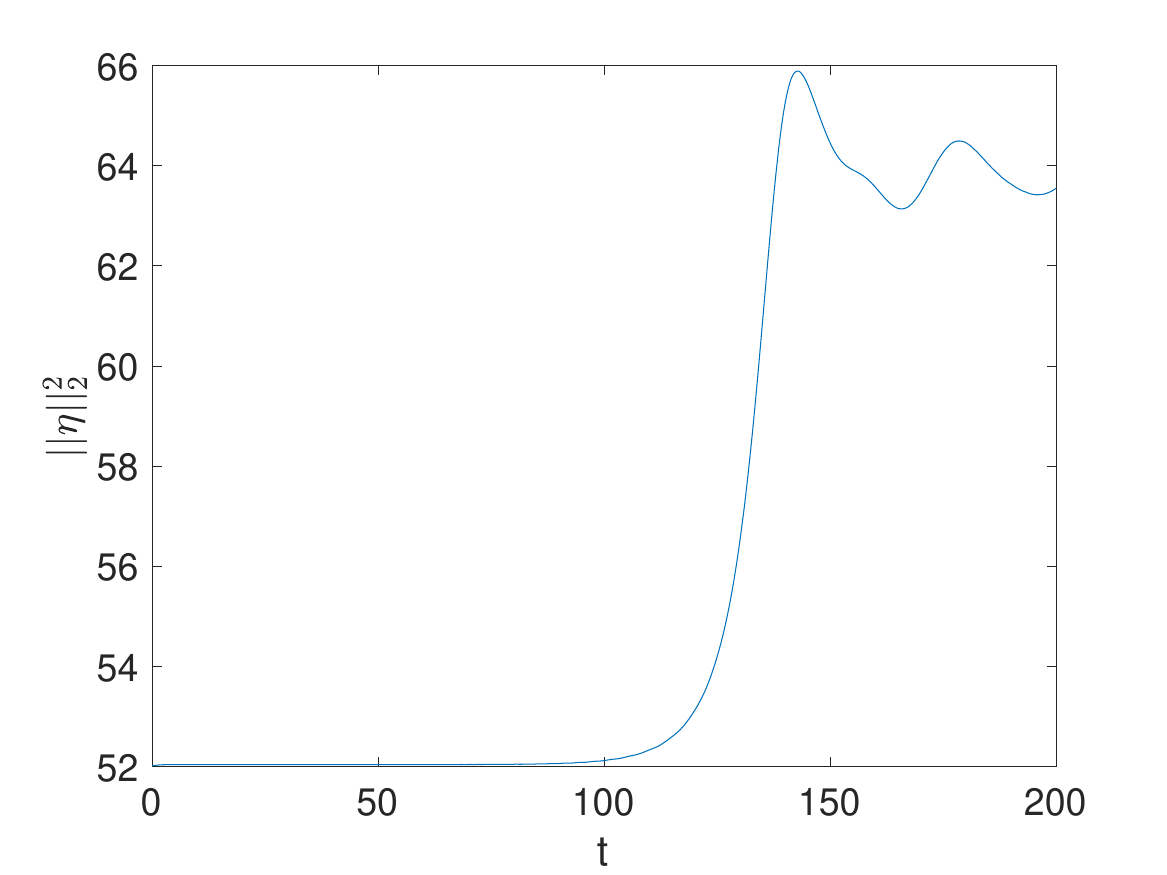}
\caption{Norms of the solution to the system KBK (\ref{Kaup-true-2D}) for the 
initial data (\ref{initrans}) with $\mu=0.01$ and $C=-0.6$ in dependence of 
time, on the 
left $||\eta||_{\infty}$, on the right $||\eta||_{2}^{2}$.} 
\label{fig:normsolc06}
\end{figure}

We now turn to $C < -0.5$, for which the profile (2) 
satisfies $1 + \eta > 0$. For $C = -0.6$, $\mu = 0.01$ and $t \leq 
200$ (we work with $L_x = 40$, $L_y = 5$, $N_x = 2^{11}$, and 
$N_y = 2^{7}$ and $N_{t}=10^{4}$), the perturbed line solitary wave first becomes 
slowly narrower and deeper and then evolves into a localised 
structure, which we identify as a lump (see 
Fig.~\ref{fig:combined_06}). After a rapid growth around $t \approx 
130$, the $L^\infty$ norm of $\eta$ stays essentially constant, at 
about $4.7$, up to small oscillations (see Fig.~\ref{fig:normsolc06}) 
around a presumably stable final state. Note that this lump does not 
satisfy the condition $1+\eta>0$ though the perturbed line solitary 
wave did. However, there is no indication of a blow-up.

The situation is very similar for larger or negative 
values of $\mu$. For $C=-0.7$, we get once more a lump for these 
perturbations which is thereforee not shown.

%
%

\begin{figure}[htb!]
    \centering

    \begin{subfigure}[b]{0.49\textwidth}
        \centering
        \includegraphics[width=\textwidth]{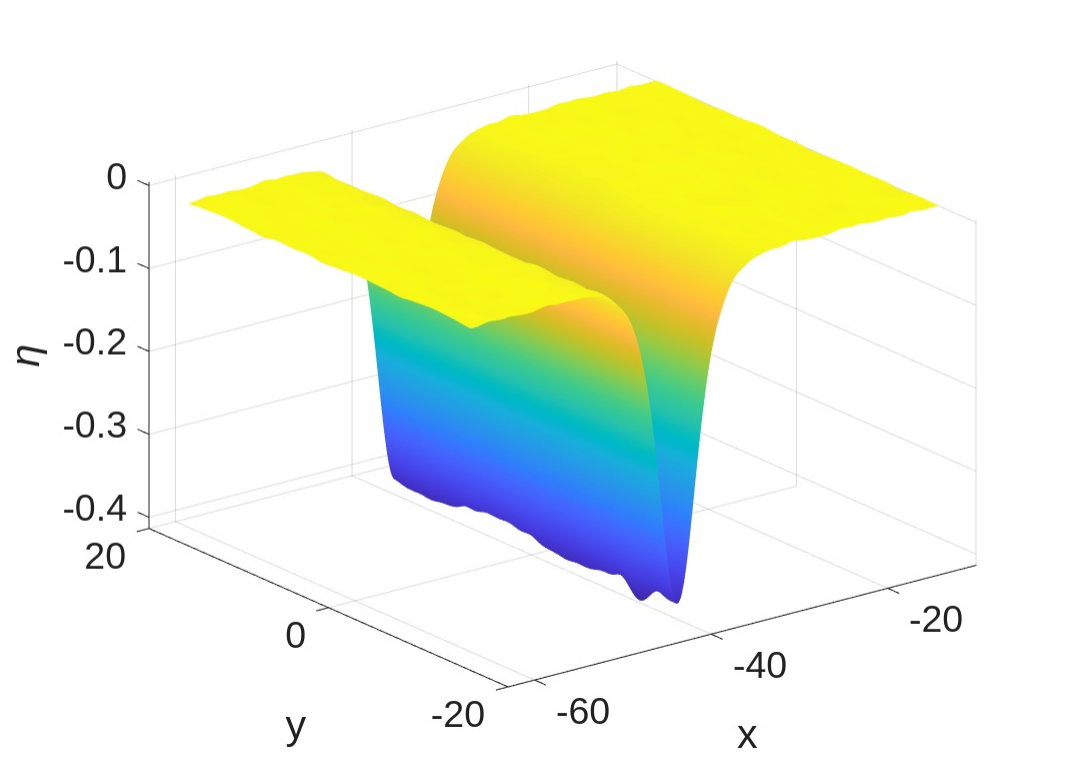}
    \end{subfigure}
    \hfill
    \begin{subfigure}[b]{0.49\textwidth}
        \centering
        \includegraphics[width=\textwidth]{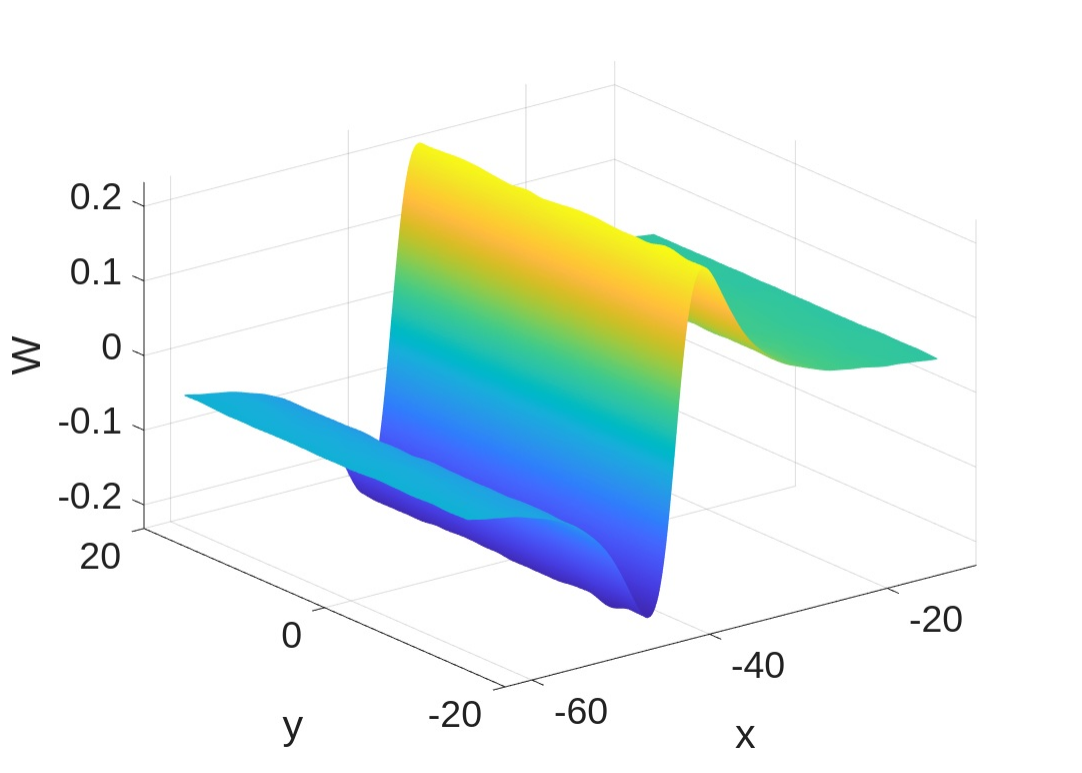}
    \end{subfigure}
\caption{Solution of the KBK system (\ref{Kaup-true-2D}) with initial data (\ref{initrans}), for $\mu = +0.1$ and $C = -0.8$ for $t=200$, on the left $\eta$, on the right $W$.}
    \label{fig:trans08}
\end{figure}

\begin{figure}[htb!]
\includegraphics[width=0.49\textwidth]{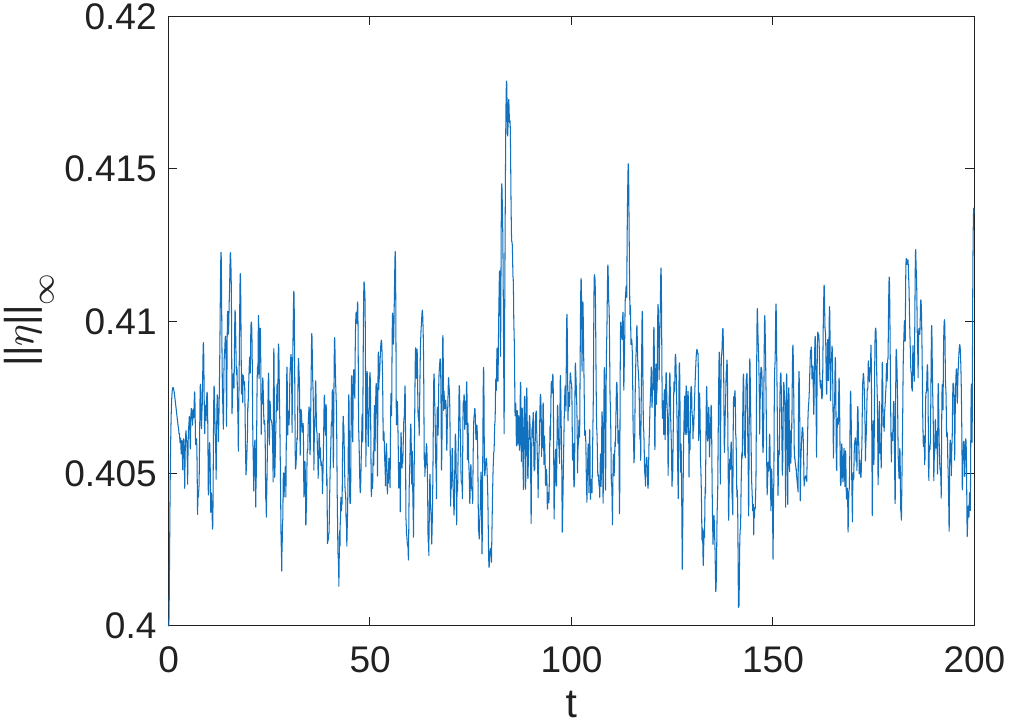}
\includegraphics[width=0.49\textwidth]{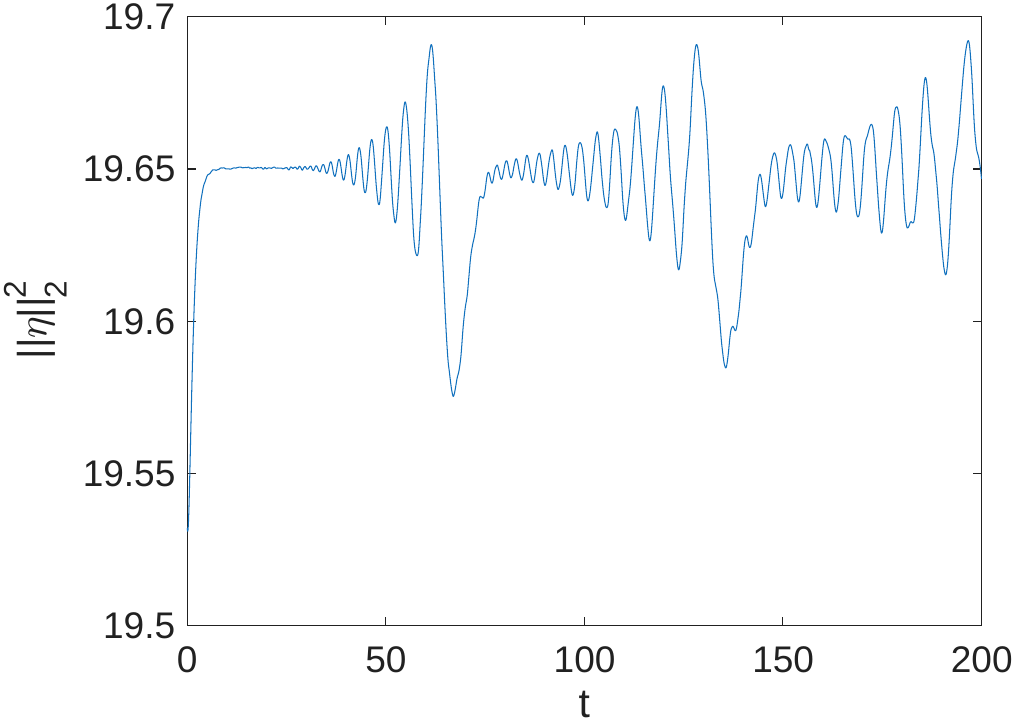}
\caption{Norms of the solution to the system KBK (\ref{Kaup-true-2D}) for the 
initial data (\ref{initrans}) with $\mu=+0.1$ and $C=-0.8$ in dependence of 
time, on the 
left $||\eta||_{\infty}$, on the right $||\eta||_{2}^{2}$.} 
\label{fig:normsolc08}
\end{figure}

For $C = -0.8$ and $\mu = +0.1$ the line solitary 
wave  is stable and only weak radiation appears (see 
Fig.~\ref{fig:trans08}). The $L^\infty$ norm stays close to the 
soliton amplitude $2(1 + C) = 0.4$ (see Fig.~\ref{fig:normsolc08}). 
The same holds for $C = -0.9$ and $\mu = \pm 0.1$. These simulations 
suggest that lumps form for $C$ just below the non-cavitation 
threshold $C = -0.5$, whereas line solitary waves are stable for $C$ closer to $-1$.

\section{Localised initial data}
In this section we study the KBK system (\ref{Kaup-true-2DF}) for 
localised initial data of the form 
\begin{equation}
	\eta(x,y,0)=\kappa_{1}\exp(-x^{2}-y^{2}),\quad 
	V(x,y,0)=\kappa_{2}\exp(-x^{2}-y^{2}),
	\label{inilocalised}
\end{equation}
where $\kappa_{1,2}$ are constants. Note that the static solution of 
section \ref{secstatic} is also localised, but of a very specific 
form. It was shown in section \ref{secpertstatic} that perturbations 
of this solution lead either to pure scattering or to a blow-up in 
finite time. Thus it was already shown that localised initial data 
can lead to a blow-up in solutions to the KBK system in 2D. 
In this section we consider Gaussian initial data with 
values of $|\kappa_{1,2}|<10$ for which we only observe scattering. 
Once more we use  $L_{x}=L_{y}=5$, 
$N_{x}=N_{y}=2^{9}$ with $N_{t}=10^{3}$ time 
steps for $t\leq1$.

\begin{figure}[htb!]
\includegraphics[width=0.49\textwidth]{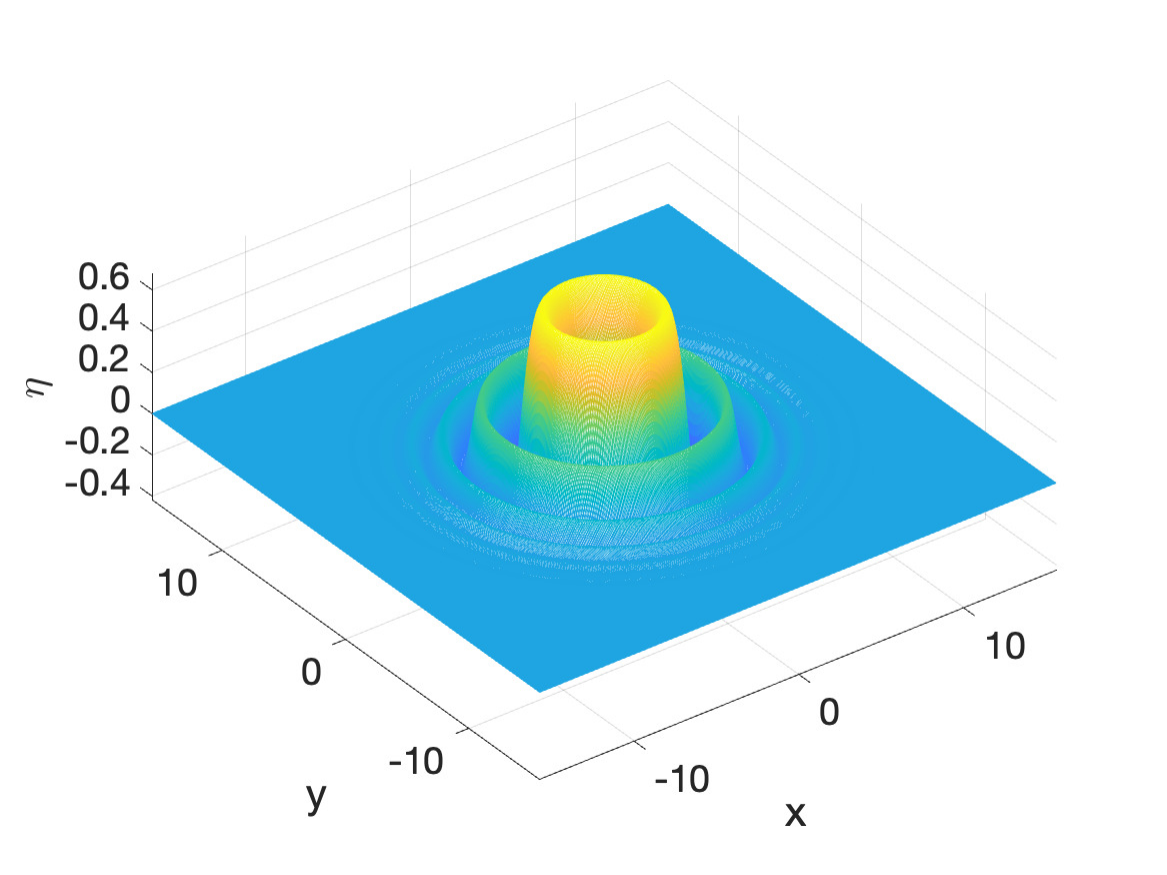}
\includegraphics[width=0.49\textwidth]{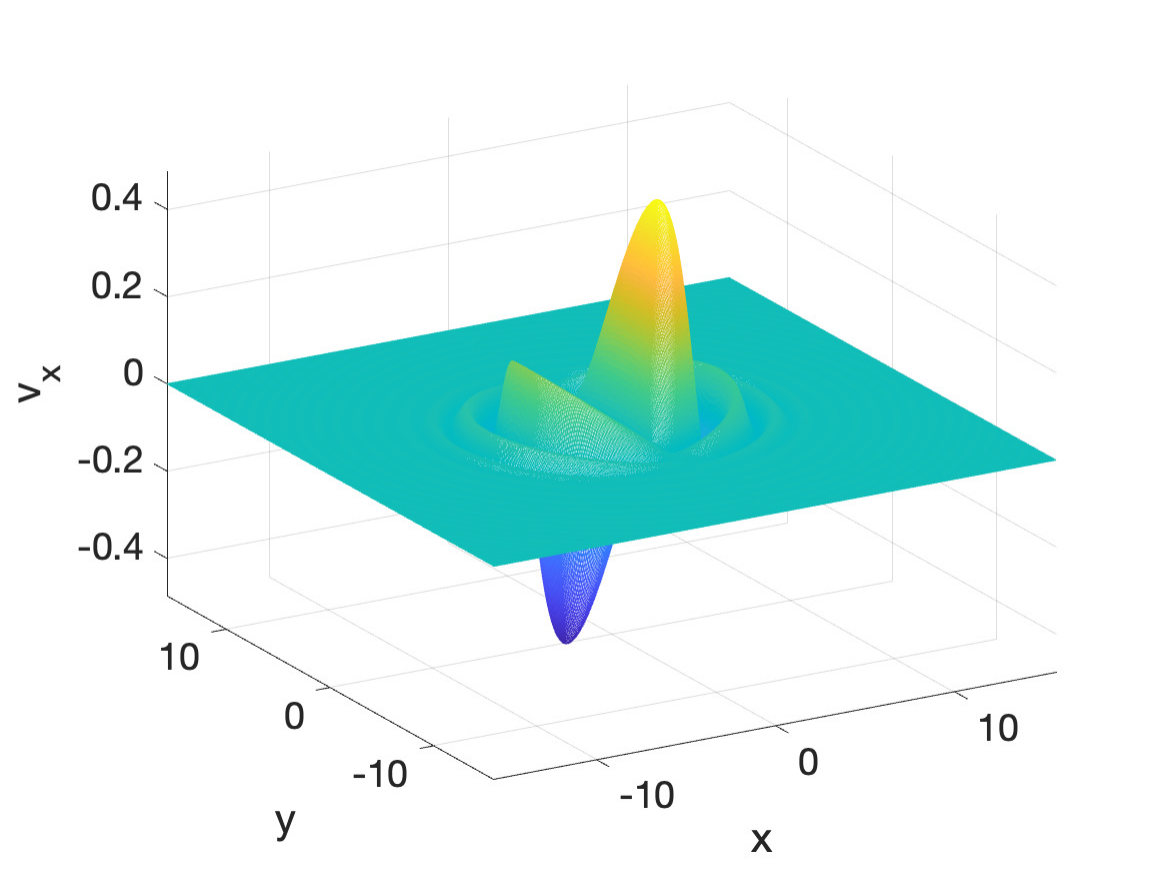}
\caption{Solution to the system KBK (\ref{Kaup-true-2D}) for the 
initial data (\ref{inilocalised}) with $\kappa_{1}=5$ and 
$\kappa_{2}=0$ for $t=1$,  on the 
left $\eta$, on the right $v_{x}$. }
\label{5gausst1}
\end{figure}

First we study the case $\kappa_{1}=5$ and $\kappa_{2}=0$. The 
solution at the final time $t=1$ is shown in Fig.~\ref{5gausst1}. 
It appears that the initial data form an annular structure that is 
simply dispersed to infinity.

\begin{figure}[htb!]
\includegraphics[width=0.49\textwidth]{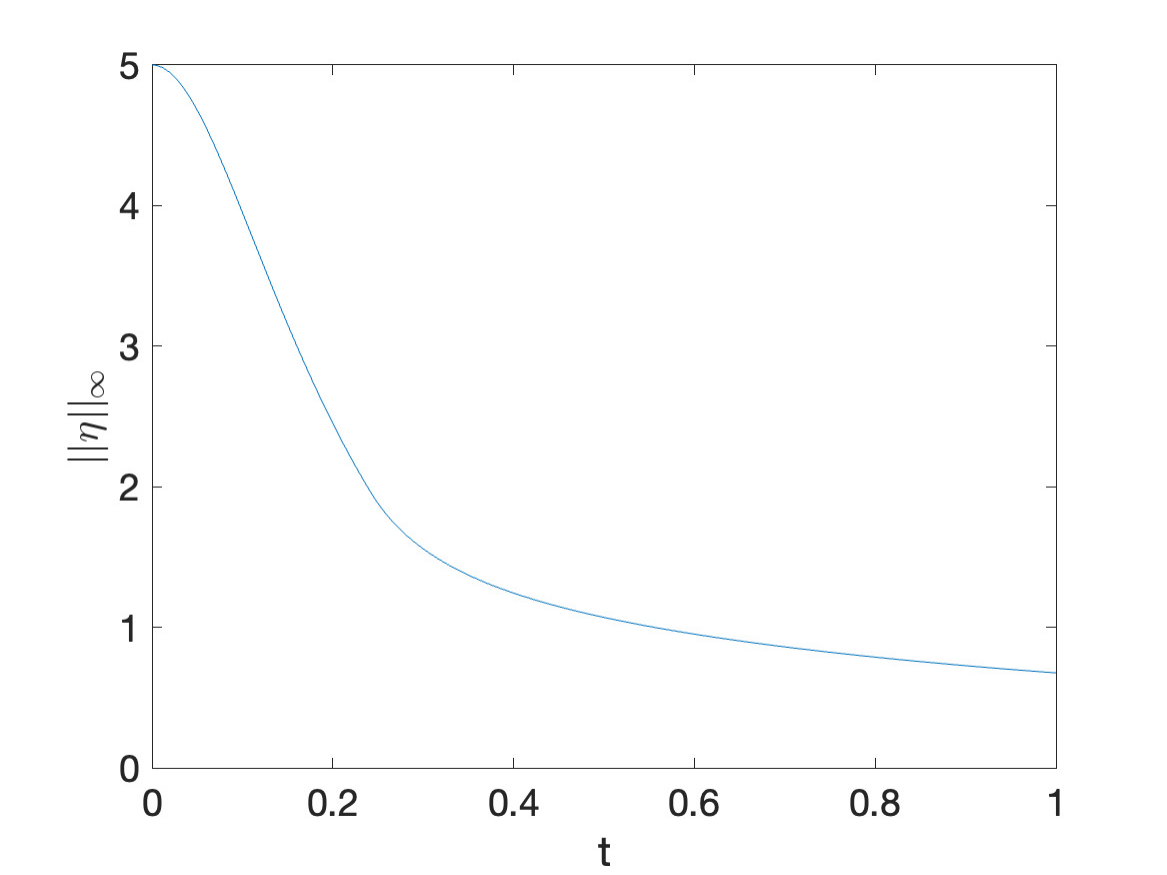}
\includegraphics[width=0.49\textwidth]{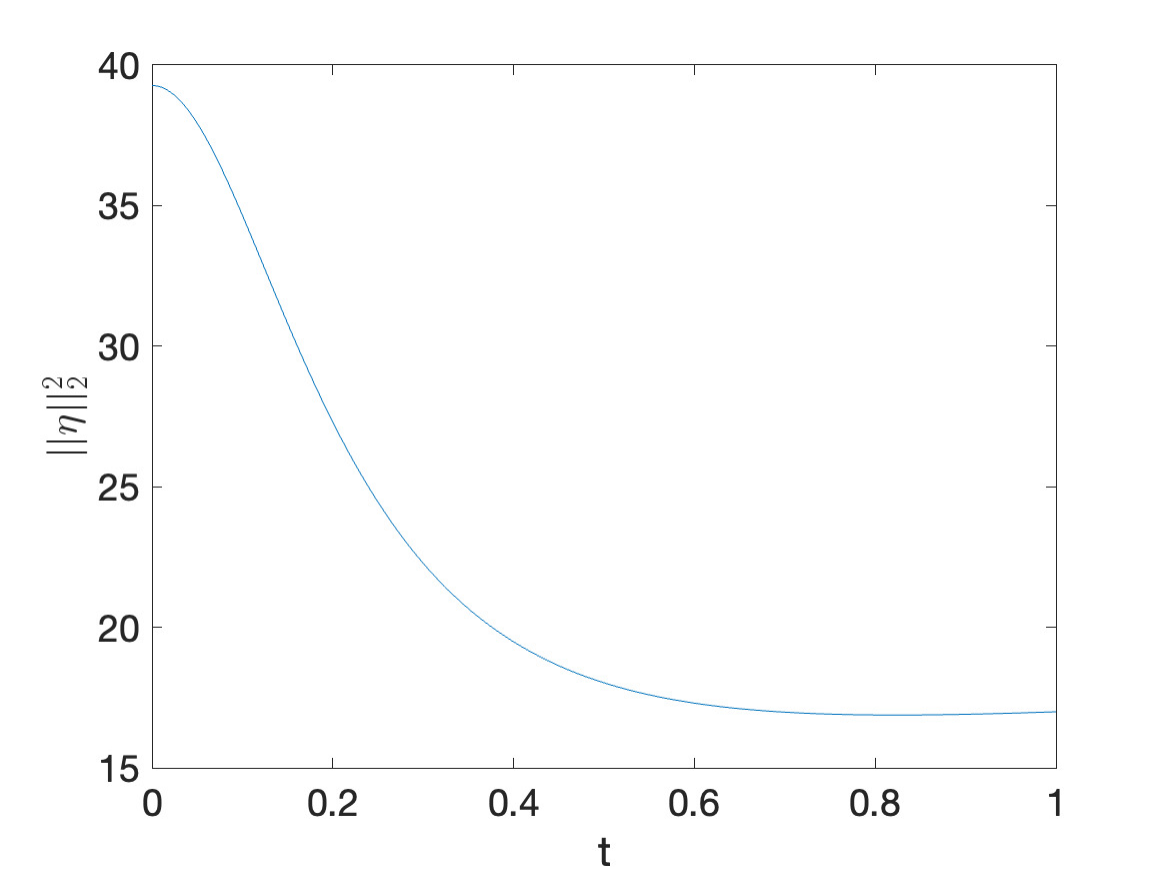}
\caption{Norms of the solution to the system KBK (\ref{Kaup-true-2D}) for the 
solution shown in Fig.~\ref{5gausst1} in dependence of 
time, on the 
left $||\eta||_{\infty}$, on the right $||\eta||_{2}^{2}$. }
\label{5gaussinf}
\end{figure}

This is confirmed by the $L^{\infty}$ and the $L^{2}$ norm of $\eta$ 
shown in Fig.~\ref{5gaussinf}. Both appear to be monotonically 
decreasing. Since we work on a torus, the radiation cannot escape to 
infinity, but there is no indication of a stable localised state 
appearing in the time evolution of these initial data. 

\begin{figure}[htb!]
\includegraphics[width=0.49\textwidth]{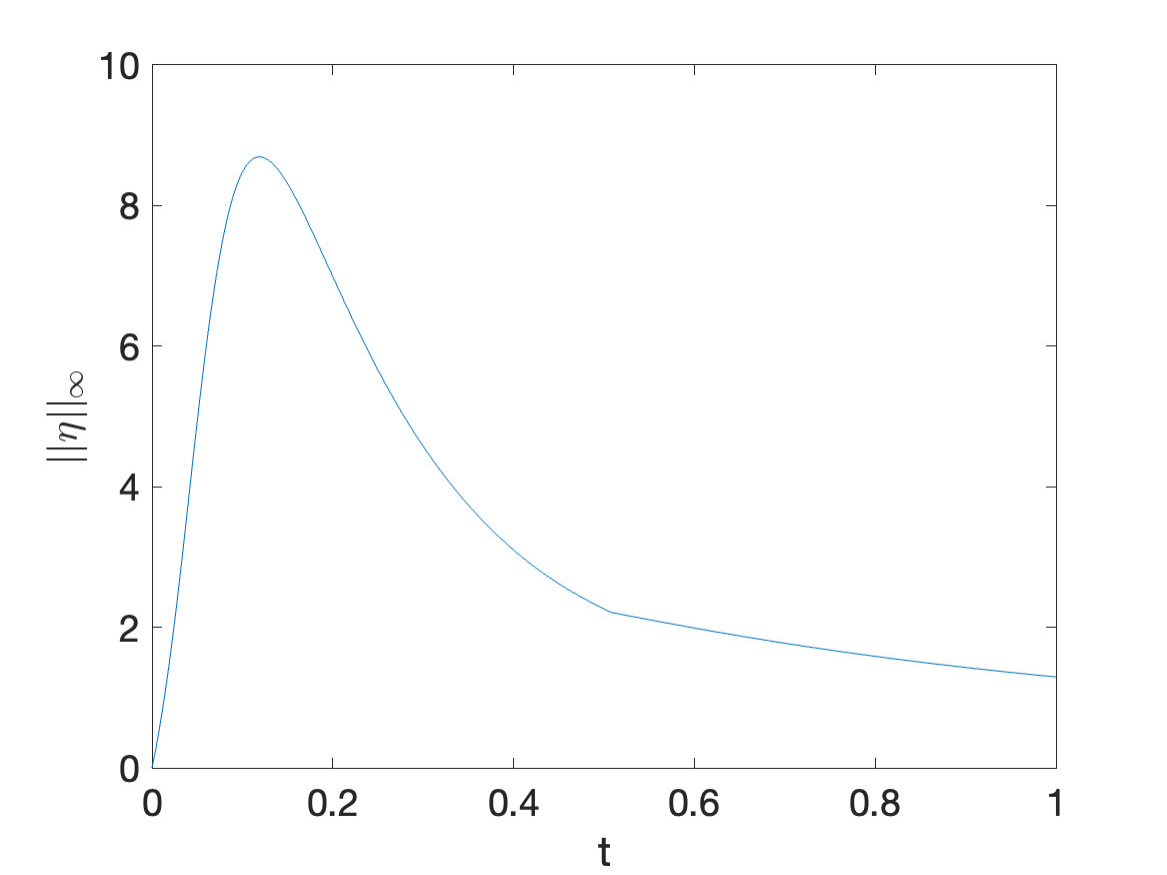}
\includegraphics[width=0.49\textwidth]{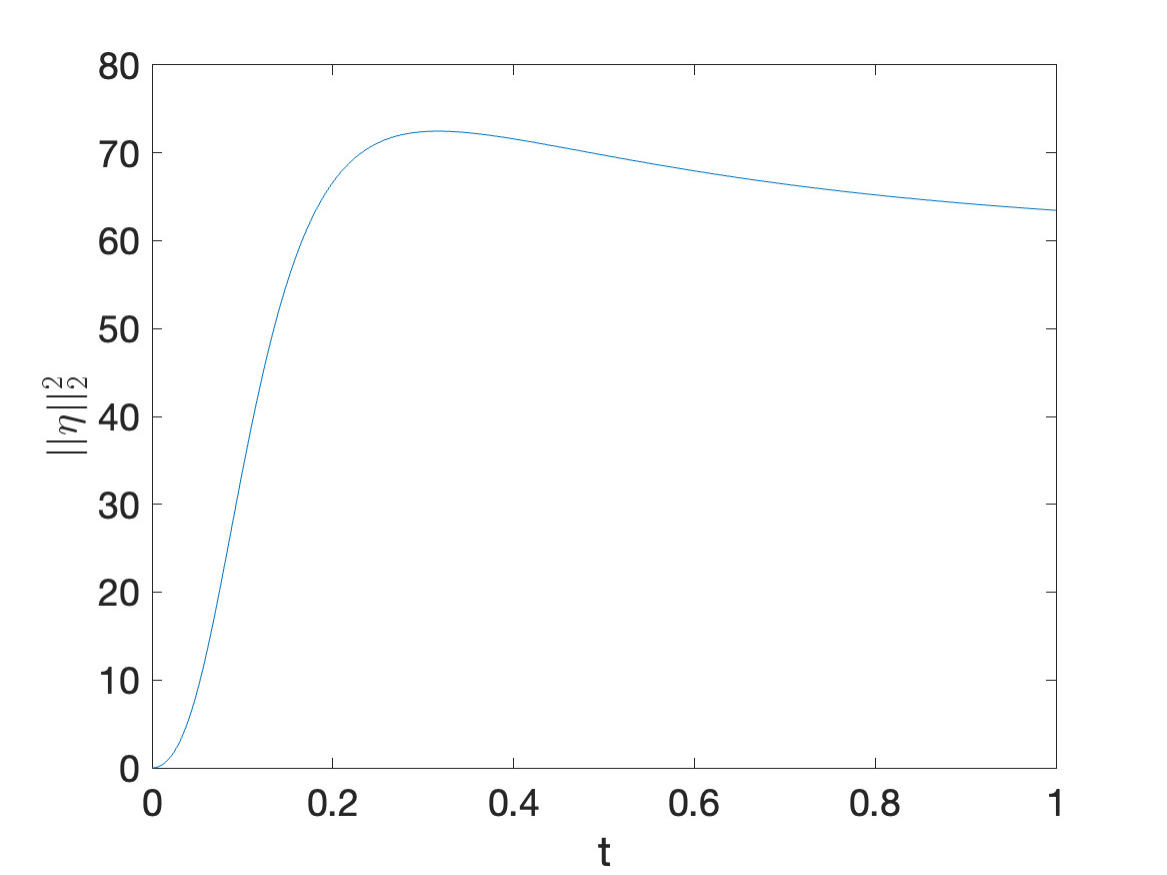}
\caption{Norms of the solution to the system KBK (\ref{Kaup-true-2D}) for the 
initial data (\ref{inilocalised}) with $\kappa_{1}=0$, $\kappa_{2}=5$ in dependence of 
time, on the 
left $||\eta||_{\infty}$, on the right $||\eta||_{2}^{2}$. }
\label{5gausswinf}
\end{figure}

The situation is somewhat different in the case $\kappa_{1}=0$ and 
$\kappa_{2}=5$ of the initial data (\ref{inilocalised}). As can be 
seen in Fig.~\ref{5gausswinf}, the solution is first growing with 
time, but will eventually disperse. 

\begin{figure}[htb!]
\includegraphics[width=0.49\textwidth]{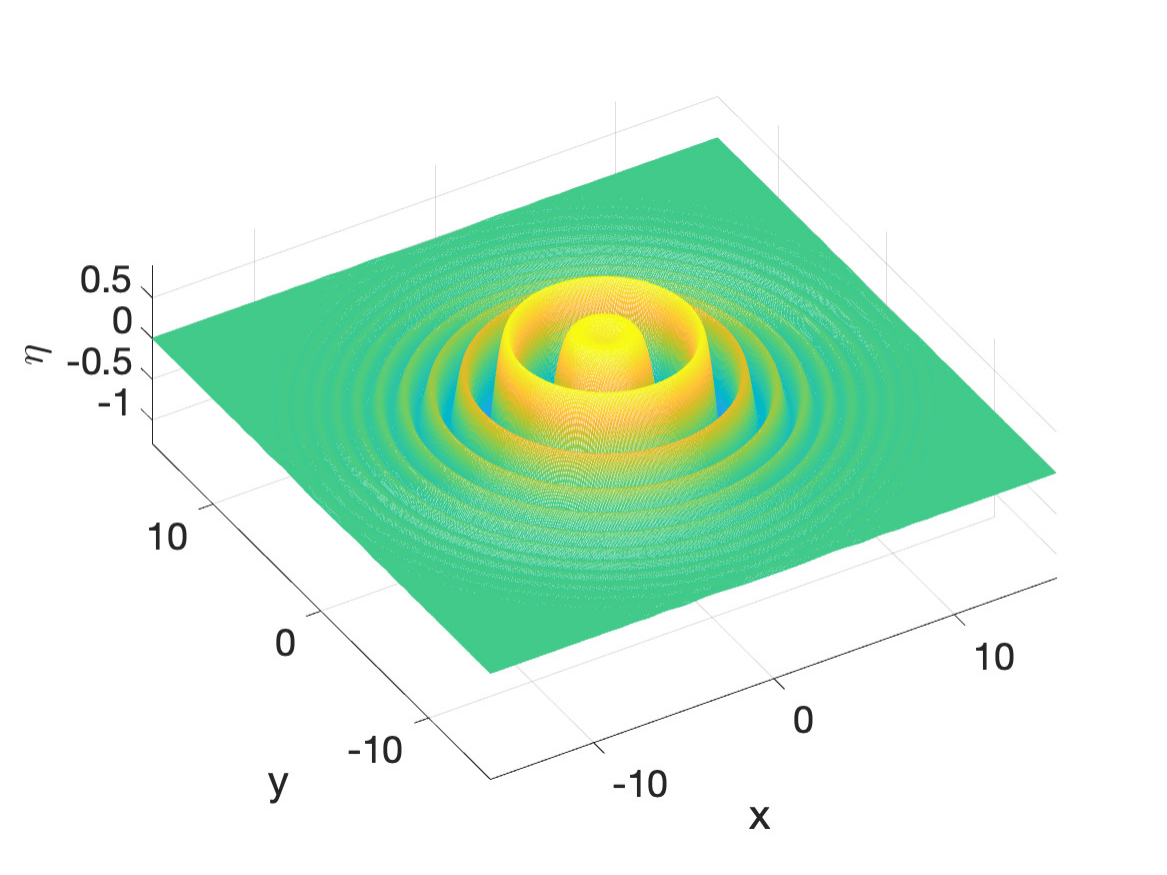}
\includegraphics[width=0.49\textwidth]{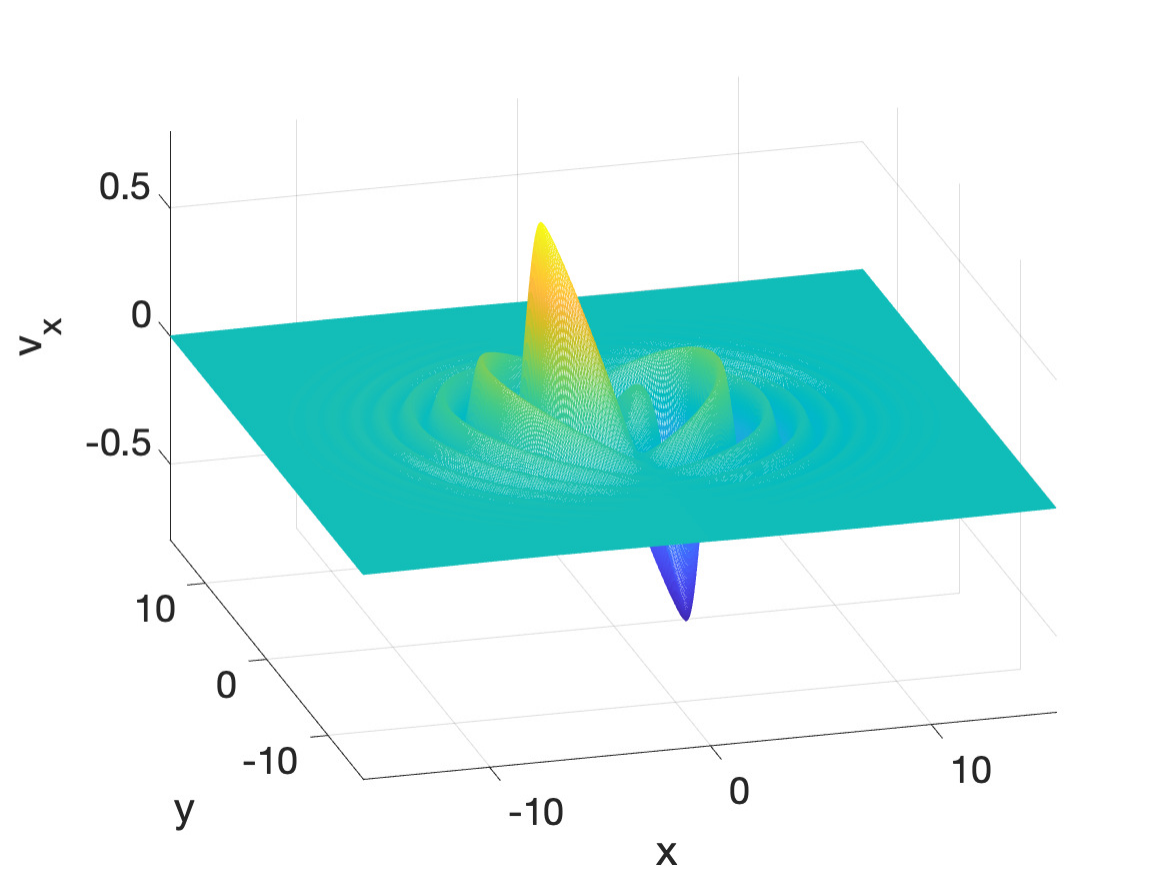}
\caption{Solution to the system KBK (\ref{Kaup-true-2D}) for the 
initial data (\ref{inilocalised}) with $\kappa_{1}=0$ and 
$\kappa_{2}=5$ for $t=1$,  on the 
left $\eta$, on the right $v_{x}$. }
\label{5gausswt1}
\end{figure}

The solution in this case at the final time is shown in 
Fig.~\ref{5gausswt1}. Thus once more no stable structure is observed 
in the time evolution. Larger values of $\kappa_{2}$ for 
$\kappa_{1}=0$ lead to stronger growth of the $L^{\infty}$ norm, but 
a qualitatively similar behavior. 

\section{Conclusion}
In this paper we have presented a  numerical study of solutions to 
the irrotational 2D KBK system. It was shown that a static solution 
localised in two spatial dimensions can be numerically constructed. This solution is unstable 
against both dispersion (for perturbations with a smaller $L^{2}$ 
norm of $\eta$ than 
the static solution) and blow-up in finite time (for perturbations 
with larger $||\eta||_{2}$). The line solitary waves, $y$-independent and thus 
infinitely extended 1D solitons, were shown to be unstable against 
blow-up and thus strongly unstable in the case where the wave profile $\eta$ initially violates the non-cavitation condition, i.e. when $C>-0.5$. This blow-up could be identified to correspond to the 
blow-up observed for $L^{2}$ supercritical NLS equations. 

The results can be summarized as follows:
\begin{enumerate}
	\item Solutions to the irrotational 2D KBK systems for localised smooth initial 
	data  can blow up at finite time $t^{*}$ as solutions to the 
	supercritical NLS equation. As per equations  
	(\ref{dyn}) and (\ref{KBKresc}) the blow up is self similar. 
	\begin{equation}
		\eta\propto \frac{U^{\infty}(X,Y)}{L^{2}}, \quad V\propto
		V^{\infty}(X,Y),\quad L\propto (t^{*}-t)^{1/2}
		\label{bu}.
	\end{equation}
	The precise criteria for initial data to blow up as well as the blow-up 
	profiles $U^{\infty}$ and $V^{\infty}$ are unknown. 
	\item The line solitary waves are strongly unstable 
	when they initially violate the non-cavitation 
	condition $\eta +1>0$, i.e. when $C>-0.5$. In these cases perturbations 
	of the line solitary wave blow up in a self similar way
	according to (\ref{bu}) even for a mass on the torus smaller than 
	the mass of the static solution.
	\item For $-1<C^{*}<C<-0.5$, the line solitary waves are unstable 
	agains the formation of lumps, solitary waves localised in two 
	dimensions. These lumps appear to be stable (only lumps violating 
	the condition $1+\eta>0$ were found). The precise value of the 
	constant $C^{*}$ is unknown. 
	\item For $-1<C<C^{*}$, the line solitary wave is stable. 
	
\end{enumerate}

\section*{Competing interests declaration}
All authors that have contributed to the  submission declare that 
they have no competing interests. 

\bibliographystyle{amsplain}

\end{document}